\documentclass[a4paper]{amsproc}

\usepackage{amssymb}
\usepackage{amscd}
\usepackage{amsxtra}
\usepackage[T1]{fontenc}
\usepackage[latin1]{inputenc}
\usepackage{enumerate}

\usepackage[pdftitle={Primes, L-functions and group representations},
  pdfauthor={Ralf Meyer},
  pdfsubject={Mathematics; MSC 11M26, 22D12, 18H10, 43A35, 58B34},
]{hyperref}
\usepackage[lite]{amsrefs}

\newcommand*{\brd}{-\hspace{0pt}}
\newcommand*{\nbd}{\nobreakdash-\hspace{0pt}}

\newcommand*{\N}{{\mathbb{N}}}
\newcommand*{\Z}{{\mathbb{Z}}}
\newcommand*{\Q}{{\mathbb{Q}}}
\newcommand*{\R}{{\mathbb{R}}}
\newcommand*{\C}{{\mathbb{C}}}

\newcommand*{\abs}[1]{\lvert#1\rvert}
\newcommand*{\norm}[1]{\lVert#1\rVert}
\newcommand*{\ket}[1]{\lvert#1\rangle}
\newcommand*{\bra}[1]{\langle#1\vert}
\newcommand*{\braket}[2]{\langle#1,#2\rangle}
\newcommand*{\ooival}[1]{\mathopen]#1\mathclose[}
\newcommand*{\ocival}[1]{\mathopen]#1\mathclose]}
\newcommand*{\coival}[1]{\mathopen[#1\mathclose[}
\newcommand*{\conj}[1]{\overline{#1}}
\newcommand*{\Cont}[1]{\complement #1}
\newcommand*{\inv}{^\times}
\newcommand*{\one}{^{1}}
\newcommand*{\blank}{{\llcorner\!\!\lrcorner}}

\newcommand*{\LF}{{\mathcal{K}}}
\newcommand*{\GF}{{K}}

\newcommand*{\MCR}{{\mathcal{O}}}
\newcommand*{\MI}{{\mathcal{P}}}

\newcommand*{\Adel}{{\mathbb{A}}}
\newcommand*{\ICL}{{\mathcal{C}}}

\newcommand*{\Irrep}{\Omega}
\newcommand*{\Places}{\mathcal{P}}

\newcommand*{\Mod}{\mathrm{Mod}}
\newcommand*{\Der}{\mathrm{Der}}
\DeclareMathOperator{\Right}{\mathbb{R}}

\newcommand*{\ima}{\mathrm{i}}
\newcommand*{\ID}{\mathrm{id}}

\newcommand*{\comp}{{\mathrm{c}}}
\newcommand*{\IN}{{\smallint\!}}
\newcommand*{\pp}{>}
\newcommand*{\mm}{<}
\newcommand*{\ppmm}{{><}}

\newcommand*{\different}{\partial}
\newcommand*{\Discriminant}{D}

\newcommand*{\Cred}{C^\ast_{\mathrm{red}}}
\newcommand*{\Mult}{\mathcal{M}}
\newcommand*{\Sch}{\mathcal{S}}
\newcommand*{\Twist}{\mathcal{T}}
\newcommand*{\CCINF}{\mathcal{D}}
\newcommand*{\CINF}{\mathcal{E}}

\newcommand*{\Lop}{\mathcal{L}}
\newcommand*{\Lopi}{\mathcal{L}^{-1}}

\newcommand*{\SUM}{\Sigma}
\newcommand*{\Sing}{R}
\newcommand*{\Hilm}{\mathcal{H}}
\newcommand*{\sumpm}{\mathcal{X}}
\newcommand*{\Fourier}{\mathfrak{F}}
\newcommand*{\EFou}{\mathsf{F}}
\newcommand*{\EFoi}{\mathsf{F}^\ast}

\newcommand*{\hot}{\mathbin{\hat{\otimes}}}
\newcommand*{\Lhot}{\mathbin{\mathbb{L}{\hat{\otimes}}}}

\newcommand*{\defeq}{\mathrel{:=}}
\newcommand*{\congto}{\mathrel{\overset{\cong}{\to}}}
\newcommand*{\boin}{\prec}
\newcommand*{\prto}{\twoheadrightarrow}
\newcommand*{\into}{\rightarrowtail}
\newcommand*{\addconv}{\dagger}
\newcommand*{\cross}{\ltimes}

\DeclareMathOperator{\supp}{supp}

\DeclareMathOperator{\tr}{tr}
\DeclareMathOperator{\vol}{vol}
\DeclareMathOperator{\RE}{Re}

\DeclareMathOperator{\Ker}{Ker}

\DeclareMathOperator{\Hom}{Hom}
\DeclareMathOperator{\End}{End}
\DeclareMathOperator{\Aut}{Aut}
\DeclareMathOperator{\smooth}{smooth}
\DeclareMathOperator{\mult}{mult}
\DeclareMathOperator{\ord}{ord}
\DeclareMathOperator{\spec}{spec}

\theoremstyle{plain}
\newtheorem{theorem}{Theorem}[section]
\newtheorem{proposition}[theorem]{Proposition}
\newtheorem{lemma}[theorem]{Lemma}
\newtheorem{corollary}[theorem]{Corollary}

\theoremstyle{definition}
\newtheorem{definition}[theorem]{Definition}

\begin{document}

\title[Primes, \(L\)-functions and group representations]{On a
  representation of the idele class group related to primes and zeros
  of \(L\)\nobreakdash-functions}

\author{Ralf Meyer}

\address{Mathematisches Institut\\
         Westfälische Wilhelms-Universität Münster\\
         Einsteinstraße~62\\
         48149 Münster\\
         Germany}

\email{rameyer@math.uni-muenster.de}

\begin{abstract}
  Let~\(\GF\) be a global field.  Using natural spaces of functions on
  the adele ring and the idele class group of~\(\GF\), we construct a
  virtual representation of the idele class group of~\(\GF\) whose
  character is equal to a variant of the Weil distribution that occurs
  in André Weil's Explicit Formula.  Hence this representation encodes
  information about the distribution of the prime ideals of~\(\GF\) and
  is a spectral interpretation for the poles and zeros of the
  \(L\)\nobreakdash-function of~\(\GF\).  Our construction is motivated by
  a similar spectral interpretation by Alain Connes.
\end{abstract}


\thanks{This research was supported by the EU-Network \emph{Quantum
  Spaces and Noncommutative Geometry} (Contract HPRN-CT-2002-00280)
  and the \emph{Deutsche Forschungsgemeinschaft} (SFB 478).}

\maketitle

\tableofcontents

\section{Introduction}
\label{sec:intro}

Let~\(\GF\) be an algebraic number field, that is, a finite algebraic
extension of the field of rational numbers~\(\Q\).  Already the
case~\(\Q\) is of great interest.  Let \(\Places(\GF)\) be the set of
places of~\(\GF\), that is, equivalence classes of dense embeddings
of~\(\GF\) into local fields.  For instance, \(\Places(\Q)\) contains the
embeddings of~\(\Q\) into~\(\R\) and into the \(p\)\nbd{}adic
integers~\(\Q_p\) for all prime numbers~\(p\).  Let~\(\Adel_\GF\) be the
adele ring of~\(\GF\), let~\(\Adel\inv_\GF\) be the idele group and let
\(\ICL_\GF\defeq \Adel\inv_\GF/\GF\inv\) be the idele class group
of~\(\GF\) (see~\cite{Weil:Basic} for these constructions).  For
instance, the group~\(\ICL_\Q\) is isomorphic to the direct product
\(\prod_p \Z_p\inv \times \R\), where~\(p\) runs through the prime numbers
and~\(\Z_p\inv\) denotes the multiplicative group of \(p\)\nbd{}adic
integers of norm~\(1\).

Let \(\Irrep(\ICL_\GF)\defeq\Hom(\ICL_\GF,\C\inv)\) be the space of
quasi-characters of~\(\ICL_\GF\).  We replace the family of
\(L\)\nbd{}functions \(L(\chi,s)\) of~\(\GF\) associated to
characters~\(\chi\) of~\(\ICL_\GF\) by a single meromorphic function
\(L_\GF\colon \Irrep(\ICL_\GF)\to\C\) which satisfies
\(L_\GF(\chi\abs{x}^s)= L(\chi,s)\) for all \(\chi\in\widehat{\ICL_\GF}\),
\(s\in\C\).  We only work with complete \(L\)\nbd{}functions, that is, we
multiply by appropriate \(\Gamma\)\nbd{}functions for the infinite
places.  We let \(\ord(\omega,L_\GF)\) be the pole order of~\(L_\GF\)
at~\(\omega\).  This is positive for poles and negative for zeros
of~\(L_\GF\).

A \emph{virtual representation} is a pair of representations~\(\pi_\pm\)
which we interpret as a formal difference \(\pi_+\ominus\pi_-\).  We
construct a spectral interpretation for the poles and zeros
of~\(L_\GF\), that is, a virtual representation~\(\pi\) of~\(\ICL_\GF\)
whose spectrum consists precisely of the poles and zeros of the
\(L\)\nbd{}function~\(L_\GF\).  The two poles occur in~\(\pi_+\), the zeros
in~\(\pi_-\).  The spectral multiplicity \(\mult(\omega,\pi)\) is
equal to the pole order \(\ord(\omega,L_\GF)\).  Therefore, the
character of~\(\pi\) is the distribution~\(\chi_\pi\) on~\(\ICL_\GF\)
defined by
\[
\chi_\pi(f) \defeq
\sum_{\omega\in\Irrep(\ICL_\GF)} \mult(\omega,\pi) \cdot
\hat{f}(\omega)
=
\sum_{\omega\in\Irrep(\ICL_\GF)} \ord(\omega,L_\GF) \cdot
\hat{f}(\omega)
\]
for \(f\in\CCINF(\ICL_\GF)\), where~\(\hat{f}\) denotes the
Fourier-Laplace transform of~\(f\).  In particular, the
representation~\(\pi\) contains a spectral interpretation for the
non-trivial zeros of the \(\zeta\)\nbd{}function of~\(\GF\), which is
Riemann's well-known function for \(\GF=\Q\).

Our construction is motivated by Alain Connes' spectral interpretation
for the zeros of \(L\)\nbd{}functions (\cite{Connes:Trace_Formula}).
There zeros on the critical line appear as above, but zeros off the
critical line appear as resonances.  Given the Explicit Formula of
André Weil, a geometric computation of the ``character'' in Connes'
situation is therefore equivalent to the Generalized Riemann
Hypothesis.  Since I have not been able to do this computation, I have
developed an alternative spectral interpretation that is no longer
directly related to the Riemann Hypothesis.  First we consider the
\(S\)\nbd{}local case, where we can work in the same Hilbert space
situation considered by Connes.  Our \(S\)\nbd{}local computations are,
in fact, essentially equivalent to Connes'.  In the global case, we
deviate from~\cite{Connes:Trace_Formula} and use smaller function
spaces instead of the Hilbert spaces in Connes' setup.  These are no
longer biased in favor of the critical line, so that all zeros and
poles appear in the same way.

The trace of a matrix can be computed spectrally as the sum of its
eigenvalues and geometrically as the sum of its diagonal entries.
Similarly, we have another more geometric formula for \(\chi_\pi(f)\).
The equivalence of the two character formulas is equivalent to the
Explicit Formula of André Weil (\cites{Weil:Explicit_Formula_first,
  Weil:Explicit_Formula}).  It contains the following ingredients.
Let~\(\Discriminant_\GF\) be the discriminant of~\(\GF\).  For any
place~\(v\) of~\(\GF\), we write \(d\inv x\) for the normalized Haar measure
on~\(\GF\inv_v\).  The normalization condition is that
\(\int_{1\le\abs{x}<\abs{u}} d\inv x=\ln{}\abs{u}\) for all
\(u\in\GF\inv_v\) with \(\abs{u}\ge1\).  Consider the distribution
\(W\in\CCINF'(\ICL_\GF)\) defined by
\begin{displaymath}
  W(f) \defeq
  \sum_{v\in\Places(\GF)} \int'_{\GF\inv_v}
  \frac{f(x)\abs{x}}{\abs{1-x}} \,d\inv x
  - f(1)\cdot \ln{}\Discriminant_\GF
\end{displaymath}
with certain principal values \(\int'\).  We prove that~\(W\) is equal to
the character of~\(\pi\) and hence to \(\sum_{\omega\in\Irrep(\ICL_\GF)}
\ord(\omega,L_\GF) \cdot \hat{f}(\omega)\).  We call~\(W\) the
\emph{(raw) Weil distribution}.  The distribution
\[
f\mapsto \hat{f}(\abs{x}^{-1/2}) + \hat{f}(\abs{x}^{1/2})
- W(f\cdot \abs{x}^{-1/2}) = \chi_{\pi_-}(f\abs{x}^{-1/2}),
\]
which is usually called Weil distribution, is called \emph{modified
  Weil distribution} here.  Our character formula yields
\begin{displaymath}
  \hat{f}(\abs{x}^{-1/2}) + \hat{f}(\abs{x}^{1/2})
  - W(f\cdot \abs{x}^{-1/2})
  =
  \sum_{\omega\in\Irrep(\ICL_\GF)}
  \mult(\omega\abs{x}^{-1/2},\pi_-).
\end{displaymath}
Hence the modified Weil distribution is positive definite if and only
if the spectrum of~\(\pi_-\) is contained in the \emph{critical line}
\(\abs{x}^{1/2}\widehat{\ICL_\GF}\).  This is equivalent to the
Generalized Riemann Hypothesis.  Since we cannot prove this
positivity, anyway, we have no reason to modify~\(W\).  The raw Weil
distribution~\(W\) appears more naturally in our approach.  In addition,
it has a good geometric interpretation in connection with a
generalized Lefschetz formula (see~\cite{Connes:Trace_Formula}).

Although we state most results only for algebraic number fields, they
also work with small modifications for global function fields.  In
that case, the representation~\(\pi\) is an admissible representation on
a vector space with a countable basis.

For algebraic number fields, each \(L\)\nbd{}function has infinitely
many zeros, so that the sum defining the character is infinite.  Thus
we need analysis in order to define and to compute the character.  We
can modify the representation~\(\pi\) so that it lives on a Banach
space, even a Hilbert space.  However, there is little point in doing
so because we cannot make it unitary once there are multiple zeros or
zeros off the critical line.  Instead, the representation~\(\pi\) lives
on a nuclear bornological vector space.  We recall some basic facts
about nuclearity and operators of order~\(0\) in
Section~\ref{sec:nuclear_summable}.  There we also define summable
representations of locally compact groups on bornological vector
spaces.  We use bornologies instead of topologies because general
representation theory works better in that context (see
also~\cite{Meyer:Smooth}).  We develop the basic theory of summable
representation, although we need only very little of it for our
spectral interpretation, because we expect these representations to be
of independent interest.  The only facts that we really use are that
the character of a summable representation can be computed both
spectrally and tracially.

The representation~\(\pi\) is constructed as follows.  One main
ingredient is the space \(\Sch(\Adel_\GF)\) of Bruhat-Schwartz functions
on~\(\Adel_\GF\) (see~\cite{Bruhat:Distributions}).  The
group~\(\Adel\inv_\GF\) acts smoothly on \(\Sch(\Adel_\GF)\) by the
representation \(\lambda_g f(x)\defeq f(g^{-1}x)\).  Since we
want a representation of the quotient group~\(\ICL_\GF\), we replace
\(\Sch(\Adel_\GF)\) by the coinvariant space
\begin{equation}  \label{eq:def_Hilm_p}
  \Hilm_+\defeq \Sch(\Adel_\GF)/\GF\inv,
\end{equation}
which is the quotient of \(\Sch(\Adel_\GF)\) by the closed subspace
generated by elements of the form \(\lambda_g f-f\) with \(g\in\GF\inv\),
\(f\in\Sch(\Adel_\GF)\).  The representation~\(\lambda\)
of~\(\Adel\inv_\GF\) on \(\Sch(\Adel_\GF)\) descends to a smooth
representation of~\(\ICL_\GF\) on~\(\Hilm_+\).

We now describe~\(\Hilm_+\) more explicitly.  Let~\(S\) be a sufficiently
large finite set of places.  Let \(\Cont{S}\defeq \Places(\GF)\setminus
S\) and let \(\MCR\inv_{\Cont{S}}\) be the maximal compact subgroup of
the \(\Cont{S}\)\nbd{}idele group \(\Adel\inv_{\Cont{S}}\).  Let
\(\Hilm_+^S\subseteq\Hilm_+\) be the subspace of elements that are
invariant under the action of \(\MCR\inv_{\Cont{S}}\).  As a smooth
representation, \(\Hilm_+\) is the direct union of the
subspaces~\(\Hilm_+^S\), so that it suffices to describe those.  The
space~\(\Hilm_+^S\) carries a representation of the \(S\)\nbd{}idele class
group \(\ICL_S\defeq \Adel\inv_S/\GF\inv_S\).

Let \(\Sch(\ICL_S)\) be the Bruhat-Schwartz algebra of~\(\ICL_S\).  Let
\(\abs{x}\colon \ICL_S\to\R\inv_+\) be the norm homomorphism
on~\(\ICL_S\).  For an interval \(I\subseteq\R\), let
\[
\Sch(\ICL_S)_I \defeq
\{f\colon \ICL_S\to\C\mid
  \text{\(f\abs{x}^\alpha \in\Sch(\ICL_S)\) for all \(\alpha\in I\)}
\}.
\] Let \(\Fourier=\Fourier_S\) be the Fourier transform on~\(\Adel_S\) and
define \(Jf(x)\defeq \abs{x}^{-1} f(x^{-1})\).  Both \(\Fourier\) and~\(J\)
are isometries on the Hilbert space \(L^2(\Adel_S,dx)=
L^2(\Adel\inv_S,\abs{x}\,d\inv x)\).  Here \(dx\) and \(d\inv x\) denote
suitably normalized Haar measures on the locally compact groups
\(\Adel_S\) and~\(\Adel\inv_S\), respectively.  The operator \(\EFou\defeq
\Fourier J\) is called the equivariant Fourier transform on~\(\Adel_S\)
because it is \(\lambda\)\nbd{}equivariant.  Its inverse is \(\EFoi\defeq
J \Fourier^\ast\).  The operator~\(\EFou\) descends to a unitary operator
on \(L^2(\ICL_S,\abs{x}\,d\inv x)\).  This is a bounded operator on
\(\Sch(\ICL_S)_I\) for any \(I\subseteq\ooival{0,\infty}\).  Similarly,
\(\EFoi\) is a bounded operator on \(\Sch(\ICL_S)_I\) for \(I\subseteq
\ooival{-\infty,1}\).  For \(I\subseteq\ooival{0,1}\) both \(\EFou\)
and~\(\EFoi\) are bounded and hence bornological isomorphisms.  We let
\[
\Twist(\ICL_S) \defeq
\Sch(\ICL_S)_{\ooival{0,\infty}} \cap
\EFou(\Sch(\ICL_S)_{\ooival{-\infty,1}})
=
\bigcap_{\alpha>0} \Sch(\ICL_S)_\alpha \cap
\bigcap_{\alpha<1} \EFou(\Sch(\ICL_S)_\alpha).
\]
We show that
\[
\Hilm_+^S \cong
\Sch(\Adel_S)/\GF\inv_S\cong
\Twist(\ICL_S).
\]
This is not completely straightforward.  Our proof uses total derived
functors and shows at the same time that the higher group homology
\(H_n(\GF\inv,\Sch(\Adel_\GF))\) for \(n\ge1\) vanishes.  This is relevant in
order to compute the cyclic type homology theories of the crossed product
\(\GF\inv\cross\Sch(\Adel_\GF)\) (see~\cite{Meyer:Adelic}).  If we defined
\(\Hilm_+^S\defeq\Twist(\ICL_S)\) right away, we could construct our spectral
interpretation without computing \(\Sch(\Adel_\GF)/\GF\inv\).  We prefer the
above definition because we want our constructions to be as canonical as
possible.

Another more trivial ingredient of~\(\pi\) is the regular representation
of~\(\ICL_\GF\) on
\begin{equation}  \label{eq:def_Hilm_m}
  \Hilm_-\defeq \Sch(\ICL_\GF)_{\ooival{-\infty,\infty}}.
\end{equation}
Roughly speaking, \(\pi\) is the difference of the representations
\(\Hilm_+\) and~\(\Hilm_-\).  However, we have to take away a common
subrepresentation \(\Hilm_+\cap\Hilm_-\) in order to get a summable
virtual representation.  We do this as follows.  The summation map
\(\SUM f(x)\defeq \sum_{a\in\GF\inv} f(ax)\) for \(x\in\ICL_\GF\),
\(f\in\Sch(\Adel_\GF)\), descends to a bounded map from~\(\Hilm_+\) to
\(\Sch(\ICL_\GF)_{\ooival{1,\infty}}\).  We shall use the space
\begin{equation}  \label{eq:def_ppmm}
  \Sch(\ICL_\GF)_\ppmm \defeq
  \Sch(\ICL_\GF)_{\ooival{1,\infty}} \oplus
  \Sch(\ICL_\GF)_{\ooival{-\infty,0}},
\end{equation}
which intentionally leaves out the critical strip.  We define
\begin{equation}  \label{eq:def_i_pm}
  \begin{alignedat}{2}
    i_+&\colon \Hilm_+ \to \Sch(\ICL_\GF)_\ppmm,
    &\qquad f&\mapsto (\SUM f, J\SUM\Fourier f),
    \\
    i_-&\colon \Hilm_- \to \Sch(\ICL_\GF)_\ppmm,
    &\qquad f&\mapsto (f,f),
  \end{alignedat}
\end{equation}
and show that both \(i_+\) and~\(i_-\) are \(\ICL_\GF\)\nbd{}equivariant
bornological embeddings.  This is non-trivial because we left out the
critical strip in \(\Sch(\ICL_\GF)_\ppmm\).  In the following, we
view~\(\Hilm_\pm\) as subrepresentations of \(\Sch(\ICL_\GF)_\ppmm\)
via~\(i_\pm\).  We let
\begin{align*}
  \Hilm^0_+ &\defeq
  (\Hilm_+ + \Hilm_-)/\Hilm_- \cong \Hilm_+/(\Hilm_+\cap\Hilm_-),
  \\
  \Hilm^0_- &\defeq
  (\Hilm_+ + \Hilm_-)/\Hilm_+ \cong \Hilm_-/(\Hilm_+\cap\Hilm_-).
\end{align*}
The regular representations of~\(\ICL_\GF\) on~\(\Hilm_\pm\) descend to
smooth representations~\(\pi_\pm\) of~\(\ICL_\GF\) on~\(\Hilm^0_\pm\).  We
call \(\pi=\pi_+\ominus\pi_-\) the \emph{global difference
  representation}.  It is the desired spectral interpretation for the
poles and zeros of~\(L_\GF\).

The Poisson Summation Formula implies that the image of~\(i_+\) is
almost contained in the image of~\(i_-\).  More precisely, we find that
\begin{displaymath}
  \Hilm_+ + \Hilm_- =
  \{(f_0,f_1)\in
    \Sch(\ICL_\GF)_\ppmm \mid
    \text{\(f_0 - f_1 = c_0 - c_1 \abs{x}^{-1}\) for some \(c_0,c_1\in\C\)}
  \}.
\end{displaymath}
This means that the representation~\(\pi_+\) is \(2\)\nbd{}dimensional and
that \(\spec\pi_+= \{1,\abs{x}\}\).  These are exactly the two poles
of~\(L_\GF\).  We identify the spectrum of~\(\pi_-\) with the set of zeros
of~\(L_\GF\) as follows.  The Fourier-Laplace transform is an isomorphism
from~\(\Hilm_-\) onto the space of holomorphic functions on
\(\Irrep(\ICL_\GF)\) whose restriction to
\(\widehat{\ICL_\GF}\abs{x}^\alpha\) is a Schwartz function for any
\(\alpha\in\R\).  The subspace \(\Hilm_-\cap\Hilm_+\subseteq\Hilm_-\) is
mapped to an ideal in this space of holomorphic functions.  It consists
of all holomorphic functions that have a zero of order at least
\(\ord(\omega,L_\GF)\) at the zeros of~\(L_\GF\) and that, in addition,
satisfy some growth conditions on \(\widehat{\ICL_\GF}\abs{x}^\alpha\).
This allows us to identify the spectrum of~\(\pi_-\) with the zero set
of~\(L_\GF\).  Hence the spectrum of~\(\pi\) is contained in the
\emph{critical strip} \(\{\omega\in\Irrep(\ICL_\GF)\mid
\RE\omega\in[0,1]\}\) and has the same symmetries
\(\omega\mapsto\abs{x}\omega^{-1}\) and \(\omega\mapsto\conj\omega\) as the
\(L\)\nbd{}function.

We prove the summability of~\(\pi\) and compute its character tracially
using certain auxiliary operators~\(P_\pm\) on the space
\(\Sch(\ICL_\GF)_\ppmm\).  The operator~\(P_-\) is a projection onto the
subspace~\(\Hilm_-\) and~\(P_+\) is approximately a projection
onto~\(\Hilm_+\).  It depends on the choice of a sufficiently large
finite set of places~\(S\), so that we also denote it by~\(P_+^S\).  Our
character computation has two ingredients.  One of them is purely
\(S\)\nbd{}local.  We discuss this piece in Sections
\ref{sec:local_trace_formula} and~\ref{sec:local_trace_formula_Sch}.
We call the resulting formula the \emph{local trace formula}.  Its
main ingredient is the Fourier transform.  Roughly speaking, it
measures the failure of the Poisson Summation Formula in the
\(S\)\nbd{}local situation.  In the \(S\)\nbd{}local case, we can actually
work in the Hilbert space \(L^2(\ICL_S,\abs{x}\,d\inv x)\), whose
Fourier-Laplace transform is the space of square-integrable functions
on the critical line \(\RE\omega=1/2\).  The local trace formula on this
space is essentially due to Alain Connes, though his formulation is
slightly different.  However, in the global situation we have to use
the much smaller space \(\Twist(\ICL_S)\) instead of
\(L^2(\ICL_S,\abs{x}\,d\inv x)\).  Not surprisingly, the trace does not
depend on which space we use.  This means that the problem of
non\brd{}critical zeros and poles does not yet arise in the
\(S\)\nbd{}local case.  However, the character of the global difference
representation has another contribution from the places outside~\(S\).
If we fix \(f\in\CCINF(\ICL_\GF)\), then this contribution
annihilates~\(f\) for sufficiently large~\(S\).  Thus we can get rid of
this contribution rather easily.  Nevertheless, it prevents us from
controlling the position of the zeros.

Our construction follows the idea of Tate's Thesis (\cite{Tate:Thesis})
somewhat further and rephrases even more of the theory of
\(L\)\nbd{}functions in terms of representation theory and the Fourier
transform.  John Tate proved the meromorphic continuation of
\(L\)\nbd{}functions and their functional equations in this fashion.  We
also get a spectral interpretation for the poles and zeros and André
Weil's Explicit Formula.  The application of \(L\)\nbd{}functions to prime
number theory only uses their poles and zeros.  The relationship between
this spectral data and prime ideals is encoded by the explicit formulas
of prime number theory.  Hence it is not surprising that we can get the
Prime Number Theorem and similar results using the representation~\(\pi\)
instead of \(L\)\nbd{}functions.  In fact, we do not even have to define
\(L\)\nbd{}functions in order to prove the classical Prime Number Theorem.
To make this point, we prove an appropriate generalization of it in the
last Section~\ref{sec:PNT}.  The argument is closely related to the
proof in~\cite{Patterson:Riemann_zeta}, although it looks rather
different.

First we show that \(\spec\pi_-\) contains no points on the boundary of
the critical strip.  Our argument is a reformulation of an argument by
Pierre Deligne (\cite{Deligne:Weil2}), which in turn is an abstract
version of the classical one by Hadamard and de la Vallée Poussin.  I am
grateful to Christopher Deninger for pointing out this reference to me.
We clarify Deligne's proof slightly by using the Bohr compactification
of~\(\ICL_\GF\) and the GNS construction.  Then we deduce the following
generalization of the Prime Number Theorem.  We let~\(S\) be a
sufficiently large finite set of places and view the set of prime ideals
of the ring~\(\GF_S\) as a discrete subset of the \(S\)\nbd{}idele class
group~\(\ICL_S\).  Let \(A\subseteq\ICL\one_S\) be an open subset and let
\(\xi\in\R_{\ge1}\).  Let \(\pi_A(\xi)\) be the number of prime ideals
contained in \(A\times[1,\xi]\subseteq\ICL_S\).  We show that
\[
\pi_A(\xi) \approx \vol(A)\cdot \frac{\xi}{\ln\xi}
\qquad\text{for \(\xi\to\infty\).}
\]
Our argument is optimized to yield this result as quickly as
possible.  We do not estimate the vertical and horizontal position of
the spectrum of~\(\pi_-\), so that we get no error terms.

While it is arguable whether replacing \(L\)\nbd{}functions completely
by representations is a good idea in the rather classical situation of
the Prime Number Theorem, a purely representation theoretic approach
should be valuable in the automorphic case.  One important feature of
our constructions is that we work with all \(L\)\nbd{}functions
simultaneously and thus need not know much about the representation
theory of~\(\ICL_\GF\).  This simplifies several statements and proofs.
For instance, the Explicit Formula for a single \(L\)\nbd{}function is
more complicated to state than the one for all \(L\)\nbd{}functions
simultaneously.  This feature should become more relevant in the
automorphic case, where we know much less about the space of
automorphic representations.  While many of our constructions carry
over without change to the automorphic case, some new difficulties
also arise.  Since this article is already sufficiently long and
complicated, I have decided to leave the automorphic case for later
and to consider only representations of the idele class group here.

\section{Nuclearity and summable representations of locally compact
  groups}
\label{sec:nuclear_summable}

We define smooth summable representations of locally compact groups on
complete convex bornological vector spaces in such a way that for
totally disconnected groups we obtain exactly the well-known class of
admissible representations.  The notion of admissible representation
that is used in~\cite{Godement-Jacquet:Zeta_Simple} for adelic groups
over algebraic number fields is inappropriate for our purposes.  The
objective of that theory is to classify irreducible representations.
For this it is convenient to replace representations by purely
algebraic objects.  Our task is to make sense of infinite sums that
occur in our trace computations.  Thus we cannot avoid doing
functional analysis.  Also, we only consider representations of the
idele class group, which is Abelian, so that the classification of
irreducible representations is trivial.

We work with (complete convex) bornological vector spaces because
they provide the best setting for studying smooth representations of
locally compact groups (see~\cite{Meyer:Smooth}).  Roughly
speaking, such spaces are direct unions of Banach spaces whereas
complete locally convex topological vector spaces are inverse limits
of Banach spaces.  In both kinds of spaces, all analysis is eventually
reduced to analysis in Banach spaces using the limit structure.  This
is easier in the bornological case because inverse limits are more
subtle than direct unions.  The bornological vector spaces that we
need are either Fréchet spaces or strict direct unions of Fréchet
spaces, that is, LF-spaces.  For Fréchet spaces, the bornological and
topological ways of doing analysis are equivalent in many respects
(see~\cite{Meyer:Born_Top}).

In order to define summable representations we first have to discuss
nuclear operators and \(p\)\nbd{}summable operators for \(p<1\).  All this
is already contained in Alexander Grothendieck's ground-breaking
mémoir~\cite{Grothendieck:Produits_tensoriels}.  We briefly recall some
important ideas from~\cite{Grothendieck:Produits_tensoriels} for the
benefit of the reader.  We need this mainly in order to compare the
spectral and tracial characters of a summable representation.  Already
for nuclear operators on Banach spaces it is not always the case that
the trace is equal to the sum of the eigenvalues.  For this we need
operators of order~\(0\).  There is no difference between nuclear
operators and operators of order~\(0\) on nuclear spaces.  Hence it is
nice to know that the underlying bornological vector space of a smooth
summable representation is necessarily nuclear and regular.

Smooth summable representations are of some interest in their own
right.  Many things that can be done with admissible
\((\mathfrak{g},K)\)\brd{}modules can also be done with smooth summable
representations.  This is desirable because the latter are more
directly related to the objects of interest.  Therefore, we develop
some general theory of smooth summable representations here, although
we do not need it for our later applications.  We also discuss a few
issues related to the classification of irreducible representations,
but we do not enter this subject seriously.

\subsection{Nuclear operators and operators of order~\(0\)}
\label{sec:bornological_nuclear}

In the following, all bornological vector spaces are tacitly assumed
complete.  Most spaces that we deal with are convex.  We also allow
non-convex spaces because the operator ideals \(\ell^p(V,W)\) for
\(p\in\ooival{0,1}\) may fail to be convex.

A \emph{disk} in a vector space~\(V\) is an absolutely convex subset
\(S\subseteq V\) such that \(\bigcap_{t>1} tS=S\).  Given a disk \(S\subseteq
V\), we let \(V_S=\C\cdot S\) be its linear span.  There is a unique norm
on~\(V_S\) with closed unit ball~\(S\).  We always equip~\(V_S\) with this
norm.  The disk~\(S\) is called \emph{complete} if~\(V_S\) is a Banach
space.  The (complete) disked hull of an arbitrary subset is the
smallest (complete) disk in~\(V\) that contains it.

Let \(V,W\) be two convex bornological vector spaces.  We equip the dual
space~\(V'\) with the equibounded bornology and let \(\ell^1(V,W)\defeq
W\hot V'\).  Since tensor products commute with direct limits, we have
\(\ell^1(V,W)=\varinjlim W_S\hot V_T'\), where \(S\) and~\(T\) run through
the complete bounded disks in~\(W\) and~\(V'\), respectively.  There is a
natural bounded linear map from~\(V\) to the Banach space \((V'_T)'\) and
\(W_S\hot V_T'\subseteq \ell^1((V'_T)',W_S)\) with equality if~\(V'_T\) is
reflexive.  Thus the study of \(\ell^1(V,W)\) for general convex
bornological vector spaces \(V\) and~\(W\) reduces to the special case
where \(V\) and~\(W\) are Banach spaces.  Therefore, the results
of~\cite{Grothendieck:Produits_tensoriels} can be translated to convex
bornological vector spaces rather easily.

We denote the element \(w\otimes l\in W\hot V'\) also in the more
suggestive ket-bra notation as \(\ket{w}\bra{l}\).  Let \(V\) and~\(W\) be
Banach spaces.  Then \(\ell^1(V,W)\) is a Banach space as well and any
\(x\in\ell^1(V,W)\) can be written as
\begin{equation}  \label{eq:nuclear_operator}
  x = \sum_{n\in\N} \lambda_n \ket{w_n}\bra{l_n}
\end{equation}
with \((\lambda_n)\in\ell^1(\N)\) and \(w_n\in W\), \(l_n\in V'\) with
\(\norm{w_n}=\norm{l_n}=1\) for all \(n\in\N\).  The norm of~\(x\) is the
infimum of \(\norm{(\lambda_n)}_{\ell^1(\N)}\) over all such
presentations of~\(x\).  If \(V\) and~\(W\) are arbitrary convex
bornological vector spaces, we have a similar representation with
bounded sequences \((w_n)\) and \((l_n)\) in \(W\) and~\(V'\), respectively.
A subset of \(\ell^1(V,W)\) is bounded if its elements can be
represented in this form for two fixed bounded sequences \((w_n)\) and
\((l_n)\) and with \((\lambda_n)\) in a bounded subset of \(\ell^1(\N)\).

We can represent elements of \(\ell^1(V,W)\) as operators \(V\to W\) by
extending the formula \(\ket{w}\bra{l}(v)\defeq w\cdot \braket{l}{v}\).
The rules \(T\circ \ket{w}\bra{l}\defeq \ket{Tw}\bra{l}\) and
\(\ket{w}\bra{l}\circ T\defeq \ket{w}\bra{l\circ T}\) define bounded
bilinear composition maps
\[
\Hom(W,X)\times \ell^1(V,W)\to \ell^1(V,X),
\qquad
\ell^1(W,X)\times \Hom(V,W)\to \ell^1(V,X),
\]
which are compatible with the representation
\(\ell^1(V,W)\to\Hom(V,W)\).  Combining these two constructions, we can
define a multiplication
\[
\ell^1(W,X)\times \ell^1(V,W)\to \ell^1(V,X),
\qquad
\ket{w}\bra{l}\circ \ket{w'}\bra{l'} \defeq
\braket{l}{w'} \cdot \ket{w}\bra{l'}.
\]
This turns \(\ell^1(V)\defeq\ell^1(V,V)\) into a convex bornological
algebra.  The trace on \(\ell^1(V)\) is defined by extending
\(\tr(\ket{w}\bra{l}) \defeq \braket{l}{w}\) to a bounded linear map
\(\ell^1(V)\to\C\).  This is an \(\End(V)\)\brd{}bimodule trace on
\(\ell^1(V)\), that is, \(\tr(AB)=\tr(BA)\) for all \(A\in\ell^1(V)\),
\(B\in\End(V)\).

To define the above algebraic structure on \(\ell^1(V)\), we only need
the canonical pairing \(V'\times V\to\C\).  More generally, if we have
any (non-zero) bounded bilinear map \(\braket{\blank}{\blank}\colon
W\times V\to\C\), then \(V\hot W\) becomes a convex bornological algebra
equipped with a canonical trace.  For instance, if~\(V\) is a smooth
representation of a locally compact group~\(G\), we may replace~\(V'\) by
the contragradient \(\tilde{V}\defeq \smooth_G V'\), where \(\smooth_G
V'\) denotes the smoothening of~\(V'\), see~\cite{Meyer:Smooth}.  We
discuss algebras of the form \(V\hot W\) in greater detail in
Section~\ref{sec:irreducible}.

The range of the representation \(i_{V,W}\colon \ell^1(V,W)\to
\Hom(V,W)\) is called the space of \emph{nuclear operators} \(V\to W\).
The main problem with nuclear operators is that the
representation~\(i_{V,W}\) need not be faithful.  The only thing we know
is that if \(i_{V,W}(T)=0\), then \(T\circ S=0\) and \(S\circ T=0\) for all
\(S\in\ell^1(X,V)\) and \(S\in\ell^1(W,X)\), respectively.  In particular,
the kernel of \(i_{V,V}\) is a nilpotent ideal in \(\ell^1(V)\).
Fortunately, the representation~\(i_{V,V}\) is faithful if~\(V\) has
Grothendieck's approximation property.  Since all the spaces we need
have this property, we may sometimes forget about the distinction
between nuclear operators and \(\ell^1(V,W)\).

If~\(V\) is a Hilbert space, then \(\ell^1(V)\) is the usual ideal of trace
class operators.  Weyl's Trace Formula asserts that the trace of a
nuclear operator on Hilbert space is equal to the sum of its
eigenvalues.  For general~\(V\), we do not even know whether the trace
vanishes on \(\Ker i_{V,V}\).  That is, it is conceivable that an element
of \(\ell^1(V)\) with non-zero trace represents the zero operator, so that
Weyl's Trace Formula fails rather dramatically.  In order to extend it
to operators on bornological vector spaces, we have to restrict
attention to operators of order~\(0\).  This idea is due to Alexander
Grothendieck (\cite{Grothendieck:Produits_tensoriels}).

Let \(p\in\ocival{0,1}\).  If \(V\) and~\(W\) are Banach spaces, we define
\(\ell^p(V,W)\) by requiring \((\lambda_n)\in\ell^p(\N)\)
in~\eqref{eq:nuclear_operator}.  It is no longer a Banach space
because its unit ball is only \(p\)\nbd{}convex.  We may characterize it
as the completed \(p\)\nbd{}convex projective tensor product \(W\hot^p
V'\) of \(W\) and~\(V'\).  We extend~\(\hot^p\) to complete convex
bornological vector spaces by the usual inductive limit recipe,
\(V\hot^p W\defeq \varinjlim V_S\hot^p W_T\), where \(S\) and~\(T\) run
through the complete bounded disks in \(V\) and~\(W\), respectively.
Using the evident maps \(V\hot^p W\to V\hot^q W\) for \(p\le q\), we form
the projective limit \(V\hot^{[0]} W\defeq \varprojlim V\hot^p W\).  Let
\(\ell^p(V,W)\defeq W\hot^p V'\) and
\[
\ell^{[0]}(V,W)\defeq W\hot^{[0]} V' =
\varprojlim \ell^p(V,W).
\]
Elements of \(\ell^{[0]}(V,W)\) are called \emph{operators of
order~\(0\)}.  The spaces \(\ell^{[0]}\) and \(\ell^p\) carry the same
algebraic structure as \(\ell^1\).  The following theorems are proved by
Grothendieck for Banach spaces.  The extension to bornological vector
spaces is easy.

\begin{theorem}  \label{the:p_summable_inject}
  The canonical map \(\ell^p(V,W)\to\Hom(V,W)\) is injective for all
  convex bornological vector spaces \(V\) and~\(W\) if
  \(p\in\ocival{0,2/3}\).
\end{theorem}

In particular, the representation \(\ell^{[0]}(V)\to\End(V)\) is
faithful, so that we can really view elements of \(\ell^{[0]}(V)\) as
operators on~\(V\).

\begin{theorem}  \label{the:composition_order_p}
  Let \(S\subseteq\ell^1(V)\) be bounded and let \(n\ge4\).  Then
  \[
  S^{\circ n} \defeq
  \{x_1\circ\dots\circ x_n\mid x_1,\dots,x_n\in S\}
  \]
  is a bounded subset of \(\ell^p(V)\) for \(p=2/(n-1)\).
\end{theorem}

If \(S\subseteq V\) and \(T\subseteq V'\) are bounded disks, then the
pairing \(V\times V'\to\C\) restricts to \(V_S\times V'_T\), so that
\(V_S\hot V'_T\) is a Banach algebra.  We also obtain a natural
representation \(V_S\hot V'_T\to \End(V_S)\).  Since the kernel of this
map is always nilpotent, it is irrelevant for spectral considerations.
It is clear that the range of \(V_S\hot V'_T\) is contained in the ideal
of compact operators on \(\End(V_S)\).  Hence elements of the algebra
\(V_S\hot V'_T\) have the same spectral properties as compact operators
on~\(V_S\).  Their spectrum is a discrete set with~\(0\) as only possible
accumulation point, and each non-zero point in the spectrum is an
eigenvalue with finite algebraic multiplicity.  The subspace \(\bigcup
\Ker (T-\lambda)^n\subseteq V\) is called the \emph{generalized
eigenspace} of~\(\lambda\).  This space is finite dimensional for
\(\lambda\neq0\) and its dimension
\[
\mult(\lambda, T) \defeq
\lim_{n\to\infty} \dim \Ker (T-\lambda)^n
\]
is called the \emph{algebraic multiplicity} of the
eigenvalue~\(\lambda\).

\begin{theorem}  \label{the:spectrum_order_zero}
  Let~\(V\) be a convex bornological vector space.  Then
  \[
  \tr T = \sum_{\lambda\in\C\inv} \mult(\lambda, T)\cdot \lambda
  \]
  for all \(T\in\ell^{2/3}(V)\).  This sum converges absolutely and
  uniformly for~\(T\) in a bounded subset of \(\ell^{2/3}(V)\).
\end{theorem}

We call~\(T\) \emph{quasi-nilpotent} if it has no non-zero eigenvalues.
These operators require special attention because they tend to create
trouble.

\begin{theorem}  \label{the:quasi_nilpotent}
  Let \(T\in\ell^{2/3}(V)\).  Then~\(T\) is quasi-nilpotent if and only if
  \(\tr(T^n)=0\) for all \(n\in\N_{\ge1}\).
\end{theorem}

\begin{theorem}  \label{the:order_zero_restrict}
  Let \(T\colon V\to V\) be an operator of order~\(0\) and let \(W\subseteq
  V\) be a \(T\)\nbd{}invariant closed subspace.  Then the
  restriction~\(T|_W\) of~\(T\) to~\(W\) and the operator~\(T|_{V/W}\) on
  \(V/W\) induced by~\(T\) are again of order~\(0\) and satisfy
  \[
  \tr T|_V =
  \tr T|_W + \tr T|_{V/W}.
  \]
  Furthermore, if~\(T\) runs through a bounded subset of
  \(\ell^{[0]}(V)\), then \(T|_W\) and \(T|_{V/W}\) run through bounded
  subsets of \(\ell^{[0]}(W)\) and \(\ell^{[0]}(V/W)\), respectively.
\end{theorem}

The definition of nuclear bornological vector spaces is dual to the
corresponding definition for topological vector spaces.

\begin{definition}  \label{def:nuclear_space}
  Let~\(V\) be a complete convex bornological vector space.  It is
  called \emph{nuclear} if every bounded linear map from a Banach
  space~\(B\) into~\(V\) is nuclear.  A sequence~\((v_n)\) in~\(V\) is called
  \emph{absolutely summable} if \((\epsilon_n^{-1} v_n)\) is bounded for
  some \((\epsilon_n)\in\ell^1(\N)\).  A subset is called
  \emph{absolutely summable} if it is contained in the complete disked
  hull of an absolutely summable sequence.
\end{definition}

See~\cite{Hogbe-Nlend-Moscatelli:Nuclear} for a discussion of
nuclearity for bornological vector spaces.  It is known that~\(V\) is
nuclear if and only if any bounded subset of~\(V\) is absolutely
summable.  If~\(V\) is nuclear, then any nuclear map into~\(V\) is already
of order~\(0\).  Thus any operator from a Banach space into a nuclear
space is of order~\(0\).

\subsection{The definition of summable representations}
\label{sec:definition_summable}

Let~\(G\) be a locally compact topological group, let~\(V\) be a
(complete) convex bornological vector space and let \(\pi\colon
G\to\Aut(V)\) be a continuous group representation with integrated form
\(\IN\pi\colon \CCINF(G)\to\End(V)\).

\begin{definition}  \label{def:summable_representation}
  We call~\(\pi\) \emph{summable} if~\(\IN\pi\) is a bounded homomorphism
  into the algebra of nuclear operators on~\(V\).
\end{definition}

\begin{proposition}  \label{pro:admissible_representations}
  Let~\(G\) be totally disconnected.  Then~\(\pi\) is smooth and summable
  if and only if the fixed point subspaces~\(V^k\) for compact open
  subgroups \(k\subseteq G\) are all finite dimensional and~\(V\) is the
  bornological direct union of these subspaces.
\end{proposition}

\begin{proof}
  Since~\(G\) is totally disconnected, \(\pi\) is smooth if and only if~\(V\)
  is the bornological direct union of the subspaces~\(V^k\)
  (see~\cite{Meyer:Smooth}).  The subspace~\(V^k\) is the range of an
  idempotent in \(\CCINF(G)\), namely, the normalized Haar measure on~\(k\).
  Summability implies that this projection is mapped to a nuclear
  projection on~\(V\).  This means that~\(V^k\) is finite dimensional.
  Conversely, since \(\CCINF(G)=\varinjlim \CCINF(G//k)\), the
  representation~\(\pi\) is summable once all~\(V^k\) are finite
  dimensional.
\end{proof}

Thus the underlying bornological vector space~\(V\) of a smooth
summable representation of a totally disconnected group automatically
carries the fine bornology\mdash that is, it is the direct union of
its finite dimensional subspaces\mdash and the representation of~\(G\)
is admissible in the usual sense
(see~\cite{Godement-Jacquet:Zeta_Simple}).

\begin{definition}  \label{def:character}
  Let \(\pi\colon G\to\Aut(V)\) be a summable representation.  Its
  \emph{character} is the distribution \(\chi_\pi\in\CCINF'(G)\) defined
  by \(\chi_\pi(f)\defeq \tr \IN\pi(f)\) for \(f\in\CCINF(G)\).
\end{definition}

\begin{theorem}  \label{the:summable_spectrum}
  Let~\(\pi\) be a summable representation.  The spectrum of \(\IN\pi(f)\)
  for \(f\in\CCINF(G)\) consists only of eigenvalues and~\(\{0\}\) and all
  eigenvalues have finite algebraic multiplicity.  We have
  \[
  \tr \IN\pi(f) =
  \sum_{\lambda\in\C\inv} \mult(\lambda, \IN\pi(f))
  \]
  with uniform absolute converge for~\(f\) in bounded subsets of
  \(\CCINF(G)\).
\end{theorem}

This follows immediately from Theorem~\ref{the:spectrum_order_zero}
and the following lemma:

\begin{lemma}  \label{lem:summable_order_zero}
  If the representation~\(\pi\) is summable, then~\(\IN\pi\) is a bounded
  map into \(\ell^{[0]}(V)\) and even into \(\smooth V\hot^{[0]}\smooth
  V'\subseteq \ell^{[0]}(V)\).
\end{lemma}

The \emph{smoothening} \(\smooth(V)\) of~\(\pi\) is defined
in~\cite{Meyer:Smooth}.

\begin{proof}
  We have to prove that~\(\IN\pi\) is a bounded map into
  \[
  \smooth(V) \hot^p \smooth(V') \cong
  (\CCINF(G) \hot_{\CCINF(G)} V) \hot^p
  (\CCINF(G) \hot_{\CCINF(G)} V')
  \]
  for all \(p\in\ocival{0,2/3}\).  Recall that the representation
  \(\ell^p(V)\to\End(V)\) is faithful for those values of~\(p\).  Fix such
  a~\(p\) and let \(n\in\N\) be such that \(n-2\ge 1+2/p\).  Fix a bounded
  subset \(S\subseteq \CCINF(G)\).  The \(n\)\nbd{}fold convolution is a
  bornological quotient map from \(\CCINF(G)^{\hot n}\) onto \(\CCINF(G)\)
  by results of~\cite{Meyer:Smooth}.  Hence~\(S\) is the image of a
  bounded subset \(S'\subseteq \CCINF(G)^{\hot n}\).  There exists a
  bounded disk \(T\subseteq \CCINF(G)\) such that all elements of~\(S'\)
  can be written as sums \(\sum_{j\in\N} \lambda_j x_{0,j} \otimes
  \dots \otimes x_{n+1,j}\) with \(x_{i,j}\in T\) for all \(i,j\) and
  \((\lambda_j)\) in a bounded subset of \(\ell^p(\N)\) (the nuclearity of
  \(\CCINF(G)\) allows us to choose \(p<1\)).
  Theorem~\ref{the:composition_order_p} implies that the compositions
  \(\IN\pi(x_{1,j})\circ\cdots\circ\IN\pi(x_{n,j})\) all lie in some
  bounded subset of \(\ell^p(V)=V\hot^p V'\).  Since
  \begin{equation}  \label{eq:INpi_equivariance}
    \IN\pi(f_1)\circ \ket{a}\bra{l}\circ \IN\pi(f_2) =
    \ket{\IN\pi(f_1)(a)}\bra{\IN\pi'(f_2)(l)},
  \end{equation}
  the compositions \(\IN\pi(x_{0,j})\circ\cdots\circ\IN\pi(x_{n+1,j})\)
  lie in a bounded subset of \(\smooth(V) \hot^p \smooth(V')\).  Since
  the latter space is \(p\)\nbd{}convex, we get the assertion.
\end{proof}

\begin{theorem}  \label{the:summable_hereditary}
  Summability is hereditary for smoothenings, contragradients and
  tensor product.  Under these operations, the character behaves as
  follows:
  \[
  \chi_V = \chi_{\smooth V},
  \qquad
  \chi_{\tilde V} = \check{\chi}_V,
  \qquad
  \chi_{V\hot W}\cong \chi_V\cdot\chi_W,
  \]
  where we define \(\check{f}(g)\defeq f(g^{-1})\) for all \(g\in G\),
  \(f\in\CCINF(G)\) and \(\braket{\check{\chi}_V}{f}\defeq
  \braket{\chi_V}{\check{f}}\).

  Let \(W\into V\prto V/W\) be an extension of continuous
  representations.  Then~\(V\) is summable if and only if both \(W\)
  and~\(V/W\) are summable.  In this case, we have \(\chi_V = \chi_W +
  \chi_{V/W}\).
\end{theorem}

\begin{proof}
  Suppose that~\(V\) is a summable representation.
  Lemma~\ref{lem:summable_order_zero} yields that~\(\IN\pi\) factors
  through \(\smooth(V)\hot^{[0]} \tilde{V}\).  The integrated forms of
  both \(\smooth(V)\) and~\(\tilde{V}\) can be expressed in terms
  of~\(\IN\pi\), which yields the summability of these representations.
  Summability is hereditary for tensor products because there are
  canonical bilinear maps \(\ell^p(V)\times\ell^p(W)\to \ell^p(V\hot
  W)\).  The character formulas are straightforward to prove.

  Consider now an extension as above.  Suppose first that~\(V\) is
  summable.  Then~\(W\) is mapped into itself by \(\IN\pi(h)\) for all
  \(h\in\CCINF(G)\).  The operators on~\(W\) and~\(V/W\) induced by~\(\IN\pi\)
  are nothing but the integrated forms of the subspace and quotient
  representations on \(W\) and~\(V/W\).  They are again of order~\(0\) if
  \(\IN\pi(h)\) is of order~\(0\) by
  Theorem~\ref{the:order_zero_restrict}.  Therefore, the
  representations on \(W\) and~\(V/W\) are summable as well.  We also
  obtain the additivity of the character.

  Suppose conversely that the representations on \(W\) and \(V/W\) are
  summable.  We claim that~\(V\) is summable as well.  It suffices to
  show that \(\IN\pi(f_1\ast f_2)\) is nuclear for all
  \(f_1,f_2\in\CCINF(G)\) because
  \(\CCINF(G)\hot_{\CCINF(G)}\CCINF(G)\cong\CCINF(G)\).  By hypothesis,
  \(\IN\pi(f_2)\) induces a nuclear operator on \(V/W\).  It is well-known
  that nuclear operators can be lifted in extensions of Banach spaces.
  This extends to bornological vector spaces because we can factor
  nuclear operators through nuclear operators on Banach spaces.  Hence
  there exists a nuclear operator \(X\colon V/W\to V\) such that \(P\circ
  X=\IN\pi(f_2)\), where~\(P\) denotes the canonical projection \(V\to
  V/W\).  We write
  \[
  \IN\pi(f_1\ast f_2) =
  \IN\pi(f_1)\circ (\IN\pi(f_2)-X\circ P) + \IN\pi(f_1)\circ X\circ P.
  \]
  By construction, the range of \(\IN\pi(f_2)-X\circ P\) is contained
  in~\(W\).  Since the representation on~\(W\) is summable,
  \(\IN\pi(f_1)\circ (\IN\pi(f_2)-X\circ P)\) is nuclear.  Since~\(X\) is
  also nuclear, we get the desired nuclearity of \(\IN\pi(f_1\ast
  f_2)\).
\end{proof}

\subsection{Summable representations of Abelian groups}
\label{sec:summable_Abelian}

Let~\(G\) be an Abelian locally compact group.  Let \(\Irrep(G)\) be the
space of \emph{quasi-characters} of~\(G\).  For any
\(\omega\in\Irrep(G)\), we let \(\C(\omega)\) be the corresponding
\(1\)\nbd{}dimensional representation of~\(G\).  The \emph{Fourier-Laplace
transform} \(f\mapsto\hat{f}\) maps \(\CCINF(G)\) to a certain algebra of
holomorphic functions on \(\Irrep(G)\).  We say that a summable
representation \(\pi\colon G\to\Aut(V)\) \emph{contains the
quasi-character~\(\omega\)} if there exists an equivariant map
\(\C(\omega)\to V\).  The \emph{spectrum} \(\spec \pi\) of~\(\pi\) is the
set of all quasi-characters contained in~\(V\).  (This definition is
only reasonable for summable representations!)

We call the representation~\(\pi\) \emph{quasi-nilpotent} if \(\IN\pi(f)\)
is quasi-nilpotent for all \(f\in\CCINF(G)\).  If~\(\pi\) is
quasi-nilpotent, then it cannot contain any quasi-character, that is,
\(\spec\pi=\emptyset\).  Conversely, if~\(\pi\) is not quasi-nilpotent,
then \(\spec\pi\neq\emptyset\).  To prove this, pick \(f\in\CCINF(G)\) for
which \(\IN\pi(f)\colon V\to V\) is not quasi-nilpotent.  That is,
\(\IN\pi(f)\) has a non-zero eigenvalue~\(\lambda\).  The
\(\lambda\)\nbd{}eigenspace of \(\IN\pi(f)\) is finite dimensional and
\(G\)\nbd{}invariant because~\(G\) is Abelian.  Hence~\(\pi\) has finite
dimensional invariant subspaces.  Standard linear algebra yields
that~\(\pi\) contains some quasi-character.

Any finite dimensional subrepresentation \(W\subseteq V\) has a
Jordan-Hölder series, whose subquotients are of the form \(\C(\omega)\)
for quasi-characters~\(\omega\).  The \emph{algebraic multiplicity}
\(\mult(\omega, \pi)\) of a quasi-character~\(\omega\) in the
representation~\(\pi\) is defined as the maximal number of times
\(\C(\omega)\) occurs as a subquotient in such a Jordan-Hölder series.
This multiplicity is finite for all \(\omega\in\Irrep(G)\) because
\[
\mult(\lambda,  \IN\pi(f)) =
\sum_{\{\omega\in\Irrep(G)\mid \hat{f}(\omega)=\lambda\}}
\mult(\omega, \pi)
\]
for all \(f\in\CCINF(G)\).  Hence the trace formula of
Theorem~\ref{the:summable_spectrum} can be rewritten as
\begin{equation}
  \label{eq:summable_character_spectral}
  \chi_\pi(f) \defeq
  \tr \IN\pi(f) =
  \sum_{\omega\in\Irrep(G)} \mult(\omega, \pi)\cdot \hat{f}(\omega).
\end{equation}

\subsection{Properties of summable representations}
\label{sec:properties_summable}

If a bornological vector space carries a smooth and summable group
representation, then it cannot be completely arbitrary.  Most
importantly, the space~\(V\) must be nuclear, so that there is no
difference between nuclear operators and operators of order~\(0\).  The
results of this section are not relevant for our number theoretic
applications because there the properties that we show to hold in
general are easy to verify directly.  Nevertheless, it seems
worthwhile to include them here together with our discussion of
summable representations.  The proof requires the following
preparatory lemma.

\begin{lemma}  \label{lem:nuclear_in_p_convex_tensor}
  Let \(X,V,W\) be convex bornological vector spaces.  Suppose~\(X\) to be
  nuclear.  Let \(\phi\colon X\to V\hot^{1/3} W\) be a bounded linear
  map and let \(S\subseteq X\) be bounded.  Then \(\phi(S)\) is contained
  in the complete disked hull of \(\{v_m\otimes w_n\mid m,n\in\N\}\)
  with absolutely summable sequences \((v_n)\) and \((w_n)\) in \(V\)
  and~\(W\), respectively.
\end{lemma}

\begin{proof}
  Let \(p=1/3\).  Since~\(X\) is nuclear, there are a bounded
  sequence~\((x_n)\) in~\(X\) and \((\epsilon_n)\in\ell^p(\N)\) such
  that~\(S\) is contained in the disked hull of \((\epsilon_n x_n)\).
  Since the sequence \(\phi(x_n)\) is bounded in \(V\hot^p W\), we have
  \(\phi(x_n)=\sum \lambda_{n,j} v_j \otimes w_j\) with \(\sum_j
  \abs{\lambda_{n,j}}^p\le1\) for all \(n\in\N\) and bounded sequences
  \((v_j)\) and \((w_j)\) in \(V\) and~\(W\), respectively.  Hence
  \(\phi(\epsilon_n x_n)\) is contained in the complete disked hull of
  the set of elementary tensors \((\lambda_{n,j}^p\epsilon_n^p
  v_j)\otimes (\lambda_{n,j}^p\epsilon_n^p w_j)\). Our construction
  guarantees that the bisequence \((\lambda^p_{n,j}\epsilon_n^p)\) is
  absolutely summable.  Therefore, the bisequences
  \((\lambda_{n,j}^p\epsilon_n^p v_j)\) and
  \((\lambda_{n,j}^p\epsilon_n^p w_j)\) in \(V\) and~\(W\) are absolutely
  summable.
\end{proof}

\begin{theorem}  \label{the:summable_nuclear}
  Let \(\pi\colon G\to\Aut(V)\) be smooth and summable.  Then~\(V\) is
  regular and nuclear.  For any smooth compact subgroup \(k\subseteq
  G\), the \(k\)\nbd{}fixed point subspace \(V^k\subseteq V\) is
  bornologically metrizable.
\end{theorem}

\begin{proof}
  If~\(V\) fails to be regular, then there is \(v\in V\) with \(l(v)=0\) for
  all \(l\in V'\).  Hence \(T(v)=0\) for all \(T\in\ell^1(V)\), so that
  \(\IN\pi(f)(v)=0\) for all \(f\in\CCINF(G)\).  This is impossible for a
  continuous representation.  For any bounded subset \(S\subseteq V\),
  there are bounded subsets \(S'\subseteq\CCINF(G)\), \(S''\subseteq V\)
  such that~\(S\) is contained in the complete disked hull of
  \(\IN\pi(S')(S'')\).  Recall that~\(\IN\pi\) is a bounded map from
  \(\CCINF(G)\) to \(V\hot^{[0]} V'\) by
  Lemma~\ref{lem:summable_order_zero}.  Since \(\CCINF(G)\) is nuclear,
  Lemma~\ref{lem:nuclear_in_p_convex_tensor} yields that \(\IN\pi(S')\)
  is contained in the complete disked hull of
  \((\ket{v_m}\bra{v_n'})_{m,n}\) for absolutely summable sequences
  \((v_m)\) and \((v_n')\) in \(V\) and~\(V'\), respectively.  It follows that
  the given set~\(S\) is absorbed by the complete disked hull of
  \((v_m)\).  That is, any bounded subset of~\(V\) is absolutely summable.
  This implies that~\(V\) is nuclear.
  
  Let \(L\subseteq G\) be a compact neighborhood of the identity in~\(G\)
  such that \(kLk=L\) and let \(\CINF_0(L//k)\subseteq\CCINF(G)\) be the
  subspace of \(k\)\nbd{}biinvariant functions supported in~\(L\).  The
  restriction of~\(\IN\pi\) to \(\CINF_0(L//k)\hot V^k\to V^k\) is a
  bornological quotient map because the linear section for \(\IN\pi\colon
  \CCINF(G)\hot V\to V\) constructed in~\cite{Meyer:Smooth} can
  evidently be chosen to map~\(V^k\) into \(\CINF_0(L//K,V^k)\).  The space
  \(\CINF_0(L//k)\) is a Fréchet space equipped with the precompact
  bornology.  Hence it is bornologically metrizable
  by~\cite{Meyer:Born_Top}.  Let \((S_n)_{n\in\N}\) be a sequence of
  bounded subsets of~\(V\).  For each \(n\in\N\), there are bounded disks
  \(S'_n\subseteq\CINF_0(L//k)\) and \(S''_n\subseteq V\) such that~\(S_n\) is
  contained in the complete disked hull of \(\IN\pi(S'_n)(S''_n)\).  Since
  \(\CINF_0(L//k)\) is metrizable, there is a bounded disk~\(S'\) that
  absorbs all~\(S'_n\).  Hence we may assume that \(S'_n=S'\) for all
  \(n\in\N\).  There are bounded disks \(T\subseteq V\), \(T'\subseteq V'\) so
  that \(\IN\pi(S')\) is contained in the complete disked hull of
  \(\ket{T}\bra{T'}\).  Hence \(S_n\subseteq \IN\pi(S')(S''_n)\subseteq
  T\braket{T'}{S''_n}\).  That is, all~\(S_n\) are absorbed by~\(T\).  This
  means that~\(V\) is bornologically metrizable.
\end{proof}

It would be nice if one could show that the spaces~\(V^k\) for a smooth
and summable representation are nuclear Fréchet spaces equipped with
the von Neumann bornology.  Although this happens in the examples that
I know of, I have no idea how to prove this in general.

\begin{proposition}  \label{pro:smooth_summable_reflexive}
  Let~\(V\) be a smooth summable representation.  Then the natural map
  from~\(V\) to its double contragradient is a bornological isomorphism.
\end{proposition}

\begin{proof}
  The integrated form of the contragradient representation~\(\tilde{V}\)
  is up to the reflection \(h\mapsto\check{h}\) the same as the
  integrated form of~\(V\) itself.  Thus the integrated form of the
  double contragradient~\(W\) is equal to the integrated form of~\(V\).
  Hence the integrated form \(\CCINF(G)\hot W\to W\) of the double
  contragradient representation factors through \(V\hot \tilde{V}\hot
  W\) and hence through the inclusion map \(V\to W\).  Since the map
  \(\CCINF(G,W)\to W\) is a bornological quotient map, so is the map
  \(V\to W\).  Since~\(V\) is regular, this map is injective, hence a
  bornological isomorphism.
\end{proof}

\subsection{Quasi-nilpotent representations}
\label{sec:quasi_nilpotent_rep}

\begin{definition}  \label{def:quasi_nilpotent_representation}
  A summable representation \(\pi\colon G\to\Aut(V)\) is called
  \emph{quasi\brd{}nilpotent} if \(\IN\pi(f)\) is quasi\brd{}nilpotent
  for all \(f\in\CCINF(G)\).
\end{definition}

If~\(G\) is totally disconnected, then the zero representation is the
only quasi-nilpotent summable continuous representation of~\(G\): in
this case, the operators \(\IN\pi(f)\) have finite rank and hence must
vanish once they are quasi-nilpotent.

\begin{proposition}  \label{pro:quasi_nilpotent_character}
  A summable representation is quasi-nilpotent if and only if its
  character vanishes.
\end{proposition}

\begin{proof}
  Since a quasi-nilpotent operator has no non-zero eigenvalues, a
  quasi-nilpotent representation must have character~\(0\).  Conversely,
  if \(\chi_\pi=0\), then \(\tr \IN\pi(f)^n=0\) for all \(f\in\CCINF(G)\),
  \(n\in\N\).  This implies that \(\IN\pi(f)\) is quasi-nilpotent by
  Theorem~\ref{the:quasi_nilpotent}.
\end{proof}

Together with Theorem~\ref{the:summable_hereditary}, this criterion
implies that quasi-nilpotence is hereditary for smoothenings,
contragradients, tensor products, subrepresentations, quotients and
extensions.

\begin{proposition}  \label{pro:split_off_quasi_nilpotent}
  Let \(\pi\colon G\to\Aut(V)\) be a summable representation.  Then
  there are closed \(G\)\nbd{}invariant subspaces \(N_0\subseteq
  V_0\subseteq V\) such that \(N_0\) and \(V/V_0\) are quasi-nilpotent and
  the quotient \(V_0/N_0\) is orthogonal to all quasi-nilpotent summable
  representations in the sense that
  \[
  \Hom_G(N,V_0/N_0)\cong 0,
  \qquad
  \Hom_G(V_0/N_0,N)\cong 0
  \]
  for all quasi-nilpotent summable representations~\(N\).
\end{proposition}

\begin{proof}
  Let \(V_0\subseteq V\) be the closed subspace that is generated by the
  generalized eigenspaces of the operators \(\IN\pi(f)\) with
  \(f\in\CCINF(G)\) for non-zero eigenvalues.  This subspace is clearly
  \(G\)\nbd{}invariant.  By construction, all the generalized eigenspaces
  of \(\IN\pi(f)\) on~\(V\) are already contained in~\(V_0\).  Hence the
  spectral computation of the trace of \(\IN\pi(f)\) yields the same
  answer for \(V_0\) and~\(V\).  Thus \(V_0\) and~\(V\) have the same
  character, so that the character of~\(V/V_0\) vanishes.  Thus~\(V/V_0\)
  is quasi-nilpotent.
  
  Suppose that \(f\colon V_0\to N\) is an equivariant map to some
  quasi-nilpotent representation~\(N\).  Then~\(f\) maps generalized
  eigenspaces of \(\IN\pi(f)\) again to generalized eigenspaces.  Since
  \(\IN\pi(f)\) is quasi-nilpotent on~\(N\), it has no non-zero
  eigenvalues and hence~\(f\) must vanish on generalized eigenspaces of
  \(\IN\pi(f)\) for non-zero eigenvalues.  Thus~\(f\) must vanish
  on~\(V_0\).  \emph{A fortiori}, any map \(V_0/N_0\to N\) vanishes,
  whatever the subspace \(N_0\subseteq V_0\) may be.
  
  We let \(W=V_0\) to simplify our notation.  We can do the above
  construction also for the contragradient representation
  \(\tilde\pi\colon G\to\Aut(\tilde{W})\) and obtain a certain subspace
  \(\tilde{W}_0\subseteq\tilde{W}\).  We let~\(N_0\) be the set of \(v\in
  W\) which are orthogonal to~\(\tilde{W}_0\), that is, \(l(v)=0\) for all
  \(l\in\tilde{W}_0\).  This is a closed \(G\)\nbd{}invariant subspace.
  Proposition~\ref{pro:smooth_summable_reflexive} implies that we have
  the usual duality between subrepresentations/quotients of~\(W\) and
  quotients/subrepresentations of~\(\tilde{W}\).  Since
  \(\tilde{N}_0\cong \tilde{W}/\tilde{W}_0\) is quasi-nilpotent, \(N_0\)
  is quasi-nilpotent as well.  The contragradient of \(V_0/N_0\) is
  isomorphic to~\(\tilde{W}_0\).  If \(f\colon N\to V_0/N_0\) is any
  \(G\)\nbd{}equivariant bounded map, then its contragradient\mdash
  which is a map from \(\tilde{W}_0\) to~\(\tilde{N}\)\mdash vanishes by
  construction of~\(\tilde{W}_0\) because~\(\tilde{N}\) is
  quasi-nilpotent.  Therefore, the map~\(f\) vanishes as well.  Thus the
  filtration \(N_0\subseteq V_0\subseteq V\) has the desired properties.
\end{proof}

In this way, we can split off the quasi-nilpotent part of a summable
representation.  I do not know whether our spectral interpretation for
the zeros of \(L\)\nbd{}functions has a quasi-nilpotent part.

\subsection{Notions of irreducibility}
\label{sec:irreducible}

A representation is called \emph{irreducible} if it has no closed
invariant subspaces.  We may want to classify the irreducible, smooth,
summable representations of a locally compact group~\(G\).  However, we
must exclude irreducible quasi-nilpotent representations.  For
example, \(\R\) has irreducible quasi-nilpotent representations because
there exist operators on infinite dimensional Banach spaces without
invariant subspaces.  However, we certainly have no hope of
classifying such operators.

\begin{lemma}[Schur's Lemma]  \label{lem:summable_Schur}
  Let \(\pi\colon G\to\Aut(V)\) be an irreducible summable
  representation that is not quasi-nilpotent.  Then any
  \(G\)\nbd{}equivariant operator \(V\to V\) must be a constant multiple
  of the identity operator.
\end{lemma}

\begin{proof}
  Let \(T\colon V\to V\) be an operator that commutes with \(\pi(G)\) and
  hence with \(\IN\pi(\CCINF(G))\).  Since~\(\pi\) is not quasi-nilpotent,
  there is \(f\in\CCINF(G)\) for which \(\IN\pi(f)\) has a non-zero
  eigenvalue~\(\lambda\).  The corresponding eigenspace is finite
  dimensional because \(\IN\pi(f)\) is nuclear.  It is
  \(T\)\nbd{}invariant because~\(T\) commutes with \(\IN\pi(f)\).  The
  restriction of~\(T\) to this subspace has an eigenvalue \(c\in\C\).
  Since~\(\pi\) is irreducible, the kernel of \(T-c\) cannot be a proper
  subspace, that is, \(T=c\cdot\ID_V\).
\end{proof}

The converse of Schur's Lemma only holds for unitary representations.
It is easy to find reducible representations of~\(\Z\) on~\(\C^2\) with
trivial commutant.

Let \(Z(G)\) be the center of the category of smooth representations or,
equivalently, the center of the multiplier algebra of \(\CCINF(G)\)
(see~\cite{Meyer:Smooth}).  It acts on any smooth representation by
central module homomorphisms.  It is equal to the Bernstein center for
totally disconnected groups.  If~\(G\) is a connected Lie group, then
\(Z(G)\) contains the center of the universal enveloping algebra of~\(G\).
The latter is equal to \(Z(G)\) for connected complex Lie groups with
trivial center and also for groups like \(\mathrm{PSl}(n,\R)\).
If~\(\pi\) is irreducible and Schur's Lemma applies, then \(Z(G)\) acts by
multiples of~\(\ID_V\).  Thus we obtain a character on \(Z(G)\), which we
call the \emph{central character} of~\(\pi\).

Recall that the integrated form of a smooth summable
representation~\(\pi\) is a bounded homomorphism \(\IN\pi\colon
\CCINF(G)\to V\hot\tilde{V}\) by Lemma~\ref{lem:summable_order_zero}.

\begin{definition}  \label{def:strongly_irreducible}
  A smooth summable representation~\(\pi\) is called \emph{strongly
  irreducible} if~\(\IN\pi\) is a bornological quotient map onto
  \(V\hot\tilde{V}\).
\end{definition}

It is well-known that irreducible admissible representations of
totally disconnected groups are strongly irreducible.  It is not so
clear what happens for general groups.  We certainly have to exclude
quasi-nilpotent representations\mdash which do not exist in the
totally disconnected case\mdash in order to get a true statement.

We use the following notion of Morita equivalence for bornological
algebras.  It is the straightforward extension of the corresponding
Banach space theory from \cites{Groenbaek:Morita,
Groenbaek:Imprimitivity}.  A bornological algebra~\(A\) is called
\emph{self-induced} if \(A\hot_A A\cong A\).  A bornological left
module~\(V\) over a self-induced bornological algebra~\(A\) is called
\emph{essential} if \(A\hot_A V\cong V\).  Two self-induced bornological
algebras \(A\) and~\(B\) are called \emph{Morita equivalent} if there
exist an essential bornological \(A,B\)\brd{}bimodule~\(P\) and an
essential bornological \(B,A\)\brd{}bimodule~\(Q\) such that \(P\hot_B
Q\cong A\) and \(Q\hot_A P\cong B\) as bimodules.  This implies that the
categories of essential modules over \(A\) and~\(B\) are isomorphic via
\(Q\hot_A\blank\) and \(P\hot_B\blank\).

\begin{lemma}  \label{lem:Morita_C}
  A bornological algebra is Morita equivalent to~\(\C\) if and only if
  it is isomorphic to the algebra \(V\hot W\) defined by some non-zero,
  bounded bilinear pairing \(W\hot V\to\C\).  The Morita equivalence is
  implemented by~\(V\) with the obvious structure of essential
  bornological \(V\hot W,\C\)\brd{}bimodule and by~\(W\) with the obvious
  structure of essential bornological \(\C,V\hot W\)\brd{}bimodule.  Any
  irreducible left or right module over \(V\hot W\) is isomorphic to \(V\)
  or~\(W\), respectively.
\end{lemma}

\begin{proof}
  It is straightforward to check that \(V\hot W\) is self-induced and
  that \(V\) and~\(W\) implement a Morita equivalence between \(\C\) and
  \(V\hot W\).  Since \(P\hot_\C Q\cong P\hot Q\), it is also evident that
  any bornological algebra Morita equivalent to~\(\C\) is of this form.
  The category of essential bornological modules over~\(\C\) is
  isomorphic to the category of bornological vector spaces.  The
  bornological vector space~\(\C\) is the unique irreducible object of
  this category.  Since the Morita equivalence gives rise to an
  isomorphism between the categories of modules, we find that \(V\)
  or~\(W\) are the unique irreducible left or right modules over \(V\hot
  W\), respectively.
\end{proof}

\begin{corollary}  \label{cor:strongly_irreducible}
  The map \(\pi\mapsto \Ker\IN\pi\) gives rise to a bijection between
  the strongly irreducible, smooth, summable representations of~\(G\)
  and the ideals in \(\CCINF(G)\) for which \(\CCINF(G)/\Ker\IN\pi\) is
  Morita equivalent to~\(\C\).
\end{corollary}

\begin{proof}
  By definition, \(\pi\) is strongly irreducible if and only if
  \(\CCINF(G)/\Ker\pi\cong V\hot\tilde{V}\).  The lemma implies that
  this is Morita equivalent to~\(\C\) and that \(V\) and~\(\tilde{V}\) as
  modules over \(\CCINF(G)/\Ker\IN\pi\) and hence as representations
  of~\(G\) are uniquely determined by \(\Ker\IN\pi\).  Conversely, if
  \(\CCINF(G)/I\) is Morita equivalent to~\(\C\) for some ideal in
  \(\CCINF(G)\), then \(\CCINF(G)/I\cong V\hot W\) for some essential left
  and right modules \(V\) and~\(W\) over \(\CCINF(G)\) together with a
  non-zero linear map \(W\hot_{\CCINF(G)/I} V\to\C\).  It is
  straightforward to see that~\(V\) is a strongly irreducible, smooth
  summable representation with associated kernel~\(I\).
\end{proof}

Thus strongly irreducible representations can be classified if we know
enough about the algebra \(\CCINF(G)\).
For some groups, one can show that any irreducible representation that
has a central character is already strongly irreducible.  Given a
character~\(\chi\) of \(Z(G)\), we let \(I_\chi\defeq \{a\in Z(G)\mid
\chi(a)=0\}\) and
\[
\CCINF(G)_\chi \defeq \CCINF(G)/(I_\chi\cdot\CCINF(G)).
\]
This is the \emph{localization of~\(\CCINF(G)\) at the central
character~\(\chi\)}.  Evidently, the category of essential modules over
\(\CCINF(G)_\chi\) is isomorphic to the category of smooth
representations of~\(G\) with central character~\(\chi\).  Thus the
question whether or not irreducible representations of~\(G\) with
central character~\(\chi\) are strongly irreducible is a question about
the bornological algebra \(\CCINF(G)_\chi\).  I expect that for
reductive groups over local fields, the algebras \(\CCINF(G)_\chi\) are
always Morita equivalent to finite dimensional algebras.  This
immediately implies that any irreducible representation with central
character~\(\chi\) is strongly irreducible.

\section{The local trace formula}
\label{sec:local_trace_formula}

The results of this section are essentially already contained
in~\cite{Connes:Trace_Formula}.  Computations of a similar nature also appear
in the work of Jean-François Burnol (\cites{Burnol:Formules_explicites,
  Burnol:Fourier_Zeta}).  We have to reformulate Connes' computations because
we treat the global situation differently.  Let~\(\GF\) be a global field and
let~\(S\) be a finite set of places that contains all infinite places.  The
\emph{local trace formula} is an \(S\)\nbd{}local analogue of the Explicit
Formula.  Roughly speaking, it measures the failure of the Poisson Summation
Formula for the Fourier transform on~\(\Adel_S\).  Let \(\ICL_S\defeq
\Adel\inv_S/\GF\inv_S\) be the \(S\)\nbd{}idele class group.  The equivariant
Fourier transform~\(\EFou\) defines a unitary operator on the Hilbert space
\(L^2(\ICL_S,\abs{x}\,d\inv x)\).  Let~\(\phi\) be a smooth function on~\(\R_+\)
with \(\phi(x)=0\) near~\(\infty\) and \(\phi(x)=1\) near~\(0\).  Let~\(\lambda\) be the
regular representation of~\(\ICL_S\) on \(L^2(\ICL_S,\abs{x}\,d\inv x)\) and let
\[
T\defeq M_\phi - \EFou\circ M_\phi\circ\EFoi
\in \End(L^2(\ICL_S,\abs{x}\,d\inv x)).
\]
We show that \(\IN\lambda(f) T\) and \(T \IN\lambda(f)\) are nuclear and that
\[
\tr(T \IN\lambda(f)) =
\tr(\IN\lambda(f) T) = \sum_{v\in S} W_v(f)
\]
with the same local summands~\(W_v\) as in the Weil distribution.  For
\(S\to\Places(\GF)\), we get the Weil distribution except for the
discriminant term.  This comes in because of the different
normalization of the Fourier transform in the global case.

The local case can be treated by two rather disjoint sets of
techniques.  In this section, we use the Hilbert space approach and
work in the space \(L^2(\ICL_S,\abs{x}\,d\inv x)\).  In the next
section, we use the much smaller space \(\Sch(\Adel_S)/\GF\inv_S\)
instead.  It makes no difference which space we use in the local trace
formula.  That is, the problem of whether there are non-critical zeros
or not does not yet arise in the \(S\)\nbd{}local case.

\subsection{Some number theoretic conventions and notations}
\label{sec:local_trace_conventions}

We begin by explaining our setup and notation for the sake of the
non-experts.  Let~\(S\) be a \emph{finite} set of places of~\(\GF\) that
contains all infinite places.  Let
\[
\Adel_S \defeq \prod_{v\in S} \GF_v,
\qquad
\Adel\inv_S \defeq \prod_{v\in S} \GF\inv_v,
\]
be the \emph{\(S\)\nbd{}adele ring} and the \emph{\(S\)\nbd{}idele group}
of~\(\GF\), respectively.  We define the norm homomorphism
\(\Adel\inv_S\to\R\inv_+\), which is a continuous group homomorphism, by
\[
\abs{(x_v)_{v\in S}} =
\abs{(x_v)_{v\in S}}_S \defeq
\prod_{v\in S} \abs{x_v}_v.
\]
Its kernel is denoted \(\Adel\one_S\).  Addition and multiplication turn
\(\Adel_S\) and~\(\Adel\inv_S\) into locally compact Abelian groups.  Let
\(\Cont{S}=\Places(\GF)\setminus S\) and
\begin{align*}
  \GF_S &\defeq
  \{x\in\GF\mid
    \text{\(\abs{x}_v\le1\) for all \(v\in\Cont{S}\)}
  \},
  \\
  \GF\inv_S &\defeq
  \{x\in\GF\inv\mid
    \text{\(\abs{x}_v=1\) for all \(v\in\Cont{S}\)}
  \}.
\end{align*}
The set~\(\GF_S\) is a subring of~\(\GF\) and~\(\GF\inv_S\) is its
multiplicative group.  It is known that \(\GF_S\) and~\(\GF\inv_S\) are
discrete and cocompact subgroups of \(\Adel_S\) and~\(\Adel\one_S\),
respectively.  The \emph{\(S\)\nbd{}idele class group} \(\ICL_S\defeq
\Adel\inv_S/\GF\inv_S\) is another locally compact Abelian group.  The
norm homomorphism descends to a continuous homomorphism
\(\abs{\blank}\colon \ICL_S\to\R\inv_+\).  Its kernel
\(\ICL\one_S\subseteq\ICL_S\) is a compact group.

The groups \(\Adel\inv_S\) and~\(\ICL_S\) act on spaces of functions on
\(\Adel\inv_S\) and~\(\ICL_S\) by the \emph{regular representation}
\begin{equation}  \label{eq:def_lambda}
  \lambda_g f(x)\defeq f(g^{-1}\cdot x).
\end{equation}
The \emph{integrated form} of~\(\lambda\) is given by
\begin{equation}  \label{eq:integrate_lambda}
  \IN\lambda(h)(f)(x)\defeq \int h(g)\cdot f(g^{-1}\cdot x)\,d\inv g
\end{equation}
with respect to some fixed Haar measure \(d\inv g\) on \(\Adel\inv_S\)
or~\(\ICL_S\), respectively.  Notice that this is just the usual formula
for the convolution of two functions.

If~\(G\) is either~\(\ICL_S\) or~\(\GF\inv_v\) for some place~\(v\), we
normalize \(d\inv x\) so that
\[
\int_{\{x\in G\mid 1\le\abs{x}<\abs{u}\}} \,d\inv x = \ln{} \abs{u}
\qquad
\text{for all \(u\in G\) with \(\abs{u}\ge1\).}
\]
We do not care about the normalization of the Haar measure
on~\(\Adel\inv_S\).

If~\(v\) is an infinite place, the normalized Haar measure on~\(\GF_v\) is
the usual Lebesgue measure.  If~\(v\) is finite, the Haar measure is
normalized by the requirement that the maximal compact subring
\begin{equation}  \label{eq:def_MR}
  \MCR_v\defeq \{x\in\GF_v\mid \abs{x}_v\le 1\}
\end{equation}
should have unit volume.  The normalized Haar measure on~\(\Adel_S\) is
the product measure \(dx=\prod_{v\in S} dx_v\).  The measure
\(\abs{x}^{-1}\,dx\) is invariant under~\(\lambda\).  Hence it is
proportional to the Haar measure \(d\inv x\) on~\(\Adel\inv_S\).  For
this, the finiteness of~\(S\) is important: otherwise the set of
\(x\in\Adel_S\) with \(\abs{x}=0\) may have non-zero measure.

We identify the Pontrjagin dual~\(\widehat{\Adel_S}\) with \(\prod_{v\in
  S} \widehat{\GF_v}\).  Let \(\psi=(\psi_v)_{v\in
  S}\in\widehat{\Adel_S}\) be such that \(\psi_v\neq0\) for all \(v\in S\).
Then the map
\[
\Adel_S\to\widehat{\Adel_S},
\qquad
y\mapsto (x\mapsto\psi(xy))
\]
is an isomorphism of locally compact Abelian groups.  Hence we can
define the \emph{Fourier transform on~\(\Adel_S\)} by
\begin{equation}  \label{eq:def_Fourier}
  \Fourier=\Fourier_\psi\colon L^2(\Adel_S,dx)\to L^2(\Adel_S,dx),
  \qquad
  \Fourier f(\xi)\defeq
  \int_{\Adel_S} f(x) \psi(x\xi)\,dx.
\end{equation}
We call~\(\psi_v\) normalized if \(\Fourier_{\psi_v}\) is unitary on
\(L^2(\GF_v,dx)\).  Of course, this depends on our normalization
of~\(dx\).  We call~\((\psi_v)\) normalized if all~\(\psi_v\) are
normalized.  Hence \(\Fourier_\psi\) is a unitary operator on
\(L^2(\Adel_S,dx)\cong L^2(\Adel\inv_S,\abs{x}\,d\inv x)\).  Recall that
the measures \(dx\) and \(\abs{x}\,d\inv x\) are proportional.  If we do
not explicitly mention~\(\psi\), then we always assume that~\(\Fourier\)
is taken with respect to some normalized character~\(\psi\).  The
adjoint of~\(\Fourier\) is given simply by
\(\Fourier^{-1}=\Fourier^\ast=\lambda_{-1}\Fourier\), that is,
\(\Fourier^\ast f(x)=\Fourier f(-x)\).

The operator
\begin{equation}  \label{eq:def_J}
  J\colon L^2(\Adel\inv_S,\abs{x}\,d\inv x) \to
  L^2(\Adel\inv_S,\abs{x}\,d\inv x),
  \qquad
  Jf(x)\defeq \abs{x}^{-1} f(x^{-1}),
\end{equation}
is unitary and satisfies \(J^2=\ID\) and
\begin{equation}
  \label{eq:J_equivariance}
  J\circ \lambda_g = \abs{g}^{-1}\cdot \lambda_{g^{-1}} \circ J,
  \qquad
  J\circ \IN\lambda(h) = \IN\lambda(Jh)\circ J.
\end{equation}
Thus \(J(f_1\ast f_2)=Jf_1\ast Jf_2\), that is, \(J\) is a homomorphism
with respect to convolution.  The Fourier transform has the same lack
of equivariance, that is,
\begin{equation}
  \label{eq:Fourier_equivariance}
  \Fourier\circ \lambda_g =
  \abs{g}^{-1}\cdot \lambda_{g^{-1}} \circ \Fourier,
  \qquad
  \Fourier\circ \IN\lambda(h) = \IN\lambda(Jh)\circ \Fourier.
\end{equation}
Hence the compositions
\begin{equation}  \label{eq:def_EFou}
    \EFou \defeq \Fourier\circ J
    \quad\text{and}\quad
    \EFoi \defeq J\circ \Fourier^\ast,
\end{equation}
are \(\lambda\)\nbd{}equivariant unitary operators on \(L^2(\Adel_S,dx)\).
We call~\(\EFou\) the \emph{equivariant Fourier transform}.

\subsection{Descending operators to idele classes}
\label{sec:descend_L2_ICL}

Consider the Hilbert spaces
\[
L^2(\ICL_S)_\alpha \defeq
\{f\colon\ICL_S\to\C\mid
  f\cdot\abs{x}^\alpha\in L^2(\ICL_S,d\inv x)
\} =
L^2(\ICL_S,\abs{x}^{2\alpha}\,d\inv x)
\]
for \(\alpha\in\R\).  We are mainly interested in \(L^2(\ICL_S)_{1/2}\).
Clearly, \(J\) is a well-defined unitary operator on
\(L^2(\ICL_S)_{1/2}\).  We want to know that \(\Fourier\) or,
equivalently, \(\EFou\) also gives rise to a unitary operator on
\(L^2(\ICL_S)_{1/2}\).  This is not a very deep fact and can be proved
by elementary means.  Actually, there is a general method for
descending operators on \(L^2(\Adel\inv_S)_\alpha\) to
\(L^2(\ICL_S)_{\alpha}\).

If we work spectrally, we can identify \(L^2(\Adel\inv_S)\) with the
Hilbert space of square-integrable functions on the Pontrjagin dual
\(\widehat{\Adel\inv_S}\).  An \(\Adel\inv_S\)\brd{}equivariant bounded
operator~\(X\) on \(L^2(\Adel\inv_S)\) corresponds to the operator of
pointwise multiplication by the bounded measurable function~\(\hat{X}\)
on \(L^2(\widehat{\Adel\inv_S})\).  Since~\(\widehat{\ICL_S}\) is a closed
subgroup of \(\widehat{\Adel\inv_S}\), descending~\(X\) to~\(\ICL_S\)
corresponds to restricting~\(\hat{X}\) to
\(\widehat{\ICL_S}\subseteq\widehat{\Adel\inv_S}\).  However,
since~\(\widehat{\ICL_S}\) has vanishing measure, we cannot restrict
measurable functions to it.  We need~\(\widehat{X}\) to be continuous.
Equivalently, \(X\) is a multiplier of the reduced \(C^\ast\)\nbd{}algebra
\(\Cred(\Adel\inv_S)\) of~\(\Adel\inv_S\).  Let us denote multiplier
algebras by~\(\Mult\).  Restriction to
\(\widehat{\ICL_S}\subseteq\widehat{\Adel\inv_S}\) corresponds to a
unital \(\ast\)\nbd{}homomorphism
\[
\Mult(\Cred(\Adel\inv_S))\to\Mult(\Cred(\ICL_S)).
\]
This is how we descend operators on \(L^2(\Adel\inv_S)\) to
\(L^2(\ICL_S)\).

In order to verify whether \(X\in\Mult(\Cred(\Adel\inv_S))\), let
\[
\mathcal{R}\defeq \Cred(\Adel\inv_S) \cap L^2(\Adel\inv_S)
\subseteq L^2(\Adel\inv_S).
\]
Then \(X\in\Mult(\Cred(\Adel\inv_S))\) if and only if
\(X(\mathcal{R})\subseteq\mathcal{R}\).  Actually, it suffices to check
that \(X(\CCINF(\Adel\inv_S))\subseteq\mathcal{R}\) because for any
\(\chi\in\widehat{\Adel\inv_S}\) there is \(f\in\CCINF(\Adel\inv_S)\) with
\(\hat{f}(\chi)\neq0\).  Hence continuity of \(\hat{X}\cdot\hat{f}\)
implies that~\(\hat{X}\) is continuous near~\(\chi\).

Since continuity of the Fourier transform is difficult to check, the
space~\(\mathcal{R}\) is difficult to describe.  However, it is clear
that
\[
L^1(\Adel\inv_S,d\inv x)\cap L^2(\Adel\inv_S,d\inv x)\subseteq
\mathcal{R}.
\]
This allows us to show that~\(\EFou\) belongs to
\(\Mult(\Cred(\Adel\inv_S))\):

\begin{lemma}  \label{lem:Fourier_unitary_ICL}
  The operators \(\EFou\) and~\(\EFoi\) descend to unitary operators on
  \(L^2(\ICL_S)_{1/2}\) that are inverse to each other and belong to the
  multiplier algebra of \(\Cred(\ICL_S)_{1/2}\).
\end{lemma}

\begin{proof}
  Since we work in \(L^2(\Adel\inv_S,\abs{x}\,d\inv x)\), we have to
  modify the definition of~\(\mathcal{R}\) accordingly.  It is easy to
  see that
  \[
  \EFou(\CCINF(\Adel\inv_S))\subseteq
  \Sch(\Adel_S)\subseteq
  L^1(\Adel\inv_S,\abs{x}^{1/2}\,d\inv x) \cap
  L^2(\Adel\inv_S,\abs{x}\,d\inv x)
  \subseteq
  \mathcal{R}_{1/2}.
  \]
  This implies \(\EFou\in\Mult(\Cred(\Adel\inv_S))\).  The
  assertions of the lemma now follow from the above discussion.
\end{proof}

Equation~\eqref{eq:def_J} also defines~\(J\) as an operator on
\(L^2(\ICL_S)_{1/2}\).  Hence \(\Fourier=\EFou J\) also descends to
\(L^2(\ICL_S)_{1/2}\).  Since \(-1\in\GF\inv_S\), we have
\(\Fourier=\Fourier^\ast\) as operators on \(L^2(\ICL_S)_{1/2}\), so that
\(\Fourier^2=\ID\).  The above constructions depend on spectral analysis
and therefore require the commutativity of~\(\ICL_S\).  A generalization
to crossed products by locally compact groups can be found
in~\cite{Meyer:Fixed}.

\subsection{The standard Fredholm module over \(\Sch(\ICL_S)\)}
\label{sec:Fredholm_ICL_S}

We let \(\Sch(\ICL_S)\) be the Bruhat-Schwartz algebra of the locally
compact Abelian group~\(\ICL_S\).  It is defined
in~\cite{Bruhat:Distributions} and discussed in greater detail below.
We describe a standard Fredholm module over \(\Sch(\ICL_S)\).  A part of
this construction is used for the local trace formula.

Multiplication by the unitary character \(\abs{x}^{\ima t}\) for
\(t\in\R\) is an automorphism of~\(\Sch(\ICL_S)\).  This defines a smooth
action of~\(\R\) by bornological algebra automorphisms on
\(\Sch(\ICL_S)\).  Its generator is the bounded derivation
\[
\partial = M_{\ln{}\abs{g}}\colon \Sch(\ICL_S)\to\Sch(\ICL_S),
\qquad
\partial(f)(g) \defeq \ln {}\abs{g}\cdot f(g).
\]
The derivation property means that
\[
\partial(f_0\ast f_1)=\partial(f_0)\ast f_1+ f_0\ast\partial(f_1)
\]
for all \(f_0,f_1\in\Sch(\ICL_S)\).  The \emph{standard trace} on
\(\Sch(\ICL_S)\) is the bounded linear functional
\[
\tau\colon \Sch(\ICL_S) \to \C,
\qquad
\tau(f)=f(1_{\ICL_S}).
\]
It satisfies \(\tau\circ\partial=0\).  We abbreviate
\(L^2(\ICL_S)\defeq L^2(\ICL_S,d\inv x)\).  Since~\(\lambda\) is a unitary
representation on \(L^2(\ICL_S)\), its integrated form extends to a
bounded \(\ast\)\nbd{}homomorphism \(\IN\lambda\colon
\Sch(\ICL_S)\to\End(L^2(\ICL_S))\).

Let \(\phi\colon \R_+\to\R_+\) be a smooth function such that
\(\phi(x)=1\) for~\(x\) close to~\(0\) and \(\phi(x)=0\) for large enough~\(x\).
We pull back~\(\phi\) to a smooth function on~\(\ICL_S\) by
\(\phi(g)\defeq\phi(\abs{g})\).  The operator~\(M_\phi\) of multiplication
by~\(\phi\) is bounded on \(L^2(\ICL_S)\).  In this section, the
smoothness of~\(\phi\) is irrelevant.  Things actually become simpler if
we choose~\(\phi\) to be the characteristic function of \(\ocival{0,1}\).
Later we want~\(M_\phi\) to operate on spaces of smooth functions and
therefore require smooth~\(\phi\).

The following notation is convenient.  Suppose that we have a linear
map \(h\colon V\to\Hom(X,Y)\).  We say that \(h(v)\) is \emph{uniformly
  nuclear} for \(v\in V\) if \(h(v)\) is nuclear for all \(v\in S\) and if
\(h(S)\) is bounded in the space of nuclear operators for all bounded
subsets \(S\subseteq V\).  For instance, a representation~\(\pi\) of a
locally compact group~\(G\) is summable if and only if \(\IN\pi(f)\) is
uniformly nuclear for \(f\in\CCINF(G)\).

\begin{lemma}  \label{lem:trace_trivial}
  Let~\(\phi\) be as above and let \(h\in\Sch(\ICL_S)\).  Then the
  operators
  \[
  [M_\phi,\IN\lambda(f)],
  \qquad
  \IN\lambda(f) M_h,
  \qquad
  M_h\IN\lambda(f)
  \]
  on \(L^2(\ICL_S)\) are uniformly nuclear for \(f\in\Sch(\ICL_S)\).
  We have
  \[
  \tr \IN\lambda(f_0) [M_\phi,\IN\lambda(f_1)] =
  \tau(-f_0 \ast \partial f_1) =
  \tau(\partial(f_0)\ast f_1)
  \]
  and \(\tr \IN\lambda(f_0)M_h = f_0(1)\cdot \int_{\ICL_S} h(x)\,d\inv
  x\) for all \(f_0,f_1\in \Sch(\ICL_S)\).
\end{lemma}

\begin{proof}
  We have
  \[
  [M_\phi,\IN\lambda(f)] h(x) =
  \int_{\ICL_S} \bigl(\phi(\abs{x})-\phi(\abs{y})\bigr) f(xy^{-1})
  \cdot h(y)\,d\inv y,
  \]
  so that the operator \([M_\phi,\IN\lambda(f)]\) has the integral
  kernel
  \[
  K_f(x,y)\defeq f(xy^{-1})\bigl(\phi(\abs{x})-\phi(\abs{y})\bigr).
  \]
  We claim that \(f\mapsto K_f\) is a bounded map from \(\Sch(\ICL_S)\)
  to \(\Sch(\ICL_S) \hot \Sch(\ICL_S) \cong \Sch(\ICL_S\times\ICL_S)\).
  This yields the uniform nuclearity of \([M_\phi,\IN\lambda(f)]\).  To
  prove the claim, observe first that~\(f\) belongs to \(\Sch(\ICL_S/k)\)
  for some elementary Abelian quotient \(\ICL_S/k\) of~\(\ICL_S\).  Hence
  \(K_f\) is a function on \((\ICL_S/k)^2\).  Thus we may replace~\(\ICL_S\)
  by the elementary Abelian group \(G\defeq \ICL_S/k\).  It is important
  that the function \(\abs{\ln{}\abs{x}}\) is a proper length function
  and hence defines rapid decay on~\(G\).  There is \(R\in\R_{\ge1}\)
  depending only on~\(\phi\) such that \(K_f(x,y)\) vanishes for
  \(\abs{x},\abs{y}\ge R\) and \(\abs{x},\abs{y}\le R^{-1}\).  This
  implies
  \[
  \max \{\abs{x},\abs{x}^{-1},\abs{y},\abs{y}^{-1}\} \le
  R\cdot \max\{\abs{xy^{-1}}, \abs{xy^{-1}}^{-1}\}
  \]
  on the support of~\(K_f\).  Hence the rapid decay of~\(f\) implies
  rapid decay of~\(K_f\).  If we apply a \(G^2\)\nbd{}equivariant
  differential operator to~\(K_f\), we do not enlarge the support and we
  retain rapid decay in the variable \(xy^{-1}\).  Hence the
  differentiated functions still have rapid decay.  This yields the
  claim.
  
  Similarly, the operator \(\IN\lambda(f_0)[M_\phi,\IN\lambda(f_1)]\)
  has the integral kernel
  \[
  K(x,y) \defeq \int_{\ICL_S} f_0(xz^{-1})(\phi(\abs{z})-\phi(\abs{y}))
  f_1(zy^{-1}) \,d\inv z,
  \]
  which belongs to \(\Sch(\ICL_S\times\ICL_S)\) as well.  The trace
  of the operator defined by such an integral kernel is equal to the
  integral of the kernel along the diagonal.  Thus
  \begin{multline*}
    \tr \IN\lambda(f_0) [M_\phi,\IN\lambda(f_1)]
    =
    \int_{\ICL_S} K(x,x)\,d\inv x
    \\ =
    \int_{\ICL_S} f_0(a) f_1(a^{-1})
    \int_{\ICL_S} \phi(\abs{z})-\phi(\abs{az}) \,d\inv z\,d\inv a.
  \end{multline*}
  We compute \(\int_{\ICL_S} \phi(\abs{z})-\phi(\abs{a z})\,d\inv z\).
  If~\(\phi\) had compact support, the left invariance of \(d\inv z\)
  would force the integral to vanish.  Therefore, we may
  replace~\(\phi\) by any function~\(\phi'\) with the same behavior at \(0\)
  and~\(\infty\).  We choose~\(\phi'\) to be the characteristic function
  of \(\ocival{0,1}\).  Thus \(\phi'(\abs{z})-\phi'(\abs{a z})\) is the
  characteristic function of the set \(\{z\in\ICL_S\mid \abs{a}^{-1}<
  \abs{z}\le 1\}\) for \(\abs{a}\ge1\).  Hence the integral is
  \(\ln{}\abs{a}\).  We also get \(\ln{}\abs{a}\) for \(\abs{a}\le1\).  This
  yields the asserted value for \(\tr \IN\lambda(f_0)
  [M_\phi,\IN\lambda(f_1)]\).
  
  The corresponding results for \(\IN\lambda(f_0)M_h\) are simpler and
  proved similarly\mdash and therefore left to the reader.
\end{proof}

Let \(F=1-2M_\phi\).  Then~\(\IN\lambda\) and~\(F\) define a
\emph{\(1\)\nbd{}summable odd Fredholm module} over \(\Sch(\ICL_S)\).
That is, \(\IN\lambda(f)\cdot (1-F^2)\) and \([\IN\lambda(f),F]\) are
uniformly nuclear operators on \(L^2(\ICL_S)\) for \(f\in \Sch(\ICL_S)\)
and \(F=F^\ast\).  If we let~\(\phi\) be the characteristic function of
\(\ocival{0,1}\), then even \(F^2=\ID\).  Lemma~\ref{lem:trace_trivial}
implies easily that the Chern-Connes character of this Fredholm module
is the cyclic \(1\)\nbd{}cocycle
\[
\chi(f_0,f_1)\defeq \tau(f_0 \ast \partial f_1).
\]
However, we do not use this fact.  Instead, we twist this Fredholm
module by~\(\EFou\) and obtain the \(0\)\nbd{}cocycle \(\chi(\EFou,\EFoi
f)\) in Lemma~\ref{lem:local_trace_I}.

\subsection{A commutator with the equivariant Fourier transform}
\label{sec:local_trace_details}

Since we want to use the equivariant Fourier transform, we replace
\(L^2(\ICL_S)\) by \(L^2(\ICL_S)_{1/2}\), on which the operators \(\EFou\)
and~\(\EFoi\) are defined.  These Hilbert spaces are isomorphic via
multiplication by \(\abs{x}^{1/2}\).  The algebra \(\Sch(\ICL_S)\) now
acts by \(\IN\lambda(f\cdot\abs{x}^{-1/2})\).  Since we do not want to
write \(\abs{x}^{-1/2}\) in all our formulas, we replace \(\Sch(\ICL_S)\)
by the \(\ast\)\nbd{}algebra
\[
\Sch(\ICL_S)_{1/2}\defeq
\Sch(\ICL_S)\cdot\abs{x}^{-1/2} =
\{f\colon \ICL_S\to\C\mid f\cdot \abs{x}^{1/2}\in \Sch(\ICL_S)\},
\]
which is isomorphic to \(\Sch(\ICL_S)\) via multiplication by
\(\abs{x}^{\pm 1/2}\).  For~\(\phi\) as above we let
\[
T\defeq M_\phi - \EFou\circ M_\phi\circ\EFoi
\colon L^2(\ICL_S)_{1/2}\to L^2(\ICL_S)_{1/2}.
\]
We prefer to choose~\(\phi\) smooth, but again it is also possible to
take the characteristic function of~\(\ocival{0,1}\).

\begin{lemma}  \label{lem:local_trace_I}
  The operators \(\IN\lambda(f) T\) and \(T \IN\lambda(f)\) on
  \(L^2(\ICL_S)_{1/2}\) are uniformly nuclear for
  \(f\in\Sch(\ICL_S)_{1/2}\) and
  \begin{equation}
    \label{eq:local_trace}
    \tr(\IN\lambda(f) T) =
    \tr(T \IN\lambda(f)) =
    \tau(\EFou \partial(\EFoi f))
  \end{equation}
  for all \(f\in\Sch(\ICL_S)_{1/2}\).
\end{lemma}

\begin{proof}
  We have \([\IN\lambda(f),T] = [\IN\lambda(f),M_\phi] - \EFou
  [\IN\lambda(f),M_\phi] \EFoi\) because \(\EFou\) and~\(\EFoi\) are
  equivariant.  Lemma~\ref{lem:trace_trivial} implies that
  \([\IN\lambda(f),M_\phi]\) and hence \([\IN\lambda(f),T]\) are uniformly
  nuclear for \(f\in\Sch(\ICL_S)_{1/2}\).  Since \(\EFou\) and \(\EFoi\) are
  bounded on \(L^2(\ICL_S)\) and inverse to each other and
  \([\IN\lambda(f),M_\phi]\) is nuclear, we get
  \begin{multline*}
    \tr {}[\IN\lambda(f),T]
    =
    \tr {}[\IN\lambda(f),M_\phi] -
    \tr \EFou [\IN\lambda(f),M_\phi] \EFoi
    \\ =
    \tr {}[\IN\lambda(f),M_\phi] -
    \tr {}[\IN\lambda(f),M_\phi] \EFoi \EFou
    =
    0.
  \end{multline*}
  Thus it suffices to check uniform nuclearity and to compute the
  trace for \(\IN\lambda(f) T\).  We rewrite this operator using the
  equivariance of \(\EFou\) and~\(\EFoi\).  We have \(\EFou\circ
  \IN\lambda(f)(h)=\EFou(f\ast h)=\EFou(f)\ast h = \IN\lambda(\EFou
  f)(h)\) for all \(f,h\in\Sch(\ICL_S)_{1/2}\).  That is,
  \begin{equation}  \label{eq:Fourier_lambda}
    \IN\lambda(f)\circ\EFou =
    \EFou\circ\IN\lambda(f) =
    \IN\lambda(\EFou f)
  \end{equation}
  on \(\Sch(\ICL_S)_{1/2}\).  Lemma~\ref{lem:Fourier_unitary_ICL} yields
  that \(\EFou f\in\Cred(\ICL_S)_{1/2}\).  Hence both sides
  of~\eqref{eq:Fourier_lambda} define the same bounded linear operator
  on \(L^2(\ICL_S)_{1/2}\).

  We compute
  \begin{displaymath}
    \IN\lambda(f) T
    =
    \IN\lambda(f)\circ M_\phi -
    \IN\lambda(f)\circ \EFou M_\phi \EFoi
    =
    -[M_\phi,\IN\lambda(f)] + [M_\phi, \IN\lambda(\EFou f)]\circ\EFoi.
  \end{displaymath}
  A finer analysis than the above
  (Proposition~\ref{pro:Fourier_estimate}) yields \(\EFou
  f\in\Sch(\ICL_S)_{1/2}\).  Since~\(\EFoi\) is a bounded operator on
  \(L^2(\ICL_S)_{1/2}\), Lemma~\ref{lem:trace_trivial} implies that
  \(\IN\lambda(f)\circ T\) is uniformly nuclear for
  \(f\in\Sch(\ICL_S)_{1/2}\).
  
  Since convolutions \(f_1\ast f_2\) are dense in \(\Sch(\ICL_S)_{1/2}\), it
  suffices to compute the trace of \(\IN\lambda(f_1\ast f_2) T\).  We have
  \begin{multline*}
    \tr \IN\lambda(f_1\ast f_2) T
    =
    - \tr \IN\lambda(f_1) [M_\phi,\IN\lambda(f_2)]
    + \tr \IN\lambda(f_1) [M_\phi, \IN\lambda(\EFou f_2)]\circ\EFoi
    \\ =
    \tau(f_1\ast \partial f_2) +
    \tr \EFoi \IN\lambda(f_1) [M_\phi, \IN\lambda(\EFou f_2)]
    \\ =
    \tau(f_1\ast \partial f_2)
    + \tau(\partial(\EFoi f_1)\ast \EFou f_2)
    =
    \tau(\EFou \partial(\EFoi f_1\ast f_2))
  \end{multline*}
  because \(\EFou\) and~\(\EFoi\) are inverse to each other.  This
  yields~\eqref{eq:local_trace}.
\end{proof}

We are going to compute \(\tau(\EFou \partial(\EFoi f))\).  Let
\(H\defeq\GF\inv_S\).  For any \(f\in\CCINF(\ICL_S)\), there is
\(\bar{f}\in\CCINF(\Adel\inv_S)\subseteq \Sch(\Adel_S)\) with
\(\sum_{x\in H} \bar{f}(gx)=f(g)\) for all \(g\in\ICL_S\).  Since
\(\int_{\Adel_S} h(x)\ln{}\abs{x}\,dx\) exists for all
\(h\in\Sch(\Adel_S)\), the Fourier transform of \(\ln{}\abs{x}\) is a
well-defined tempered distribution on~\(\Adel_S\).  Let~\(\addconv\)
denote the convolution on the additive group~\(\Adel_S\).  Since the
Fourier transform intertwines pointwise multiplication and convolution
and \(\EFou\circ\partial\circ\EFoi=-\Fourier
M_{\ln{}\abs{x}}\Fourier^\ast\), we get
\begin{displaymath}
  \tau(\EFou \partial(\EFoi f)) =
  -\sum_{a\in H} (\Fourier M_{\ln {}\abs{x}} \Fourier^{-1} \bar{f})(a)
  =
  -\sum_{a\in H} \bigl((\Fourier \ln {}\abs{x}) \addconv
  \bar{f}\bigr)(a),
\end{displaymath}
We have \(\ln{}\abs{x}_S = \sum_{v\in S} \ln{}\abs{x}_v\).  This implies
\[
\Fourier(\ln{}\abs{x}_S)=
\sum_{v\in S} \pi_v^\ast \Fourier(\ln{}\abs{x}_v),
\]
where \(\braket{\pi_v^\ast D}{f}\defeq D(f|_{0\times\GF_v})\) for any
\(f\in\Sch(\Adel_S)\), \(D\in\Sch'(\GF_v)\).  Thus convolution with
\(\pi_v^\ast D\) only moves in the direction of~\(\GF_v\) and leaves the
other coordinates fixed.

\begin{lemma}[\cite{Connes:Trace_Formula}*{page~50}]
  \label{lem:Fourier_ln_norm}
  The Fourier transform \(P_v\defeq \Fourier(\ln{}\abs{x}_v)\) satisfies
  \begin{equation}  \label{eq:Fourier_ln_norm}
    \braket{P_v}{\lambda_a f}= \braket{P_v-\ln{}\abs{a} \delta_0}{f}
  \end{equation}
  for all \(f\in\Sch(\GF_v)\).  We have
  \[
  \braket{P_v}{f} = -\int_{\GF\inv_v} f(x) \,d\inv x
  \]
  for all \(f\in\Sch(\GF_v)\) with \(f(0)=0\).  Thus~\(P_v\) is a principal
  value for the above integral.  It is determined uniquely by the
  normalization condition \(\Fourier P_v(1)=0\).
\end{lemma}

\begin{proof}
  The relation \(\ln{}\abs{a^{-1}x} = \ln{}\abs{x}-\ln{}\abs{a}\) for
  \(a\in\GF\inv_v\) easily implies~\eqref{eq:Fourier_ln_norm}.  If we
  restrict to functions that vanish at~\(0\), the multiple of~\(\delta_0\)
  vanishes, so that~\(P_v\) is \(\lambda\)\nbd{}invariant on this
  subspace.  By the uniqueness of Haar measure, we get \(P_v(f) =
  c\cdot \int_{\GF\inv_v} f(x)\,d\inv x\) for all
  \(f\in\CCINF(\GF\inv_v)\) with some constant~\(c\).  This extends to all
  \(f\in\Sch(\GF_v)\) with \(f(0)=0\).  This can be proved directly.  It
  also follows rather easily from Lemma~\ref{lem:Sch_local} and a
  spectral description of the space \(\Twist(\GF\inv_v)\) as in
  Section~\ref{sec:global_spectral_analysis}.
  
  If \(f\in\CCINF(\GF_v)\) is constant in a neighborhood of~\(0\), then
  \(f-\lambda_a f\in\CCINF(\GF\inv_v)\) for all \(a\in\GF\inv_v\).
  Equation~\eqref{eq:Fourier_ln_norm} yields
  \[
  c\cdot \int_{\GF\inv_v} f(x) - f(a^{-1}x) \,d\inv x =
  \braket{P_v}{f-\lambda_a f} =
  \ln{}\abs{a} f(0).
  \]
  Computing the integral on the left hand side as in the proof of
  Lemma~\ref{lem:trace_trivial}, we obtain \(c=-1\).  Thus~\(P_v\) is a
  principal value for \(-\int_{\GF\inv_v} f(x)\,d\inv x\).  So far we
  have determined~\(P_v\) up to the addition of a constant multiple of
  the Dirac measure~\(\delta_0\).  This is fixed by the additional
  condition \(\Fourier P_v(1)=\ln{}\abs{1}=0\).
\end{proof}

We get
\begin{multline*}
  ((\Fourier \ln{}\abs{x}_v) \addconv \bar{f}\bigr)(a)
  =
  \braket{P_v}{x\mapsto \bar{f}(a-x)}
  \\ =
  \braket{P_v}{x\mapsto \bar{f}(a\cdot (1-a_v^{-1}x))}
  =
  \braket{P_v}{x\mapsto \bar{f}(a\cdot (1-x))}
  - \ln{}\abs{a}_v \bar{f}(1).
\end{multline*}
The terms \(\ln{}\abs{a}_v\) are canceled by the summation over~\(S\)
because \(\sum_{v\in S} \ln{}\abs{a}_v = \ln{}\abs{a}_S=0\) for
\(a\in H\).  The summation over~\(H\) replaces~\(\bar{f}\) by~\(f\) again.
Therefore,
\begin{multline*}
  - \sum_{v\in S} \sum_{a\in H}
  (\Fourier \ln{}\abs{x}_v) \addconv \bar{f}(a)
  =
  - \sum_{v\in S} \braket{P_v}{x\mapsto f(1-x)}
  \\ =
  \sum_{v\in S} \int'_{\GF\inv_v} f(1-x) \,d\inv x
  =
  \sum_{v\in S} \int'_{\GF\inv_v} \frac{f(x)\abs{x}}{\abs{1-x}}
  \,d\inv x,
\end{multline*}
where we use \(\abs{1-x}\,d\inv(1-x)= \abs{x}\,d\inv x\) and pull back
\(f\in\Sch(\ICL_S)_{1/2}\) to a function on~\(\GF\inv_v\) using the
evident homomorphism \(\GF\inv_v\to\ICL_S\).  Of course, the last
principal value differs from the first one because we have shifted the
singularity of the integral from~\(0\) to~\(1\).  We still just get the
ordinary integral if \(f(1)=0\).  The normalization of the principal
value at~\(v\) only depends on the local field~\(\GF_v\).

\begin{lemma}  \label{lem:normalize_principal_value}
  For a finite local field~\(\LF\), the principal value is normalized by
  \[
  \int'_{\MCR\inv_\LF} \frac{\abs{x}}{\abs{1-x}} \,d\inv x=0.
  \]
\end{lemma}

\begin{proof}
  Let \(\MI\subseteq\MCR\) be the maximal ideal of elements with
  \(\abs{x}<1\).  Reversing the above computation, we obtain
  \begin{multline*}
    \int'_{\MCR\inv_\LF} \frac{\abs{x}}{\abs{1-x}} \,d\inv x =
    \int'_{\LF\inv} \frac{1_{\MCR\inv}\abs{x}}{\abs{1-x}} \,d\inv x =
    \int'_{\LF\inv} 1_{1-\MCR\inv} \,d\inv x
    \\ =
    \int'_{\LF\inv} 1_{\MCR} - 1_{1+\MI} \,d\inv x =
    - \braket{P_\LF}{1_{\MCR}} -
    \ln(q_\LF) (q_\LF-1)^{-1}
  \end{multline*}
  because \(1+\MI\subseteq\LF\inv\) has volume \(\ln(q_\LF)
  (q_\LF-1)^{-1}\) for our normalization of the Haar measure.  Since
  the Fourier transform of~\(1_{\MCR}\) is again~\(1_{\MCR}\), we also get
  \begin{displaymath}
    \braket{P_\LF}{1_{\MCR}}
    =
    \int_{\{x\in\LF\mid \abs{x}\le 1\}} \ln{}\abs{x} \,dx
    =
    -\ln(q_\LF) (q_\LF-1)^{-1}.
  \end{displaymath}
  This yields the assertion.
\end{proof}

The normalization of the principal values at the infinite places is
more painful to describe explicitly.  Alain Connes checks
in~\cite{Connes:Trace_Formula}*{page 91\ndash 96} that the above
principal values are the same as those used by André Weil in his
Explicit Formula (\cites{Weil:Explicit_Formula,
Weil:Explicit_Formula_first}).  Putting everything together, we obtain
\begin{equation}  \label{eq:local_trace_formula}
  \tr(\IN\lambda(f) T) =
  \sum_{v\in S} \int'_{\GF\inv_v} \frac{f(x)\abs{x}}{\abs{1-x}}
  \,d\inv x
\end{equation}
for all \(f\in\Sch(\ICL_S)\) with the same principal values as
in~\cites{Weil:Explicit_Formula, Weil:Explicit_Formula_first}.  We
call~\eqref{eq:local_trace_formula} the \emph{local trace formula}.

What happens if we take the Fourier transform with respect to a
non-normalized character~\(\psi\)?  Of course, then we should also
renormalize the Haar measure so that~\(\Fourier\) is again unitary on
\(L^2(\Adel_S,dx)\).  We have \(\psi(x)=\psi_0(\different\cdot
x)\) for some \(\different\in\Adel\inv_S\) and a normalized
character~\(\psi_0\).  The Fourier transforms for \(\psi\) and~\(\psi_0\)
differ by
\[
\Fourier_\psi f =
\abs{\different}_S^{-1/2} \Fourier(\lambda_\different f).
\]
Therefore,
\[
\braket{\Fourier_\psi \ln{}\abs{x}_S}{f}=
\braket{\Fourier \ln{}\abs{\different^{-1}x}_S}{f} =
\braket{\Fourier \ln{}\abs{x}_S}{f}
+\ln{}\abs{\different}_S^{-1} \cdot f(0).
\]
Since convolution with the Dirac measure~\(\delta_0\) is the identity
map, we obtain
\begin{equation}  \label{eq:local_trace_formula_psi}
  \tr(
\IN\lambda(f) T_\psi) =
  \sum_{v\in S} \int'_{\GF\inv_v} \frac{f(x)\abs{x}}{\abs{1-x}}
  \,d\inv x
  - f(1)\cdot \ln{}\abs{\different}^{-1}_S.
\end{equation}

This variant of~\eqref{eq:local_trace_formula} becomes important when
we pass to the global situation.

\section{The local trace formula without Hilbert spaces}
\label{sec:local_trace_formula_Sch}

Again we let~\(\GF\) be a global field and we let~\(S\) be a finite set of
places of~\(\GF\) containing all infinite places.  In order to pass to
the global case, we have to replace \(L^2(\ICL_S)_{1/2}\) by a smaller
space that allows us to work outside the critical strip.  This is the
purpose of the space \(\Twist(\ICL_S) \boin L^2(\ICL_S)_{1/2}\).  The
notation~\(\boin\) means that \(\Twist(\ICL_S)\subseteq
L^2(\ICL_S)_{1/2}\) and that the inclusion map is bounded.  The space
\(\Twist(\ICL_S)\) consists of those functions \(f\colon \ICL_S\to\C\) for
which \(f\cdot\abs{x}^\alpha\) is a Bruhat-Schwartz function for
\(\alpha>0\) and \((\EFoi f) \cdot\abs{x}^\alpha\) is a Bruhat-Schwartz
function for \(\alpha<1\).  The operator~\(T\) of the previous section
restricts to a bounded operator on \(\Twist(\ICL_S)\) and the operators
\(\IN\lambda(f)T\) and \(T\IN\lambda(f)\) are uniformly nuclear
operators \(L^2(\ICL_S)_{1/2}\to\Twist(\ICL_S)\) for
\(f\in\CCINF(\ICL_S)\).  Therefore, they are uniformly nuclear as
operators on \(L^2(\ICL_S)_{1/2}\) and \(\Twist(\ICL_S)\) and have the
same trace on both spaces.

If~\(\GF\) is a global function field, then there is a variant
\(\Twist_\comp(\ICL_S)\) of \(\Twist(\ICL_S)\) that should be used
instead.  With \(\Twist(\ICL_S)\) we can prove the meromorphic extension
of \(L\)\nbd{}functions to \(\Irrep(\ICL_S)\).  With
\(\Twist_\comp(\ICL_S)\) we even get that they are rational functions.

The definition of the spaces \(\Twist(\ICL_S)\) and
\(\Twist_\comp(\ICL_S)\) is rather \emph{ad hoc}.  There is a more
natural way to get them.  The coinvariant space
\(\Sch(\Adel_S)/\GF\inv_S\) is, by definition, the quotient of
\(\Sch(\Adel_S)\) by the closed linear span of elements of the form
\(\lambda_a f-f\) with \(a\in\GF\inv_S\), \(f\in\Sch(\Adel_S)\).  We show
that \(\Sch(\Adel_S)/\GF\inv_S\) is isomorphic to \(\Twist(\ICL_S)\) for
algebraic number fields and to \(\Twist_\comp(\ICL_S)\) for global
function fields.  We also prove that the higher group homology
\(H_n(\GF\inv_S,\Sch(\Adel_S))\) for \(n\ge1\) vanishes.  This is
necessary in order to compute the Hochschild and cyclic homology of
the crossed product \(\GF\inv\cross\Sch(\Adel_\GF)\).

We remark that very little of this is needed in order to construct the
global difference representation and compute its character.  We may
just take the space \(\Twist(\ICL_S)\) or \(\Twist_\comp(\ICL_S)\) and be
happy with it.  The only additional fact that we need is the ``maximum
principle'' of Lemma~\ref{lem:maximum_principle}.  We mainly bring
\(\Sch(\Adel_S)/\GF\inv_S\) into play because we want the global
difference representation to be as natural as possible.

In order to compute \(H_n(\GF\inv_S,\Sch(\Adel_S))\), we proceed in two
steps.  We only explain this in the case where~\(\GF\) is an algebraic
number field.  Define \(\Twist(\Adel\inv_S)\) in the same way as
\(\Twist(\ICL_S)\).  Let \(\Lhot_{\CCINF(\Adel\one_S)}\) denote the total
left derived functor of the balanced tensor product functor
\(\hot_{\CCINF(\Adel\one_S)}\).  Let \(\Sch(\Adel\one_S)\) be the
Bruhat-Schwartz algebra of~\(\Adel\one_S\).  We show in the first step
that
\[
\Sch(\Adel\one_S)\Lhot_{\CCINF(\Adel\one_S)} \Sch(\Adel_S) =
\Twist(\Adel\inv_S).
\]
Roughly speaking, this means that from the point of view of those
representations of~\(\Adel\inv_S\) whose restriction to~\(\Adel\one_S\) is
tempered, there is no difference between \(\Sch(\Adel_S)\) and
\(\Twist(\Adel\inv_S)\).  I expect an analogue of this result in the
automorphic case, but I have not yet checked the details.  In the
second step we use that the canonical representation of~\(\Adel\inv_S\)
on \(\CCINF(\ICL_S)\) integrates to a (right) module structure over
\(\Sch(\Adel\one_S)\) because~\(\ICL\one_S\) is compact.  An easy
computation then shows
\begin{displaymath}
  \C(1) \Lhot_{\CCINF(\GF\inv_S)} \Sch(\Adel_S)
  \cong
  \Twist(\ICL_S).
\end{displaymath}
This yields the assertions about group homology and the coinvariant
space as the homology spaces of the chain complex \(\C(1)
\Lhot_{\CCINF(\GF\inv_S)} \Sch(\Adel_S)\).  Our computation requires
some general facts about homology for smooth convolution algebras.  We
need that \(V\Lhot_{A(G)} W \cong V\Lhot_{\CCINF(G)} W\) if \(V\) and~\(W\)
are right modules over \(A(G)\), where \(G=\Adel\inv_S\) and \(A(G)\) is a
certain weighted Schwartz algebra.  We call \(A(G)\)
\emph{isocohomological} if this happens.

\subsection{Weighted Bruhat-Schwartz algebras on locally compact
  Abelian groups}
\label{sec:Sch_weighted}

First we recall some facts about the Bruhat-Schwartz algebra \(\Sch(G)\)
of a locally compact Abelian group~\(G\)
(see~\cite{Bruhat:Distributions}).  We need the full generality of
Bruhat's definition for the adelic groups we are dealing with.  If~\(G\)
is elementary Abelian, \(\Sch(G)\) is defined as usual.  It is a nuclear
Fréchet space.  We equip it with the von Neumann bornology to get a
bornological vector space.  For arbitrary~\(G\), we consider pairs of
subgroups \(k\subseteq U\subseteq G\) such that~\(k\) is compact and \(U/k\)
is elementary Abelian.  We call \(U/k\) a \emph{proper elementary
  Abelian subquotient of~\(G\)}.  We can pull back \(\Sch(U/k)\) to a
space of functions on~\(G\).  The proper elementary Abelian subquotients
of~\(G\) form a directed set and the spaces \(\Sch(U/k)\) form an
inductive system with injective structure maps.  We let \(\Sch(G)\) be
its bornological direct union.  That is, a subset of \(\Sch(G)\) is
bounded if and only if it is bounded in \(\Sch(U/k)\) for some proper
elementary Abelian subquotient \(U/k\).  Since \(\Sch(G)\) for
arbitrary~\(G\) looks locally (that is, on any bounded subset) like
\(\Sch(G)\) for elementary Abelian~\(G\), it inherits many properties from
the elementary Abelian case.

The space \(\Sch(G)\) is a nuclear, complete, convex bornological vector
space for all~\(G\).  Both pointwise multiplication and convolution are
bounded bilinear maps on \(\Sch(G)\).  The group~\(G\) itself and its
Pontrjagin dual~\(\hat{G}\) act on \(\Sch(G)\) by the regular
representation~\(\lambda\) defined as in~\eqref{eq:def_lambda} and by
pointwise multiplication
\begin{equation}  \label{eq:def_act_multiplication}
  \alpha_\chi f(x)\defeq \chi(x)\cdot f(x)
\end{equation}
for all \(\chi\in\hat{G}\), \(f\in\Sch(G)\).  These two representations
are smooth in the sense of~\cite{Meyer:Smooth}.  For direct
products of groups, we get \(\Sch(G_1\times G_2)\cong
\Sch(G_1)\hot\Sch(G_2)\).  The Fourier transform for~\(G\) is an
isomorphism
\[
\Fourier\colon \Sch(G)\congto \Sch(\hat{G}).
\]
It intertwines convolution and pointwise multiplication and the
representations of \(G\) and~\(\hat{G}\) on \(\Sch(G)\) and \(\Sch(\hat{G})\)
described above.  If~\(G\) is elementary Abelian, then these assertions
are well-known if we take into account that bornological and
topological analysis are equivalent in this case by
\cites{Meyer:Born_Top, Meyer:Smooth}.  The extension to
arbitrary~\(G\) is straightforward using the direct union definition of
\(\Sch(G)\).

Given two bornological vector spaces \(V,W\) we write \(V\boin W\) if
\(V\subseteq W\) and the inclusion map \(V\to W\) is bounded.  In most
cases, both \(V\) and~\(W\) will be spaces of functions on some group~\(G\),
so that the set-theoretical condition \(V\subseteq W\) is meaningful.
Suppose that we are given a set \((V_i)_{i\in I}\) of bornological
vector spaces and a bornological vector space~\(W\) such that \(V_i\boin
W\).  In our applications, we can take~\(W\) to be the set of functions
on some group~\(G\).  Then the intersection \(\bigcap_{i\in I} V_i\) is a
vector subspace of~\(W\).  By convention, we equip \(\bigcap_{i\in I}
V_i\) with the \emph{intersection bornology}, which consists of subsets
that are bounded in~\(V_i\) for all \(i\in I\).

We let~\(G\) be a locally compact Abelian group equipped with a positive
quasi-character, that is, a continuous group homomorphism
\(\abs{\blank}\colon G\to\R\inv_+\).  Let
\[
G\one\defeq \{g\in G\mid \abs{g}=1\}.
\]
For each \(s\in\C\), the function \(x\mapsto\abs{x}^s\) is a
quasi-character on~\(G\).  It is a unitary character if~\(s\) is purely
imaginary.  Multiplication by a quasi-character is an algebra
homomorphism with respect to convolution.  Hence
\[
\Sch(G)_s \defeq \{f\colon G\to\C\mid f\cdot\abs{x}^s \in \Sch(G)\}
\]
is a convolution algebra as well.  The algebra \(\Sch(G)_s\) only
depends on the real part of~\(s\) because \(\Sch(G)\) is closed under
multiplication by characters.

Let \(I_+\) and~\(I_-\) be two intervals in~\(\R\).  We let \(I_\cap\defeq
I_+\cap I_-\) and \(I_\cup\defeq I_-\cup I_+\).  We always assume in this
situation that \(I_\cap\neq\emptyset\)\mdash so that~\(I_\cup\) is an
interval as well\mdash and that \(\sup I_+\ge \sup I_-\).  In the
following definition, we also assume that~\(I_+\) is bounded below and
not above and that~\(I_-\) is bounded above and not below.
Thus~\(I_\cap\) is bounded and \(I_\cup=\R\).  We mainly need the
intervals \(I_+=\ooival{0,\infty}\), \(I_-=\ooival{-\infty,1}\),
\(I_\cup=\R\) and \(I_\cap=\ooival{0,1}\) in this section.  In the global
case, we have to work with \(I_>\defeq \ooival{1,\infty}\) and
\(I_<\defeq\ooival{-\infty,0}\) instead to avoid the critical strip
\([0,1]\).

\begin{definition}  \label{def:weighted_Sch}
  For any non-empty interval \(I\subseteq\R\), let
  \[
  \Sch(G)_I \defeq \bigcap_{s\in I} \Sch(G)_s.
  \]
  Let \(K\subseteq\R\inv_+\) be a non-empty closed interval.  We let
  \(\Sch(G;K)_I\subseteq \Sch(G)_I\) be the subspace of functions
  supported in~\(K\), equipped with the subspace bornology.  If
  \(I_+,I_-,I_\cup,I_\cap\) are as above, we abbreviate \(\Sch(G)_+\defeq
  \Sch(G)_{I_+}\), etc., and we let
  \begin{alignat*}{2}
    \Sch_\comp(G)_\cup &\defeq
    \varinjlim_{n\to\infty} \Sch(G;[1/n,n])_\cup,
    &\qquad
    \Sch_\comp(G)_+ &\defeq
    \varinjlim_{n\to\infty} \Sch(G;\ocival{0,n})_+,
    \\
    \Sch_\comp(G)_\cap &\defeq
    \Sch(G)_\cap,
    &\qquad
    \Sch_\comp(G)_- &\defeq
    \varinjlim_{n\to\infty} \Sch(G;\coival{1/n,\infty})_-.
  \end{alignat*}
\end{definition}

The variants \(\Sch_\comp(G)_I\) are mainly used in the global function
field case.  They allow us to remember the compact support properties
of Schwartz functions in that case.  If~\(G\one\) is compact, then
\(\Sch_\comp(G)_\cup\) is nothing but the space \(\CCINF(G)\) of compactly
supported smooth functions on~\(G\).

It is clear that \(\Sch_\comp(G)_I\boin \Sch(G)_I\) for all~\(I\) and that
\(\Sch(G)_{I'}\boin\Sch(G)_I\) and
\(\Sch_\comp(G)_{I'}\boin\Sch_\comp(G)_I\) if \(I\subseteq I'\).  In the
following, we denote inclusion maps of this form by~\(\iota\).

Let \(\alpha\le\beta\).  We claim that
\begin{equation}  \label{eq:Sch_weighted_convex}
  \Sch(G)_{[\alpha,\beta]} =
  \Sch(G)_\alpha\cap\Sch(G)_\beta =
  \{f\colon G\to\C \mid f\cdot (\abs{x}^\alpha+\abs{x}^\beta)
  \in \Sch(G)\}
\end{equation}
as bornological vector spaces.  It is evident that we get a true
statement if we write~\(\boin\) instead of~\(=\)
in~\eqref{eq:Sch_weighted_convex}.  To prove the claim, we let~\(S\) be
bounded in the third space, that is, \(S\cdot
(\abs{x}^\alpha+\abs{x}^\beta)\) is bounded in \(\Sch(G)\).  For any
\(s\in [\alpha,\beta]\), the function \(w_s\defeq
\abs{x}^s/(\abs{x}^\alpha+\abs{x}^\beta)\) is a bounded smooth
function on~\(G\) whose derivatives remain bounded.  Therefore,
multiplication by~\(w_s\) is a bounded operator on \(\Sch(G)\) and \(S\cdot
\abs{x}^s\) is bounded in \(\Sch(G)\).  This implies the claim.

We can write any interval~\(I\) as an increasing union of a sequence of
compact intervals \([\alpha_n,\beta_n]\).  Thus
\[
\Sch(G)_I = \bigcap_{n\in\N} \Sch(G)_{\alpha_n}\cap\Sch(G)_{\beta_n}.
\]
Therefore, \(\Sch(G)_I\) is a nuclear, locally multiplicatively convex
Fréchet algebra if~\(G\) is elementary Abelian.  Moreover, the
representations of \(G\) and~\(\hat{G}\) on \(\Sch(G)_I\) described above
are smooth in this case.  For arbitrary~\(G\), the space \(\Sch(G)_I\) is
the direct union of the corresponding spaces \(\Sch(U/k)_I\) for proper
elementary Abelian subquotients of~\(G\).  Hence it is again a nuclear,
complete, convex bornological algebra and the representations of \(G\)
and~\(\hat{G}\) on it are smooth.  This also holds for the space
\(\Sch_\comp(G)_I\), we leave the verification to the reader.

\begin{lemma}  \label{lem:weighted_Sch_extension}
  Let \(I_+,I_-\subseteq\R\) be as above.  There is a bornological
  extension
  \[
  \Sch(G)_\cup
  \overset{\bigl(
    \begin{smallmatrix} \iota \\ \iota \end{smallmatrix}
  \bigr)}\into
  \Sch(G)_+ \oplus \Sch(G)_-
  \overset{(\iota, -\iota)}\prto
  \Sch(G)_\cap,
  \]
  which has a \(G\one\)\nbd{}equivariant bounded linear section and
  a similar extension with~\(\Sch_\comp\) instead of~\(\Sch\).
\end{lemma}

\begin{proof}
  It is clear that \(\Sch(G)_\cup = \Sch(G)_+\cap\Sch(G)_-\) as
  bornological vector spaces.  This means that the sequence is
  bornologically exact at \(\Sch(G)_\cup\) and
  \(\Sch(G)_+\oplus\Sch(G)_-\).  It remains to construct a bounded
  linear section for \((\iota,-\iota)\).  This implies bornological
  exactness at \(\Sch(G)_\cap\).  We assume that \(\sup I_-\le \sup I_+\).
  
  Let \(\phi\in\CCINF(\R_+)\) be as before and pull it back to a
  \(G\one\)\nbd{}invariant function on~\(G\) by \(\phi(g)\defeq
  \phi(\abs{g})\).  For any \(G\)\nbd{}equivariant differential
  operator~\(D\), the function \(D(\phi)\colon G\to\C\) is bounded and
  supported in the region \(\abs{x}\le R\) for some \(R\in\R\).  Hence
  \(\phi\cdot\abs{x}^s\) and \((1-\phi)\cdot\abs{x}^{-s}\) are multipliers
  of \(\Sch(G)\) for all \(s\ge0\).  That is, pointwise multiplication
  with these functions is a bounded linear operator on \(\Sch(G)\).  It
  follows that the linear map
  \[
  \sigma\colon \Sch(G)_\cap\to\Sch(G)_+\oplus\Sch(G)_-,
  \qquad f\mapsto (\phi\cdot f,-(1-\phi)\cdot f),
  \]
  is well-defined and bounded.  It is the desired section for
  \((\iota,-\iota)\).  By construction, \(\phi f\) vanishes for
  \(\abs{x}\to\infty\) and \((1-\phi)f\) vanishes for \(\abs{x}\to0\).
  Hence the same map works for \(\Sch_\comp\) instead of~\(\Sch\) as well.
\end{proof}

The operator~\(J\) defined as in~\eqref{eq:def_J} evidently is a
bornological isomorphism
\[
\Sch(G)_I\congto\Sch(G)_{1-I},
\qquad
\Sch_\comp(G)_I\congto \Sch_\comp(G)_{1-I}
\]
for all intervals~\(I\) and \(1-I=\{1-x\mid x\in I\}\).

Whenever we prove the nuclearity of an operator, we eventually reduce
the problem to the following elementary lemma:

\begin{lemma}  \label{lem:mother_of_nuclearity}
  Suppose that the norm homomorphism \(G\to\R\inv_+\) is a proper map of
  topological spaces.  Let \(I\subseteq\R\) be an interval and let
  \(\alpha\in I\) be arbitrary.  Let~\(\phi\) be as above and let
  \(h\in\Sch(G)_I\).  The operators \([\IN\lambda(f),M_\phi]\)
  \(\IN\lambda(f)M_h\) and \(M_h\IN\lambda(f)\) extend to uniformly
  nuclear operators from \(L^2(G)_\alpha\) to \(\Sch(G)_I\) for
  \(f\in\Sch(G)_I\).
\end{lemma}

\begin{proof}
  We describe these operators by integral kernels as in the proof of
  Lemma~\ref{lem:trace_trivial}.  The crucial observation is that
  under our hypothesis on~\(f\), the integral kernel belongs to
  \begin{multline*}
    \Sch(G)_I\hot\Sch(G)_{-I} \cong
    \bigcap_{\alpha,\beta\in I} \Sch(G)_\alpha\hot\Sch(G)_{-\beta}
    \\ \cong
    \{f\colon G^2\to\C\mid
      \text{\(f\cdot\abs{x}^\alpha\abs{y}^{-\beta}\in\Sch(G^2)\) for all
        \(\alpha,\beta\in I\)
      }
    \}.
  \end{multline*}
  This is shown as in the proof of Lemma~\ref{lem:trace_trivial}
  and easily implies the assertion.
\end{proof}

\subsection{Estimates on Schwartz functions}
\label{sec:equivariant_Fourier}

Let \(G\defeq\Adel\inv_S\).  We compare the space \(\Sch(\Adel_S)\) of
Bruhat-Schwartz functions on the additive group with the weighted
spaces of Bruhat-Schwartz functions \(\Sch(G)_I\) constructed above.

\begin{proposition}  \label{pro:Fourier_estimate}
  There are bounded injective maps
  \begin{alignat*}{2}
    \iota &\colon \Sch(\Adel_S)\to\Sch(G)_+,
    &\qquad
    f&\mapsto f|_G,
    \\
    \iota\EFoi &\colon \Sch(\Adel_S)\to\Sch(G)_-,
    &\qquad
    f&\mapsto J(\Fourier^\ast(f)|_G).
  \end{alignat*}
  The operator~\(\EFou\) defines a bounded linear operator on
  \(\Sch(G)_I\) if \(I\subseteq\ooival{0,\infty}\).  The operator~\(\EFoi\)
  defines a bounded linear operator on \(\Sch(G)_I\) if
  \(I\subseteq\ooival{-\infty,1}\).  Thus~\(\EFou\) is an automorphism of
  \(\Sch(G)_I\) if \(I\subseteq\ooival{0,1}\).
  
  The same assertions hold for \(\Sch_\comp\) instead of~\(\Sch\) if~\(\GF\)
  is a global function field.
\end{proposition}

\begin{proof}
  First we have to describe \(\Sch(G)_I\) explicitly.  The group~\(G\) is
  compactly generated.  Hence \(G/k\) is elementary Abelian whenever it
  is a Lie group.  Thus \(\Sch(G)\) is the direct union of the spaces
  \(\Sch(G/k)\), where~\(k\) runs through the smooth compact subgroups
  of~\(G\) in the terminology of~\cite{Meyer:Smooth}.  We define
  \[
  \norm{x} = \max \{\abs{x}_v,\abs{x}^{-1}_v\mid v\in S\}
  \qquad\text{for \(x\in G\).}
  \]
  This function measures rapid decay on \(\Sch(G)_I\).  Let \(RD(G)_+\) be
  the space of all functions \(f\colon G\to\C\) for which \(f(x)\cdot
  \abs{x}^\alpha (1+\ln{}\norm{x})^\beta \in C_0(G)\) for all
  \(\alpha>0\), \(\beta\in\R\).  A set~\(S\) of functions on~\(G\) is bounded
  in \(RD(G)_+\) if the set of functions
  \(f(x)\cdot\abs{x}^\alpha(1+\ln{}\norm{x})^\beta\) for \(f\in S\) is
  bounded in \(C_0(G)\) for all \(\alpha>0\), \(\beta\in\R\).  It is
  straightforward to see that \(\Sch(G)_+\) is the smoothening of the
  regular representation of~\(G\) on \(RD(G)_+\)
  (see~\cite{Meyer:Smooth} for the definition of the smoothening).
  
  We extend \(\abs{x}^\alpha (1+\ln{}\norm{x})^\beta\) to a function
  on~\(\Adel_S\) by~\(0\) on \(\Adel_S\setminus \Adel\inv_S\).  If
  \(\norm{x}\) is large and \(\abs{x}_v\) remains bounded for all \(v\in
  S\), then \(\abs{x}\approx \norm{x}^{-1}\), so that \(\abs{x}^\alpha
  (1+\ln{}\norm{x})^\beta\) is small.  Hence we have a continuous
  function on~\(\Adel_S\).  It also has polynomial growth for
  \(x\to\infty\), so that \(\Sch(\Adel_S)\boin RD(G)_+\).  The regular
  representation of~\(G\) on \(\Sch(\Adel_S)\) is easily seen to be
  smooth.  Hence the universal property of the smoothening yields
  \(\Sch(\Adel_S)\boin \Sch(G)_+\).  Since~\(\Fourier\) is a bornological
  isomorphism on \(\Sch(\Adel_S)\) and~\(J\) is a bornological isomorphism
  between \(\Sch(G)_\pm\), we also obtain that \(\iota\EFoi\) is a bounded
  linear map from \(\Sch(\Adel_S)\) to \(\Sch(G)_-\).

  Let \(I\subseteq\ooival{-\infty,1}\) be an interval.  Since the left
  regular representation on \(\Sch(G)_I\) is smooth, \(\Sch(G)_I\) is an
  essential module over \(\CCINF(G)\) (see~\cite{Meyer:Smooth}).
  This means that
  \[
  \CCINF(G)\hot_{\CCINF(G)} \Sch(G)_I \cong \Sch(G)_I.
  \]
  Hence the following defines a \(G\)\nbd{}equivariant bounded linear
  map:
  \begin{multline*}
    \Sch(G)_I \overset{\cong}\to
    \CCINF(G)\hot_{\CCINF(G)} \Sch(G)_I
    \subseteq
    \Sch(\Adel_S)\hot_{\CCINF(G)} \Sch(G)_I
    \\ \overset{\iota\EFoi\hot\ID}\longrightarrow
    \Sch(G)_I\hot_{\CCINF(G)} \Sch(G)_I
    \overset{\ast}\to
    \Sch(G)_I.
  \end{multline*}
  The last map above is the convolution on~\(G\).  The above composition
  sends \(f_1\ast f_2\) with \(f_1,f_2\in\CCINF(G)\) to \(\EFoi (f_1)\ast
  f_2=\EFoi(f_1\ast f_2)\).  Hence it extends~\(\EFoi\) on \(\CCINF(G)\).
  Being bounded, it must agree with the restriction of~\(\EFoi\) on
  \(L^2(\Adel_S,dx)\) to \(\Sch(\Adel_S)\).  A similar argument shows
  that~\(\EFou\) restricts to a bounded linear operator on \(\Sch(G)_I\)
  if \(I\subseteq\ooival{0,\infty}\).  Therefore, if
  \(I\subseteq\ooival{0,1}\), then both \(\EFoi\) and~\(\EFou\) are bounded
  on \(\Sch(G)_I\).
  
  If~\(\GF\) is a global function field, then
  \(\Sch(\Adel_S)=\CCINF(\Adel_S)\).  Thus functions in \(\Sch(\Adel_S)\)
  vanish for \(\abs{x}\gg1\).  This allows us to prove the corresponding
  statements for \(\Sch_\comp\).
\end{proof}

We let
\[
\Twist(G)\defeq
\Sch(G)_+\cap \EFou(\Sch(G)_-)
=
\bigcap_{\alpha>0} \Sch(G)_\alpha \cap
\bigcap_{\alpha<1} \EFou(\Sch(G)_\alpha).
\]
Thus \(\Twist(G)\boin L^2(G,\abs{x}\,d\inv x)\) and a subset~\(T\) of
\(L^2(G,\abs{x}\,d\inv x)\) is bounded in \(\Twist(G)\) if and only if
\(T\cdot\abs{x}^\alpha\) is bounded in \(\Sch(G)\) for \(\alpha>0\) and
\(\EFoi(T)\cdot\abs{x}^\alpha\) is bounded in \(\Sch(G)\) for \(\alpha<1\).
Since \(\EFoi=J\Fourier^\ast\) is a bornological isomorphism on
\(\Sch(G)_\alpha\) for \(\alpha\in\ooival{0,1}\), these two conditions
hold if and only if \(T\abs{x}^\alpha\) and \(\Fourier(T)\abs{x}^\alpha\)
are bounded in \(\Sch(G)_{1/2}\) for all \(\alpha\ge0\).  It is clear
that~\(G\) acts smoothly on \(\Twist(G)\) by the left regular
representation.  If we view \(\Sch(\Adel_S) \subseteq
L^2(\Adel_S,dx)\cong L^2(G,\abs{x}\,d\inv x)\), then
Proposition~\ref{pro:Fourier_estimate} implies \(\Sch(\Adel_S)\boin
\Twist(G)\).

If~\(\GF\) is a global function field, we also define
\[
\Twist_\comp(G) \defeq
\Sch_\comp(G)_+ \cap \EFou(\Sch_\comp(G)_-).
\]
We have \(\Sch(\Adel_S)\boin \Twist_\comp(G)\boin\Twist(G)\) in this
case.  We do \emph{not} define \(\Twist_\comp(G)\) if~\(\GF\) is an
algebraic number field because the above definition would just
produce~\(0\).  Whenever we make statements about \(\Twist_\comp\), we
tacitly assume~\(\GF\) to be a global function field.

\begin{lemma}  \label{lem:Twist_extension}
  There is a \(G\)\nbd{}equivariant bornological extension
  \begin{equation}  \label{eq:Twist_extension}
    \Twist(G)
    \overset{\bigl(
      \begin{smallmatrix} \iota \\ \iota\EFoi \end{smallmatrix}
    \bigr)}\into
    \Sch(G)_+ \oplus \Sch(G)_-
    \overset{(\iota, -\EFou\iota)}\prto
    \Sch(G)_\cap,
  \end{equation}
  which has a bounded \(G\one\)\nbd{}equivariant linear section and
  similarly for \(\Twist_\comp(G)\).
\end{lemma}

\begin{proof}
  Bornological exactness at \(\Twist(G)\) and in the middle is trivial.
  We only have to construct a bounded linear section for the map to
  \(\Sch(G)_\cap\).  Let~\(\phi\) be as in the proof of
  Lemma~\ref{lem:weighted_Sch_extension} and define
  \[
  \sigma\colon \Sch(G)_\cap \to \Sch(G)_+ \oplus \Sch(G)_-,
  \qquad
  \sigma(f)\defeq (\phi\cdot f,-\EFoi((1-\phi)\cdot f)).
  \]
  As in the proof of Lemma~\ref{lem:weighted_Sch_extension},
  multiplication by \(\phi\) and \(1-\phi\) are bounded linear maps from
  \(\Sch(G)_\cap\) to \(\Sch(G)_\pm\), respectively.  The map~\(\EFoi\) is a
  bounded linear operator on \(\Sch(G)_-\) by
  Proposition~\ref{pro:Fourier_estimate}.  Therefore, \(\sigma\) is a
  well-defined bounded linear map.  It is clear that
  \((\iota,-\EFou\iota)\circ\sigma=\ID\) and that~\(\sigma\) is
  \(G\one\)\nbd{}equivariant.
\end{proof}

Since~\(\Fourier\) is an automorphism of \(\Sch(\Adel_S)\), we expect it
to be a bornological isomorphism on \(\Twist(G)\) and \(\Twist_\comp(G)\)
as well.  This is indeed the case because
\begin{align*}
  \Fourier(\Sch(G)_+) = \Fourier J(\Sch(G)_-) = \EFou(\Sch(G)_-),
  \\
  \Fourier^\ast\circ\EFou(\Sch(G)_-) = J(\Sch(G)_-) = \Sch(G)_+.
\end{align*}

Let \(H\defeq\GF\inv_S\), so that \(G/H=\ICL_S\).  Since~\(H\) is a
cocompact subgroup of~\(G\one\), the regular representation of~\(G\one\)
on \(\CCINF(\ICL\one_S)\) is \(\Sch(G\one)\)\brd{}tempered.  We let
\[
\Twist(\ICL_S) \defeq
\CCINF(\ICL\one_S)\hot_{\Sch(G\one)} \Twist(G).
\]
It is easy to see that \(\Sch(G)_I\cong
\Sch(G\one)\hot\Sch(\R\inv_+)_I\) and similarly for \(\ICL_S\).  Hence
\[
\CCINF(\ICL\one_S) \hot_{\Sch(G\one)} \Sch(G)_I \cong \Sch(\ICL_S)_I
\]
for all intervals~\(I\).  Proposition~\ref{pro:Fourier_estimate}
implies that \(\EFoi\) and~\(\EFou\) induce bounded operators on
\(\Sch(\ICL_S)_I\) for \(I\subseteq\ooival{-\infty,1}\) and
\(I\subseteq\ooival{0,\infty}\), respectively.  Since the extension in
Lemma~\ref{lem:Twist_extension} splits \(G\one\)\nbd{}equivariantly,
\(\Twist(G)\) is a projective module over \(\Sch(G\one)\) and we obtain a
\(\ICL_S\)\nbd{}equivariant linearly split bornological extension
\begin{equation}  \label{eq:Twist_convariants}
  \Twist(\ICL_S)
  \overset{\bigl(
    \begin{smallmatrix} \iota \\ \iota\EFoi \end{smallmatrix}
  \bigr)}\into
  \Sch(\ICL_S)_+ \oplus \Sch(\ICL_S)_-
  \overset{(\iota, -\EFou\iota)}\prto
  \Sch(\ICL_S)_\cap.
\end{equation}
As a result, \(\Twist(\ICL_S) \cong \Sch(\ICL_S)_+\cap
\EFou(\Sch(\ICL_S)_-)\).  Analogous results hold for \(\Twist_\comp(G)\)
instead of \(\Twist(G)\).
Since~\(\EFoi\) is bounded on \(\Sch(\ICL_S)_-\), we have
\(\Sch(\ICL_S)_-\boin \EFou(\Sch(\ICL_S)_-)\) and hence
\[
\Sch(\ICL_S)_\cup =
\Sch(\ICL_S)_-\cap\Sch(\ICL_S)_+ \boin
\Twist(\ICL_S).
\]
Similarly, we get \(\Sch_\comp(\ICL_S)_\cup\boin\Twist_\comp(\ICL_S)\).

One can show that \(\Sch(\ICL_S)_\cup\) is a closed subspace of
\(\Twist(\ICL_S)\).  The quotient representation of~\(\ICL_S\) on
\(\Twist(\ICL_S)/\Sch(\ICL_S)_\cup\) is a spectral interpretation for
the poles of the \(S\)\nbd{}local \(L\)\nbd{}function for~\(\Adel_S\).
However, these poles do not have much geometric significance.  They
all lie in the region \(\RE\omega\le 0\).  The functional equation
relates them to the trivial zeros of infinite Euler products.  For
instance, if we take just the local field~\(\R\), then we obtain the
trivial zeros of the \(\zeta\)\nbd{}function.  Hence it does not seem
worthwhile to explore this spectral interpretation further.

Let~\(\phi\) and \(T\defeq M_\phi- \EFou\circ M_\phi\circ \EFoi = M_\phi
- \Fourier\circ M_{\check\phi} \circ \Fourier^{-1}\) be as before.

\begin{proposition}  \label{pro:T_on_Twist}
  The operator~\(T\) on \(L^2(\ICL_S)_{1/2}\) restricts to a bounded
  linear operator \(\Twist(\ICL_S)\to\Twist(\ICL_S)\) for algebraic
  number fields and \(\Twist_\comp(\ICL_S)\to\Twist_\comp(\ICL_S)\) for
  global function fields.

  The operators \(T\IN\lambda(f)\) and \(\IN\lambda(f)T\) are uniformly
  nuclear operators \(L^2(\ICL_S)_{1/2}\to\Twist(\ICL_S)\) for
  \(f\in\CCINF(\ICL_S)\).  The same holds for \(\Twist_\comp(\ICL_S)\)
  if~\(\GF\) is a global function field.

\end{proposition}

\begin{proof}
  Let \(\check\phi(x)\defeq \phi(1/x)\), then \(JM_\phi
  J=M_{\check\phi}\).  Recall that \(M_{1-\phi}\) and \(M_{\check\phi}\)
  are bounded operators from \(\Sch(\ICL_S)_+\) to \(\Sch(\ICL_S)_\cup\).
  Since \(\Sch(\ICL_S)_\cup\boin\Twist(\ICL_S)\boin\Sch(\ICL_S)_+\), the
  operators \(M_{\check\phi}\) and \(M_{1-\phi}\) and hence also
  \(M_\phi=\ID-M_{1-\phi}\) are bounded on \(\Twist(\ICL_S)\).
  Since~\(\Fourier\) is bounded on \(\Twist(\ICL_S)\) as well, \(T\) is
  bounded on \(\Twist(\ICL_S)\).
  
  We have \([\IN\lambda(f),T] = [\IN\lambda(f),M_\phi] + \Fourier
  [\IN\lambda(Jf),M_{1-\check\phi}] \Fourier^{-1}\).  The
  operators~\(\Fourier^{\pm1}\) are bounded on \(L^2(\ICL_S)_{1/2}\) and
  \(\Twist(\ICL_S)\).  And \([\IN\lambda(f),M_\phi]\) and
  \([\IN\lambda(Jf),M_{\check\phi}]\) are uniformly nuclear operators
  from \(L^2(\ICL_S)_{1/2}\) to \(\Sch(\ICL_S)_\cup\boin\Twist(\ICL_S)\)
  by Lemma~\ref{lem:mother_of_nuclearity}.  Hence \([\IN\lambda(f),T]\)
  is a uniformly nuclear operator from \(L^2(\ICL_S)_{1/2}\) to
  \(\Twist(\ICL_S)\) for \(f\in\CCINF(\ICL_S)\).
  
  For \(\IN\lambda(f)\circ T\), we first simplify the problem.  Since
  \(\Twist(\ICL_S)\) is nuclear and \(L^2(\ICL_S)_{1/2}\) is a Banach
  space, uniform nuclearity and uniform boundedness are equivalent in
  our situation.  A map into \(\Twist(\ICL_S)\) is bounded if and only
  if it is bounded as a map to \(\Sch(\ICL_S)_+\) and
  \(\EFou(\Sch(\ICL_S)_-)\).  Since \(\Fourier T\Fourier^{-1}= -T\), we have
  \(\Fourier \IN\lambda(f) T\Fourier^{-1} = -\IN\lambda(Jf) T\).  The
  map~\(\Fourier\) gives bornological isomorphisms
  \[
  \Sch(\ICL_S)_+\cong \EFou(\Sch(\ICL_S)_-),
  \qquad
  L^2(\ICL_S)_{1/2}\cong L^2(\ICL_S)_{1/2}.
  \]
  Hence \(\IN\lambda(f)T\) is bounded as a map to
  \(\EFou(\Sch(\ICL_S)_-)\) if (and only if) it is bounded as a map to
  \(\Sch(\ICL_S)_+\).  As in the proof of Lemma~\ref{lem:local_trace_I},
  we write
  \[
  \IN\lambda(f) T =
  -[M_\phi,\IN\lambda(f)] + [M_\phi, \IN\lambda(\EFou f)]\circ\EFoi.
  \]
  We have \(f\in\Sch(\ICL_S)_+\) and \(\EFou f\in\Sch(\ICL_S)_+\) as
  well by Proposition~\ref{pro:Fourier_estimate}.  Since~\(\EFoi\) is
  bounded on \(L^2(\ICL_S)_{1/2}\), the result follows from
  Lemma~\ref{lem:mother_of_nuclearity}.  The same argument works for
  \(\Twist_\comp(\ICL_S)\).
\end{proof}

\begin{corollary}  \label{cor:T_on_Twist}
  The operators \(\IN\lambda(f)T\) and \(T\IN\lambda(f)\) are uniformly
  nuclear for \(f\in\CCINF(\ICL_S)\) as operators on
  \(L^2(\ICL_S)_{1/2}\), \(\Twist(\ICL_S)\) and \(\Twist_\comp(\ICL_S)\) (if
  defined) and have the same traces as operators on these spaces.
\end{corollary}

\begin{proof}
  More generally, consider bornological vector spaces \(V\boin W\) and
  \(X\in\ell^1(W,V)= V\hot W'\).  Since \(V\boin W\), we can view~\(X\) as
  an element of \(\ell^1(W)\) or \(\ell^1(V)\).  Since the pairings
  \(V\times V'\to\C\) and \(W\times W'\to\C\) are compatible, the traces
  of~\(X\) in \(\ell^1(V)\) and \(\ell^1(W)\) agree on elementary tensors
  and hence everywhere.
\end{proof}

\subsection{Homology for smooth convolution algebras}
\label{sec:homology_smooth_convolution}

A \emph{smooth convolution algebra} on a locally compact group~\(G\) is
a complete convex bornological vector space \(A(G)\) of functions
on~\(G\) which is a bornological algebra under convolution, contains
\(\CCINF(G)\) as a dense subspace and has the property that the left
and right regular representations of~\(G\) on \(A(G)\) are smooth.  Thus
\(A(G)\) is an essential bimodule over \(\CCINF(G)\) in the sense
of~\cite{Meyer:Smooth}, that is,
\begin{equation}  \label{eq:AG_essential}
  \CCINF(G)\hot_{\CCINF(G)} A(G)\hot_{\CCINF(G)} \CCINF(G) \cong A(G).
\end{equation}
It follows that the approximate identities for \(\CCINF(G)\subseteq
A(G)\) constructed in~\cite{Meyer:Smooth} are approximate identities
for \(A(G)\) as well.  Thus we can speak of essential modules over
\(A(G)\) as in~\cite{Meyer:Smooth}.  Any module over \(A(G)\) becomes a
module over \(\CCINF(G)\) by restriction.  Since \(\CCINF(G)\subseteq
A(G)\) is dense, the module structure over \(\CCINF(G)\) determines the
\(A(G)\)\brd{}module structure uniquely.  A module over \(A(G)\) is
essential if and only if it is essential as a module over \(\CCINF(G)\)
because of~\eqref{eq:AG_essential}.  If~\(V\) is an essential module
over \(A(G)\), then the module structure is the integrated form of a
smooth representation of~\(G\) on~\(V\).  The \(A(G)\)\brd{}module
homomorphisms are the same as the \(G\)\nbd{}equivariant maps.  Hence:

\begin{theorem}  \label{the:AG_tempered_representations}
  The category of essential \(A(G)\)\brd{}modules is isomorphic to a full
  subcategory of the category of smooth representations of~\(G\).
\end{theorem}

We call such smooth representations \emph{\(A(G)\)\brd{}tempered}.  The
above remarks yield \(\Hom_{A(G)}(V,W)\cong\Hom_{\CCINF(G)}(V,W)\) for
all essential modules \(V\) and~\(W\) over \(A(G)\).
Equation~\eqref{eq:AG_essential} implies (with little extra work) that
\begin{align*}
  \Hom_{A(G)}(A(G),W) &\cong
  \Hom_{\CCINF(G)}(\CCINF(G),W),
  \\
  A(G)\hot_{A(G)} W &\cong
  \CCINF(G)\hot_{\CCINF(G)} W,
\end{align*}
if~\(W\) is a module over \(A(G)\).  That is, the smoothening and
roughening functors for modules over \(\CCINF(G)\) and \(A(G)\) are
compatible.

The functor \(V\mapsto A(G)\hot_{\CCINF(G)} V\) maps modules over
\(\CCINF(G)\) to essential modules over \(A(G)\) because \(A(G)\hot_{A(G)}
A(G)\cong A(G)\).  If \(V\) and~\(W\) are essential modules over
\(\CCINF(G)\) and \(A(G)\), respectively, then
\begin{multline*}
  \Hom_{A(G)}(A(G)\hot_{\CCINF(G)} V, W)
  \cong
  \Hom_{\CCINF(G)}(V, \Hom_{A(G)}(A(G),W))
  \\ \cong
  \Hom_{\CCINF(G)}(V, \Hom_{\CCINF(G)}(\CCINF(G),W))
  \cong
  \Hom_{\CCINF(G)}(V,W).
\end{multline*}
That is, the functor \(A(G)\hot_{\CCINF(G)}\blank\) is left adjoint to
the embedding of the category of \(A(G)\)\brd{}tempered smooth
representations in the category of all smooth representations.  We
call this the \emph{\(A(G)\)\brd{}temperation functor}.  It follows
immediately from the adjointness that \(A(G)\hot_{\CCINF(G)} V\cong V\)
if~\(V\) is already \(A(G)\)\brd{}tempered.

Before we can discuss derived functors, we have to describe our exact
category structure.  For the purposes of this article, the most
convenient choice is to work with \emph{linearly split extensions}.
Let~\(A\) be a complete convex bornological algebra with an approximate
identity and let~\(\Mod_A\) be the category of essential, complete
convex bornological left \(A\)\nbd{}modules.  That is, an object
of~\(\Mod_A\) is a complete convex bornological vector space~\(V\) with a
bounded homomorphism \(A\to\End(V)\) such that \(A\hot_A V\cong V\) or,
equivalently, the map \(A\hot V\to V\) is a bornological quotient map.
We remark that adjoint associativity yields an equivalence between
bounded linear maps \(A\to\End(V)\) and \(A\hot V\to V\), so that we can
describe \(A\)\nbd{}modules using either kind of structure.

An extension in~\(\Mod_A\) is called \emph{linearly split} if it splits
as an extension of bornological vector spaces.  It is evident that
this defines an exact category in the sense of Daniel Quillen
(\cite{Quillen:Higher}).  Hence we can form the derived
category~\(\Der_A\) of~\(\Mod_A\) in the usual fashion
(see~\cites{Keller:Handbook, Neeman:Derived_Exact}).  We define the
\emph{free essential module} on a complete convex bornological vector
space~\(V\) by \(A\hot V\).  \emph{We assume that free essential modules
are projective with respect to linearly split extensions.}  This holds
for \(\CCINF(G)\) and also for \(A(G)\), using the same proof as
in~\cite{Meyer:Smooth}.  Hence we can compute derived functors
using free essential bimodule resolutions.  The well-known bar
construction provides a natural such resolution for any essential
module.

Let \(F\colon \Mod_A\to\Mod_B\) be a functor.  Its total derived functor
is a functor \(\Der_A\to\Der_B\), that is, its values are equivalence
classes of chain complexes.  For instance, let~\(V\) be a right
\(A\)\nbd{}module.  Then we have a functor \(V\hot_A\blank\colon
\Mod_A\to\Mod_0\), where~\(\Mod_0\) is the category of complete convex
bornological vector spaces with no further structure.  By our choice
of exact category structure, the only extensions in~\(\Mod_0\) are the
direct sum extensions, so that~\(\Der_0\) is just the homotopy category
of chain complexes of complete convex bornological vector spaces.  We
denote the derived functor of \(V\hot_A\blank\) by \(V\Lhot_A\blank\).
Its values are homotopy classes of chain complexes.

Consider again the case \(\CCINF(G)\).  Let \(P\to\CCINF(G)\) be a
(linearly split) resolution of \(\CCINF(G)\) by free essential
\(\CCINF(G)\)\nbd{}bimodules.  If~\(V\) is an essential left
\(\CCINF(G)\)\nbd{}module, then \(P\hot_{\CCINF(G)} V\) is a linearly
split resolution of~\(V\) by free essential left
\(\CCINF(G)\)\nbd{}modules.  Hence this resolution can be used to
compute derived functors.  Something similar happens if~\(V\) is a chain
complex.  Then \(P\hot_{\CCINF(G)} V\) is a bicomplex, its total complex
is a projective chain complex and the natural map from this total
complex to~\(V\) is a quasi-isomorphism.  Thus morphisms in the derived
category are given by
\[
\Right\Hom_{\CCINF(G)}(V,W)\cong
\Hom_{\CCINF(G)}(P\hot_{\CCINF(G)} V,W),
\]
where \(\Hom_{\CCINF(G)}\) denotes the space of homotopy classes of
bounded \(\CCINF(G)\)\nbd{}linear chain maps \(V\to W\).  The right derived
functor \(\Right\Hom_{\CCINF(G)}(V,W)\) gives the morphism \(V\to W\) in
the derived category.  Similarly, we have
\[
W\Lhot_{\CCINF(G)} V\cong W\hot_{\CCINF(G)} P\hot_{\CCINF(G)} V.
\]

Now let \(A(G)\) be a smooth convolution algebra on~\(G\).  We have
observed above that \(\Mod(A(G))\) is a full subcategory of
\(\Mod(\CCINF(G))\).  Although \(A(G)\) usually is not flat as a module
over \(\CCINF(G)\), the following still holds in the cases we need:

\begin{proposition}  \label{pro:isocohomological}
  Let \(A(G)\) be a smooth convolution algebra on a locally compact
  group~\(G\).  Then the following are equivalent:
  \begin{enumerate}[(i)]
  \item \(A(G)\Lhot_{\CCINF(G)} A(G)\cong A(G)\);

  \item \(V\Lhot_{\CCINF(G)} W\cong V\Lhot_{A(G)} W\) for any chain
    complexes of essential right and left \(A(G)\)\brd{}modules \(V\) and~\(W\);
    
  \item \(\Right\Hom_{A(G)}(V,W)\cong \Right\Hom_{\CCINF(G)}(V,W)\) for
    any chain complexes of essential left \(A(G)\)\brd{}modules \(V\)
    and~\(W\).

  \end{enumerate}
\end{proposition}

We call the smooth convolution algebra \(A(G)\) \emph{isocohomological}
if it satisfies these equivalent conditions.

\begin{proof}
  Let~\(P\) be as above and define \(C\defeq A(G)\hot_{\CCINF(G)}
  P\hot_{\CCINF(G)} A(G)\).  Evidently, this is a complex of free
  essential \(A(G)\)\brd{}bimodules.  We have just observed that \(C\cong
  A(G)\Lhot_{\CCINF(G)} A(G)\).  Thus~\(C\) is a free bimodule resolution
  of \(A(G)\) if and only if the natural map \(C\to A(G)\) is a homotopy
  equivalence of chain complexes, that is, an isomorphism in \(\Der_0\).
  If this is the case, then~\(C\) can be used to compute derived
  functors on the category of \(A(G)\)\brd{}modules.  This yields (ii)
  and~(iii).  The converse implications are straightforward.
\end{proof}

\subsection{Weighted Schwartz algebras on \(\R\) and~\(\Z\) are
  isocohomological}
\label{sec:homology_R}

We equip \(G=\R\) with the norm homomorphism \(\abs{x}\defeq \exp(x)\) and
use this to define the weighted Schwartz algebras \(\Sch(\R)_I\) for
intervals \(I\subseteq\R\).  The following is a resolution of
\(\CCINF(\R)\) by free essential \(\CCINF(G)\)\nbd{}bimodules:
\[
0 \to
\CCINF(\R)\hot\CCINF(\R)
\overset{D}\longrightarrow
\CCINF(\R)\hot\CCINF(\R) \overset{\ast}\longrightarrow
\CCINF(\R).
\]
Here~\(\ast\) denotes convolution and \(D\defeq
(d/dx)\otimes\ID-(\ID\otimes d/dx)\).  Let \(V\) and~\(W\) be two smooth
representations of~\(\R\).  We view them as a right and left essential
module over \(\CCINF(\R)\) in the standard way
(see~\cite{Meyer:Smooth}).  We use the above free resolution to
identify \(V\Lhot_{\CCINF(\R)} W\) with the chain complex of
bornological vector spaces
\begin{equation}  \label{eq:R_Lhot}
  V\hot W \overset{D_\ast}\longrightarrow V\hot W.
\end{equation}
The boundary map~\(D_\ast\) is given by \(D_\ast\defeq
D_V\otimes\ID_W-\ID_V\otimes D_W\), where \(D_V\) and~\(D_W\) are the
generators of the representations of~\(\R\) on \(V\) and~\(W\),
respectively.

\begin{lemma}  \label{lem:weighted_Sch_Lhot}
  For any two intervals \(I,I'\subseteq\R\), we have a natural
  isomorphism \(\Sch(\R)_I\Lhot_{\CCINF(\R)} \Sch(\R)_{I'} \cong
  \Sch(\R)_{I\cap I'}\) with \(\Sch(G)_\emptyset\defeq\{0\}\) by
  convention.
\end{lemma}

\begin{proof}
  We can easily identify \(\Sch(\R)_I\hot\Sch(\R)_{I'}\) with the space
  \[
  \Sch(\R^2)_{I\times I'}\defeq
  \{f\colon \R^2\to\C \mid
    \text{\(f\cdot \exp(\alpha x+ \beta y)\in\Sch(\R^2)\) 
    for all \((\alpha,\beta)\in I\times I'\)}
  \}.
  \]
  The operator~\(D_\ast\) in~\eqref{eq:R_Lhot} is the directional
  derivative operator \(D_\ast f= \partial_x f-\partial_y f\).  Thus it
  is the generator of the \(1\)\nbd{}parameter group \(\tau_t
  f(x,y)\defeq f(x+t,y-t)\).  It is evident that~\(D_\ast\) is injective,
  that is, there are no \(\tau\)\nbd{}invariant functions in
  \(\Sch(\R^2)_{I\times I'}\).
  
  If \(I\cap I'\neq\emptyset\), then \(\Sch(\R)_I\boin\Sch(\R)_{I\cap
    I'}\) and \(\Sch(\R)_{I'}\boin\Sch(\R)_{I\cap I'}\).  Since the
  latter is a convolution algebra, the convolution defines a bounded
  bilinear map
  \[
  \mu \colon \Sch(\R^2)_{I\times I'} \cong
  \Sch(\R)_I\hot\Sch(\R)_{I'} \to \Sch(\R)_{I\cap I'},
  \qquad
  \mu f(x)\defeq \int_\R f(x-t,t)\,dt.
  \]
  We have to show that
  \[
  \Sch(\R^2)_{I\times I'} \overset{D_\ast}\longrightarrow
  \Sch(\R^2)_{I\times I'} \overset{\mu}\longrightarrow
  \Sch(\R)_{I\cap I'}
  \]
  is a linearly split bornological extension.  We can construct a
  section for~\(\mu\) by composing the section for \((\iota,\iota)\colon
  \Sch(\R)_I\oplus\Sch(\R)_{I'}\to \Sch(\R)_{I\cap I'}\) defined in the
  proof of Lemma~\ref{lem:weighted_Sch_extension} with the sections
  \(\Sch(\R)_I\to \Sch(\R)_I\hot\CCINF(\R)\subseteq \Sch(\R^2)_{I\times
    I'}\) and \(\Sch(\R)_{I'}\to\CCINF(\R)\hot\Sch(\R)_{I'}\subseteq
  \Sch(\R^2)_{I\times I'}\) for the convolution map defined in
  \cite{Meyer:Smooth}*{Proposition~4.7}.  Thus it remains to prove
  that~\(D_\ast\) is a bornological isomorphism onto the kernel
  of~\(\mu\).  It is clear that~\(D_\ast\) is an injective bounded map
  into \(\Ker\mu\), so we merely have to find \(\sigma\colon
  \Ker\mu\to\Sch(\R^2)_{I\times I'}\) with \(D_\ast\circ\sigma=\ID\).
  Similarly, if \(I\cap I'=\emptyset\), then we have to find such a
  map~\(\sigma\) defined on all of \(\Sch(\R^2)_{I\times I'}\).
  Since~\(D_\ast\) involves differentiation, integration provides
  obvious candidates for its inverse:
  \begin{align*}
    \sigma_+ f(x,y) &\defeq -\int_0^\infty f(x+t,y-t) \,dt,
    \\
    \sigma_- f(x,y) &\defeq \int_{-\infty}^0 f(x+t,y-t) \,dt.
  \end{align*}
  Formally, both operators satisfy \(D_\ast\circ\sigma_\pm=\ID\).
  Moreover,
  \[
  \sigma_-(f)(x,y)-\sigma_+(f)(x,y) = \mu(f)(x+y),
  \]
  so that both maps coincide on \(\Ker\mu\).  It is easy to see
  that \(\sigma_+\) and~\(\sigma_-\) are bounded on
  \(\Sch(\R^2)_{(\alpha,\beta)}\) for \(\alpha\ge\beta\) and
  \(\alpha\le\beta\), respectively.  If \(I\cap I'=\emptyset\), then
  either \(\alpha<\beta\) for all \((\alpha,\beta)\in I\cap I'\) or
  \(\alpha>\beta\) for all \((\alpha,\beta)\in I\cap I'\).  Hence either
  \(\sigma_+\) or~\(\sigma_-\) is a bounded inverse for~\(D_\ast\).  For
  \(I\cap I'\neq\emptyset\), the maps~\(\sigma_\pm\) agree on \(\Ker\mu\).
  The restriction \(\sigma_\pm|_{\Ker\mu}\) is bounded with respect to
  the bornology from \(\Sch(\R^2)_{\alpha,\beta}\) for any
  \(\alpha,\beta\in\R\).  Thus \(\sigma_\pm|_{\Ker\mu}\) is a bounded
  inverse for~\(D_\ast\).
\end{proof}

Now let \(G=\Z\subseteq\R\).  We use the same norm homomorphism as above
to define the spaces \(\Sch_\comp(\Z)_I\).  There is a similar explicit
chain complex that computes homology for~\(\Z\) and we can use an
analogous argument to obtain
\[
\Sch_\comp(\Z)_I\Lhot_{\CCINF(\Z)} \Sch_\comp(\Z)_{I'} \cong
\Sch_\comp(\Z)_{I\cap I'}.
\]
We leave the details to the reader.

\begin{corollary}  \label{cor:weighted_Sch_isocohomological}
  The smooth convolution algebras \(\Sch(\R)_I\) and \(\Sch_\comp(\Z)_I\)
  are isocohomological for all intervals~\(I\).
\end{corollary}

\begin{theorem}  \label{the:Sch_isocohomological}
  The Bruhat-Schwartz algebra \(\Sch(G)\) is isocohomological for any
  second countable locally compact Abelian group~\(G\).
\end{theorem}

\begin{proof}
  Corollary~\ref{cor:weighted_Sch_isocohomological} yields that
  \(\Sch(\R)\) and \(\Sch(\Z)\) are isocohomological.  For compact groups,
  we have \(\Sch(K)=\CCINF(K)\), so that there is nothing to prove at
  all.  Let \(G=G_1\times G_2\) and suppose that \(\Sch(G_1)\) and
  \(\Sch(G_2)\) are isocohomological.  Then
  \begin{multline*}
    \Sch(G)\Lhot_{\CCINF(G)} \Sch(G)
    \cong
    \Sch(G_1)\hot\Sch(G_2)\Lhot_{\CCINF(G_1)\hot\CCINF(G_2)}
    \Sch(G_1)\hot\Sch(G_2)
    \\ \cong
    (\Sch(G_1)\Lhot_{\CCINF(G_1)} \Sch(G_1))\hot
    (\Sch(G_2)\Lhot_{\CCINF(G_2)} \Sch(G_2))
    \\ \cong
    \Sch(G_1)\hot\Sch(G_2)
    \cong
    \Sch(G),
  \end{multline*}
  that is, \(\Sch(G)\) is isocohomological as well.  Any elementary
  Abelian group is a direct product \(K\times \Z^n\times \R^m\) with
  compact~\(K\).  Hence we get the assertion for elementary Abelian
  groups.  Finally, for arbitrary~\(G\) both \(\CCINF(G)\) and \(\Sch(G)\)
  are direct unions of the corresponding algebras \(\CCINF(U/k)\) and
  \(\Sch(U/k)\) on proper elementary Abelian subquotients \(U/k\) and
  these subspaces are bornological direct summands.  If we compute
  \(\Sch(G)\Lhot_{\CCINF(G)}\Sch(G)\) using, say, the bar resolution, we
  see easily that it is the direct union of the corresponding
  complexes for proper elementary Abelian subquotients of~\(G\) and that
  the latter are direct summands of the complex for~\(G\) itself.  Hence
  the contractibility of the subcomplexes implies the contractibility
  of the limit, so that \(\Sch(G)\) is also isocohomological.
\end{proof}

\subsection{The coinvariant space of \(\Sch(\Adel_S)\)}
\label{sec:S_local_coinvariants}

\begin{lemma}  \label{lem:Sch_local}
  If~\(\LF\) is an infinite local field, then the map
  \(\Sch(\LF)\to\Twist(\LF\inv)\) is a bornological isomorphism.
  If~\(\LF\) is a finite local field, then the map
  \(\Sch(\LF)\to\Twist_\comp(\LF\inv)\) is a bornological isomorphism.
\end{lemma}

\begin{proof}
  We have already observed that \(\Sch(\LF)\boin\Twist(\LF\inv)\) in the
  infinite case and \(\Sch(\LF)\boin\Twist_\comp(\LF\inv)\) in the
  finite case.  It remains to show that any bounded subset~\(T\) of
  \(\Twist(\LF\inv)\) or \(\Twist_\comp(\LF\inv)\) comes from a bounded
  subset of \(\Sch(\LF)\).  Suppose first that~\(\LF\) is infinite, that
  is, either \(\R\) or~\(\C\).  Then \(T\cdot\abs{x}^s\) and
  \(\Fourier(T)\cdot\abs{x}^s\) are bounded subsets of \(L^2(\LF,dx)\) for
  all \(s\ge0\).  It is well-known that this implies that~\(T\) is a
  bounded subset of \(\Sch(\LF)\).  Similarly, for finite~\(\LF\) the
  functions in \(T\) and \(\Fourier(T)\) have a common compact support
  in~\(\LF\).  This implies that~\(T\) is bounded in \(\Sch(\LF)\).
\end{proof}

Lemma~\ref{lem:Sch_local} provides us with a linearly split
extension
\begin{equation}  \label{eq:Sch_infinite_extension}
  \Sch(\LF)
  \overset{\bigl(
    \begin{smallmatrix} \iota \\ \iota\EFoi \end{smallmatrix}
  \bigr)}\into
  \Sch(\LF\inv)_+ \oplus \Sch(\LF\inv)_-
  \overset{(\iota, -\EFou\iota)}\prto
  \Sch(\LF\inv)_\cap
\end{equation}
for infinite~\(\LF\) and a similar extension for finite~\(\LF\).
Equivalently, we have a quasi-isomorphism from
\(\Sch(\LF)\) into the chain complex
\[
\Sch(\LF\inv)_+\oplus\Sch(\LF\inv)_-
\overset{(\iota, -\EFou\iota)}\longrightarrow
\Sch(\LF\inv)_\cap
\]
supported in degrees \(0,-1\).

\begin{lemma}  \label{lem:Sch_local_reduce}
  There is a natural isomorphism
  \begin{equation}  \label{eq:Sch_local_reduce_infinite}
    \Sch(\LF\inv)_+\Lhot_{\CCINF(\LF\inv)} \Sch(\LF) \cong \Sch(\LF)_+
  \end{equation}
  for an infinite local field~\(\LF\) and
  \begin{equation}  \label{eq:Sch_local_reduce_finite}
    \Sch_\comp(\LF\inv)_+\Lhot_{\CCINF(\LF\inv)} \Sch(\LF) \cong
    \Sch_\comp(\LF)_+
  \end{equation}
  for a finite local field~\(\LF\).
\end{lemma}

\begin{proof}
  In the infinite case, we have
  \(\LF\inv=\LF\one\times\R\inv_+\cong\LF\one\times\R\) and hence
  \(\Sch(\LF\inv)_I\cong\Sch(\LF\one)\hot\Sch(\R)_I\).  Since~\(\LF\one\)
  is compact, we get a natural isomorphism
  \(\Sch(\LF\inv)_I\Lhot_{\CCINF(\LF\inv)} V\cong
  \Sch(\R)_I\Lhot_{\CCINF(\R)} V\) for all~\(V\).  Using also
  Lemma~\ref{lem:weighted_Sch_Lhot} we compute
  \begin{multline*}
    \Sch(\LF\inv)_+\Lhot_{\CCINF(\LF\inv)} \Sch(\LF) \cong
    \Sch(\R)_+\Lhot_{\CCINF(\R)} \Sch(\LF)
    \\ \cong
    \Sch(\R)_+\Lhot_{\CCINF(\R)}
    [\Sch(\R)_+   \hot\Sch(\LF\one)\oplus
     \Sch(\R)_-   \hot\Sch(\LF\one)
     \overset{(\iota, -\EFou\iota)}\longrightarrow
     \Sch(\R)_\cap\hot\Sch(\LF\one)
     ]
    \\ \cong
    [\Sch(\R)_+   \hot\Sch(\LF\one)\oplus
     \Sch(\R)_\cap\hot\Sch(\LF\one)
     \overset{(\iota, -\EFou)}\longrightarrow
     \Sch(\R)_\cap\hot\Sch(\LF\one)
    ]
    \\ \cong
    [\Sch(\LF\inv)_+\oplus\Sch(\LF\inv)_\cap
     \overset{(\iota, -\EFou)}\longrightarrow
    \Sch(\LF\inv)_\cap],
  \end{multline*}
  where ``\(\cong\)'' denotes homotopy equivalence.
  Proposition~\ref{pro:Fourier_estimate} yields that~\(\EFou\) is an
  isomorphism on \(\Sch(\LF\inv)_\cap\).  Hence we can replace the last
  chain complex simply by \(\Sch(\LF\inv)_+\).  The proof in the finite
  case is similar.
\end{proof}

Now we consider a finite set of places~\(S\) of a global field~\(\GF\).
Then \(\Sch(\Adel_S)\) is the completed tensor product of the spaces
\(\Sch(\GF_v)\) for \(v\in S\) and similarly for \(\CCINF(\Adel\inv_S)\).
We abbreviate \(G\defeq \Adel\inv_S\).

\begin{lemma}  \label{lem:S_local_temperation}
  Let \(I\subseteq \ooival{0,\infty}\), then we have a natural
  isomorphism
  \[
  \Sch(G)_I\Lhot_{\CCINF(G)} \Sch(\Adel_S) \cong \Sch(G)_I.
  \]
  If~\(\GF\) is a global function field, then also
  \[
  \Sch_\comp(G)_I\Lhot_{\CCINF(G)} \Sch(\Adel_S) \cong
  \Sch_\comp(G)_I.
  \]
\end{lemma}

\begin{proof}
  Let \(A(G)\) be the tensor product of the algebras \(\Sch(\GF\inv_v)_+\)
  for infinite places \(v\in S\) and \(\Sch_\comp(\GF\inv_v)_+\) for
  finite places \(v\in S\).  Thus \(A(G)\) is a smooth convolution algebra
  on~\(G\).  Since \(I\subseteq I_+\), the convolution algebra \(A(G)\) is
  contained in \(\Sch(G)_I\).  Thus \(\Sch(G)_I\) is an essential module
  over \(A(G)\).  Therefore, \(\Sch(G)_I \hot_{\CCINF(G)} W \cong
  \Sch(G)_I \hot_{A(G)} A(G) \hot_{\CCINF(G)} W\) for all essential
  \(\CCINF(G)\)\nbd{}modules~\(W\).  The \(A(G)\)\brd{}temperation functor
  evidently maps free modules again to free modules.  Therefore, we
  may replace~\(\hot\) by~\(\Lhot\) in the above statement.  Equations
  \eqref{eq:Sch_local_reduce_infinite}
  and~\eqref{eq:Sch_local_reduce_finite} imply \(A(G) \Lhot_{\CCINF(G)}
  \Sch(\Adel_S) \cong A(G)\).  Since this is a free essential
  \(A(G)\)\brd{}module, we get
  \begin{multline*}
    \Sch(G)_I \Lhot_{\CCINF(G)} \Sch(\Adel_S)
    \cong
    \Sch(G)_I \Lhot_{A(G)} A(G) \Lhot_{\CCINF(G)} \Sch(\Adel_S)
    \\ \cong
    \Sch(G)_I \Lhot_{A(G)} A(G)
    \cong
    \Sch(G)_I.
  \end{multline*}
  The global function field case is similar.
\end{proof}

\begin{theorem}  \label{the:Sch_Adel_temperation}
  There are natural isomorphisms
  \[
  \Sch(G\one) \Lhot_{\CCINF(G\one)} \Sch(\Adel_S) \cong \Twist(G)
  \]
  for algebraic number fields and
  \[
  \Sch(G\one) \Lhot_{\CCINF(G\one)} \Sch(\Adel_S) \cong
  \Twist_\comp(G\one)
  \]
  for global function fields.  These are isomorphisms in
  \(\Der(\CCINF(G))\).
\end{theorem}

\begin{proof}
  We only write down the proof for algebraic number fields.
  Lemma~\ref{lem:S_local_temperation} shows
  \(\Sch(G)_I\Lhot_{\CCINF(G)} \Sch(\Adel_S)\cong\Sch(G)_I\) for
  \(I\subseteq\ooival{0,\infty}\).  Using the Fourier transform and~\(J\),
  we obtain such an isomorphism also for
  \(I\subseteq\ooival{-\infty,1}\), compare
  Proposition~\ref{pro:Fourier_estimate}.  The linearly split
  extension
  \[
  \Sch(G)_\cup \into \Sch(G)_+\oplus \Sch(G)_- \prto \Sch(G)_\cap
  \]
  of Lemma~\ref{lem:weighted_Sch_extension} gives rise to the
  following distinguished triangle in \(\Der_0\):
  \begin{multline*}
    \Sch(G)_\cup \Lhot_{\CCINF(G)} \Sch(\Adel_S) \longrightarrow
    \Sch(G)_+    \Lhot_{\CCINF(G)} \Sch(\Adel_S) \oplus
    \Sch(G)_-    \Lhot_{\CCINF(G)} \Sch(\Adel_S) \\ \longrightarrow
    \Sch(G)_\cap \Lhot_{\CCINF(G)} \Sch(\Adel_S) \longrightarrow^{[1]}
    \Sch(G)_\cup \Lhot_{\CCINF(G)} \Sch(\Adel_S).
  \end{multline*}
  Plugging in our results for \(\Sch(G)_I\Lhot_{\CCINF(G)}
  \Sch(\Adel_S)\) for \(I_+\), \(I_-\) and \(I_\cap\), we obtain a
  distinguished triangle
  \begin{multline*}
    \Sch(G)_\cup \Lhot_{\CCINF(G)} \Sch(\Adel_S) \longrightarrow
    \Sch(G)_+ \oplus
    \Sch(G)_- \overset{(\iota,-\EFou\iota)}\longrightarrow
    \Sch(G)_\cap \\ \longrightarrow^{[1]}
    \Sch(G)_\cup \Lhot_{\CCINF(G)} \Sch(\Adel_S).
  \end{multline*}
  Lemma~\ref{lem:Twist_extension} asserts that \((\iota,-\EFou\iota)\)
  is a split surjection with kernel \(\Twist(G)\).  Hence \(\Sch(G)_\cup
  \Lhot_{\CCINF(G)} \Sch(\Adel_S)\cong\Twist(G)\).
  
  Let \(v\in S\) be an infinite place and split the norm homomorphism
  \(G\to\R\inv_+\) by \(\R\inv_+\to\GF\inv_v\to\Adel\inv_S\).  This
  induces an isomorphism
  \(\Sch(G)_\cup\cong\Sch(G\one)\hot\Sch(\R\inv_+)_\cup\).  The action
  of~\(\R\inv_+\) only involves the place~\(v\).  It is easy to see that
  \(\Sch(\GF_v)\cong\Twist(\GF\inv_v)\) and hence \(\Sch(\Adel_S)\) are
  essential modules over \(\Sch(\R\inv_+)_\cup\).  Therefore,
  Corollary~\ref{cor:weighted_Sch_isocohomological} yields
  \begin{multline*}
    \Twist(G) \cong
    \Sch(G)_\cup \Lhot_{\CCINF(G)} \Sch(\Adel_S)
    \cong
    \Sch(G\one) \Lhot_{\CCINF(G\one)}
    (\Sch(\R\inv_+)_\cup \Lhot_{\CCINF(\R\inv_+)} \Sch(\Adel_S))
    \\ \cong
    \Sch(G\one) \Lhot_{\CCINF(G\one)}
    (\Sch(\R\inv_+)_\cup \Lhot_{\Sch(\R\inv_+)_\cup} \Sch(\Adel_S))
    \cong
    \Sch(G\one) \Lhot_{\CCINF(G\one)} \Sch(\Adel_S).
  \end{multline*}
  Finally, the isomorphism \(\Sch(G\one)\Lhot_{\CCINF(G\one)}
  \Sch(\Adel_S) \cong\Twist(G)\) is just the composition
  \begin{displaymath}
    \Sch(G\one)\Lhot_{\CCINF(G\one)} \Sch(\Adel_S) \congto
    \Sch(G\one)\hot_{\CCINF(G\one)} \Sch(\Adel_S)
    \overset{\boin}\to
    \Sch(G\one)\hot_{\CCINF(G\one)} \Twist(G) \congto
    \Twist(G).
  \end{displaymath}
  All these (chain) maps are \(G\)\nbd{}equivariant with respect to the
  natural representations of~\(G\).  Hence we have an isomorphism in the
  derived category \(\Der(\CCINF(G))\).
\end{proof}

Suppose that~\(S\) contains all infinite places of~\(\GF\) so that
\(H\defeq\GF\inv_S\) is defined.  Let \(\ICL_S\defeq \Adel\inv_S/H\).
The group cohomology of~\(\GF\inv_S\) with coefficients in a
representation~\(V\) of~\(\Adel\inv_S\) is the homology of the chain
complex
\begin{equation}
  \label{eq:subgroup_homology}
  \C(1) \Lhot_{\CCINF(H)} V \cong
  \CCINF(\ICL_S) \Lhot_{\CCINF(G)} V,
\end{equation}
where \(\C(1)\) denotes the trivial representation of~\(\GF\inv_S\)
(see~\cite{Meyer:Smooth}).  Since \(\CCINF(\ICL_S)\) is a bimodule
over \(\CCINF(\ICL_S)\), the right hand side belongs to
\(\Der(\CCINF(\ICL_S))\), that is, it is a chain complex of smooth
representations of~\(\ICL_S\) in a natural way.

\begin{corollary}  \label{cor:Sch_local_coinvariants}
  \(\C(1) \Lhot_{\CCINF(H)} \Sch(\Adel_S)\) is isomorphic in
  \(\Der(\CCINF(\ICL_\GF))\) to \(\Twist(\ICL_S)\) for algebraic number
  fields and to \(\Twist_\comp(\ICL_S)\) for global function fields.
\end{corollary}

\begin{proof}
  We only write down the proof for algebraic number fields.
  Since~\(\ICL\one_S\) is compact, the space \(\CCINF(\ICL\one_S)\) is an
  essential module over \(\Sch(G\one)\).  Hence
  \begin{multline*}
    \C(1) \Lhot_{\CCINF(H)} \Sch(\Adel_S)
    \cong
    \CCINF(\ICL\one_S) \Lhot_{\CCINF(G\one)} \Sch(\Adel_S)
    \\ \cong
    \CCINF(\ICL\one_S) \Lhot_{\Sch(G\one)} \Sch(G\one)
    \Lhot_{\CCINF(G\one)} \Sch(\Adel_S)
    \cong
    \CCINF(\ICL\one_S) \Lhot_{\Sch(G\one)} \Twist(G).
  \end{multline*}
  We already computed \(\CCINF(\ICL\one_S)\hot_{\Sch(G\one)} \Twist(G)\)
  in Section~\ref{sec:equivariant_Fourier}.  The argument there also
  applies to~\(\Lhot\).
\end{proof}

Taking the homology of this chain complex, we obtain
\(H_n(\GF\inv_S,\Sch(\Adel_S))\cong0\) for \(n\ge1\) and
\(H_0(\GF\inv_S,\Sch(\Adel_S))\cong\Twist(\ICL_S)\) or
\(\Twist_\comp(\ICL_S)\), respectively.

\subsection{Schwartz functions on finite local fields}
\label{sec:Sch_local_finite}

In order to pass to the global case, we need finer information about
\(\Sch(\LF)\) for a finite local field.  Let
\[
\MCR\defeq \{x\in\LF\mid \abs{x}_\LF\le1\}
\]
be the maximal compact subring of~\(\LF\) and let \(\xi_0\defeq 1_\MCR
\in\Sch(\LF)\) be the characteristic function of~\(\MCR\).  Let
\[
\MCR\inv=\LF\one\defeq \{x\in\LF\inv\mid \abs{x}_\LF=1\}
\]
be the maximal compact subgroup of~\(\LF\one\).  Let
\(\Sch(\LF)^{\MCR\inv}\subseteq\Sch(\LF)\) be the subspace of all
functions that are invariant under~\(\MCR\inv\).  It consists of
functions that only depend on \(\abs{x}_\LF\).  An element
\(\varpi\in\MCR\) is called a \emph{uniformizer} if
\(\abs{\varpi}_\LF=q_\LF^{-1}\).  Hence~\(\varpi\) generates the maximal
ideal \(\MI\subseteq\MCR\) and \(\LF\inv\cong\MCR\inv\times\varpi^\Z\).

\begin{lemma}  \label{lem:Sch_local_finite}
  Let~\(\LF\) be a finite local field.
  \begin{enumerate}[(1)]
  \item The functions \(\xi_n\defeq
    \lambda_\varpi^n(\xi_0)=1_{\varpi^n\MCR}\) for \(n\in\Z\) form a
    basis of \(\Sch(\LF)^{\MCR\inv}\).  That is, \(\Sch(\LF)^{\MCR\inv}\)
    is the free \(\C[\Z]\)\nbd{}module with generator~\(\xi_0\).

  \item There is a \(\LF\inv\)\nbd{}equivariant isomorphism
    \(\Sch(\LF)\cong\CCINF(\LF\inv)\) that sends~\(1_{\MCR}\)
    to~\(1_{\MCR\inv}\).

  \end{enumerate}
\end{lemma}

\begin{proof}
  We claim that the functions \((\xi_j)\) are linearly independent.
  Consider a finite linear combination \(f=\sum c_j\xi_j\).  Unless
  all~\(c_j\) vanish, there is a smallest~\(j\) with \(c_j\neq0\).  Then
  \(f(x)=c_j\neq0\) for \(\abs{x}=q_\LF^{-j}\), which is the desired
  linear independence.  Next we claim that \(\Sch(\LF)^{\MCR\inv}\) is
  spanned by the elements \((\xi_j)\).  Let \(f\in\Sch(\LF)^{\MCR\inv}\).
  Since~\(f\) is locally constant and compactly supported, there is
  \(a\in\N\) such that~\(f\) is constant on \(\varpi^a\MCR\)\brd{}cosets.
  There is a smallest \(n\in\Z\) such that \(f|_{\varpi^n\MCR\inv}\neq0\).
  Subtracting the corresponding multiple of~\(\xi_n\), we get a new
  function~\(f^{(n-1)}\) that is supported in \(\varpi^{n-1}\MCR\).
  Moreover, since~\(n\) is maximal, we must have \(n\le a\), so
  that~\(f^{(n-1)}\) is also constant on \(\varpi^a\MCR\)\brd{}cosets.
  Proceeding in this fashion, we obtain \(c_j\in\C\) such that
  \(f-\sum_{j=n}^a c_j\xi_j\) is supported in \(\varpi^{a+1}\MCR\) and
  constant on \(\varpi^a\MCR\)\brd{}cosets.  Hence this difference
  vanishes, proving the claim.
  
  Thus the functions \((\xi_j)\) form a basis for
  \(\Sch(\LF)^{\MCR\inv}\).  Therefore, there is a unique
  \(\LF\inv\)\nbd{}equivariant isomorphism
  \(\Sch(\LF)^{\MCR\inv}\cong\CCINF(\LF\inv)^{\MCR\inv}\) that
  sends~\(1_\MCR\) to \(1_{\MCR\inv}\).  It is straightforward to see that
  a function in \(\Sch(\LF)\) that vanishes at~\(0\) already lies in
  \(\CCINF(\LF\inv)\).  If \(\chi\in\widehat{\MCR\inv}\) is a non-trivial
  character, then no \(\chi\)\nbd{}homogeneous function can vanish
  at~\(0\).  Hence the \(\chi\)\nbd{}homogeneous subspaces of \(\Sch(\LF)\)
  and \(\CCINF(\LF\inv)\) coincide.  Putting all homogeneous subspaces
  together, we get the assertion.
\end{proof}

\subsection{A maximum principle}
\label{sec:maximum_principle}

In the global case, we have to control things that happen in the
critical strip from outside.  The technical results of this section
are needed for this purpose.  We need some generalities about
multipliers of smooth convolution algebras.  Since~\(\Adel\inv_S\) is
commutative, multipliers of \(\CCINF(\Adel\inv_S)\) are just
\(\Adel\inv_S\)\brd{}equivariant bounded linear operators on
\(\CCINF(\Adel\inv_S)\).  Multipliers act on any essential module over
\(\CCINF(\Adel\inv_S)\), that is, on any smooth representation
of~\(\Adel\inv_S\) by~\cite{Meyer:Smooth}.  As we observed in
Section~\ref{sec:homology_smooth_convolution}, we have
\[
\Hom_{A(G)}(A(G),A(G)) \cong \Hom_{\CCINF(G)}(\CCINF(G),A(G))
\]
for any smooth convolution algebra \(A(G)\).  This implies that a
distribution \(h\in\CCINF'(G)\) is a multiplier of \(A(G)\) if and only
if~\(h\) multiplies \(\CCINF(G)\) into \(A(G)\) both on the left and right.
Hence \(A_1(G)\subseteq A_2(G)\) implies
\(\Mult(A_1(G))\subseteq\Mult(A_2(G))\).

\begin{lemma}  \label{lem:auxiliary_multiplier}
  For any \(\alpha>0\), there exists a multiplier
  \(h\in\Mult(\CCINF(\Adel\inv_S))\) that is invertible on
  \(\Sch(\Adel\inv_S)_+\) and restricts to a bounded map from
  \(\Sch(\Adel_S)\) to \(\Sch(\Adel\inv_S)_{\ooival{1-\alpha,\infty}}\).
\end{lemma}

\begin{proof}
  Since \(\bigotimes_{v\in S} \Sch(\GF\inv_v)_I\boin
  \Sch(\Adel\inv_S)_I\) for any interval~\(I\), we get a stronger
  assertion if we replace \(\Sch(\Adel\inv_S)_I\) by this tensor product
  in the statement of the lemma.  Since \(\CCINF(\Adel\inv_S)\) and
  \(\Sch(\Adel_S)\) are completed tensor products of the corresponding
  spaces for~\(\GF_v\), \(v\in S\), we can get the multiplier
  for~\(\Adel_S\) by tensoring together corresponding multipliers for
  the local fields~\(\GF_v\).  Therefore, we only have to construct
  \(h\in\Mult(\CCINF(\LF\inv))\) for a local field~\(\LF\) that induces an
  invertible operator on \(\Sch(\LF\inv)_+\) and a bounded operator
  \(\Sch(\LF)\to\Sch(\LF\inv)_{\ooival{1-\alpha,\infty}}\).

  For finite~\(\LF\), Lemma~\ref{lem:Sch_local_finite} describes an
  isomorphism \(\Sch(\LF)\cong\CCINF(\LF\inv)\).  Restricting it to
  \(\CCINF(\LF\inv)\subseteq\Sch(\LF)\), we obtain a multiplier~\(h\) of
  \(\CCINF(\LF\inv)\).  It evidently has the required properties for
  any~\(\alpha\).  Let \(\LF=\R\).  The differential operator
  \(x^n(d/dx)^n\) on \(\Sch(\R)\) is dilation invariant and hence
  restricts to an equivariant differential operator on~\(\R\inv\).  This
  is a multiplier~\(h_n\) of \(\CCINF(\R\inv)\).  It is a bounded map from
  \(\Sch(\R)\) into the subspace of Schwartz functions that have a zero
  of order at least~\(n\) at the origin.  The latter space is evidently
  contained in \(\Sch(\R\inv)_{\ooival{-n,\infty}}\).  It is easy to
  write~\(h_n\) as a polynomial function of \(h_1=x\cdot d/dx= d/d\ln x\).
  The spectrum of \(d/d\ln x\) on \(\Sch(\R\inv)_{\ooival{0,\infty}}\) is
  contained in the region \(\{s\in\C\mid \RE s>0\}\).  This implies
  that~\(h_n\) is invertible on \(\Sch(\R\inv)_{\ooival{0,\infty}}\).  For
  \(\LF=\C\), a similar argument shows that the differential operators
  \((x\conj{x})^n (\partial^2/\partial x\,\partial\conj{x})^n\) work.
\end{proof}

\begin{lemma}  \label{lem:maximum_principle}
  For any \(\alpha\ge0\), the map
  \[
  \Twist(\ICL_S) \to
  \Sch(\ICL_S)_{\ooival{\alpha,\infty}} \cap
  \EFou(\Sch(\ICL_S)_{\ooival{-\infty,1-\alpha}}),
  \qquad
  f\mapsto (f,f),
  \]
  is a bornological embedding.
\end{lemma}

\begin{proof}
  Let \(\beta>\alpha\) and let~\(T\) be a subset of \(\Twist(\ICL_S)\) that
  is bounded in \(\Sch(\ICL_S)_{\coival{\beta,\infty}}\) and
  \(\EFou(\Sch(\ICL_S)_{\ocival{-\infty,1-\beta}})\).  We claim that
  then~\(T\) is bounded in \(\Twist(\ICL_S)\).  This claim establishes the
  lemma.  Since
  \[
  \Sch(\ICL_S)_{[1/2,\beta]}=\Sch(\ICL_S)_{1/2}\cap \Sch(\ICL_S)_\beta,
  \qquad
  \EFou(\Sch(\ICL_S)_{[1-\beta,1/2]})=
  \EFou(\Sch(\ICL_S)_{1/2})\cap \EFou(\Sch(\ICL_S)_\beta)
  \]
  and~\(\EFou\) is a bornological isomorphism on \(\Sch(\ICL_S)_{1/2}\),
  it actually suffices to prove that~\(T\) is bounded in
  \(\Sch(\ICL_S)_{1/2}\).

  Let~\(h\) be the multiplier of Lemma~\ref{lem:auxiliary_multiplier}.
  The same argument as in the proof of
  Proposition~\ref{pro:Fourier_estimate} shows that \(h\circ\EFou\) on
  \(L^2(\Adel\inv_S)_{1/2}\) extends to a bounded linear operator on
  \(\Sch(\Adel\inv_S)_I\) for any interval \(I\subseteq
  \ooival{1-\alpha,\infty}\).  Furthermore, \(h\circ\EFou\) descends to a
  bounded linear operator on \(\Sch(\ICL_S)_I =
  \CCINF(\ICL\one_S)\hot_{\Sch(\Adel\one_S)} \Sch(\Adel\inv_S)_I\).
  For the same reason, \(h\) descends to a bounded multiplier of
  \(\CCINF(\ICL_S)\) and hence acts as a bounded operator on any
  representation of~\(\ICL_S\).  By hypothesis, \(\EFoi(T)\) is bounded in
  \(\Sch(\ICL_S)_{1-\beta}\) and~\(T\) is bounded in \(\Sch(\ICL_S)_\beta\).
  Hence \(h\EFou\bigl(\EFoi(T)\bigr)=h(T)\) is bounded in
  \[
  \Sch(\ICL_S)_{1-\beta} \cap \Sch(\ICL_S)_\beta
  =
  \Sch(\ICL_S)_{[1-\beta,\beta]}
  \boin
  \Sch(\ICL_S)_{1/2}.
  \]
  By construction, \(h\) acts as a bornological isomorphism on
  \(\Sch(\Adel\inv_S)_{1/2}\) and hence on \(\Sch(\ICL_S)_{1/2}\).
  Thus~\(T\) is bounded in \(\Sch(\ICL_S)_{1/2}\), which is exactly what
  we need.
\end{proof}

\section{The global difference representation}
\label{sec:global_difference}

In this section, we construct the global difference representation and
prove some basic facts about it.  We show that its character is equal
to the Weil distribution and that its spectral multiplicity function
is equal to the pole order function for the \(L\)\nbd{}function.  This
yields André Weil's Explicit Formula.  We construct some further
algebraic structure on the global difference representation in
Section~\ref{sec:symmetries} and briefly mention in
Section~\ref{sec:global_spectral_analysis} how our approach contains
the meromorphic continuation of \(L\)\nbd{}functions and their
functional equations.

\subsection{More on adeles and ideles}
\label{sec:global_adeles}

We already recalled the basic conventions regarding adeles and ideles
in the \(S\)\nbd{}local case in
Section~\ref{sec:local_trace_conventions}.  Here we add some more
remarks about the global case.  The main nuance between the
\(S\)\nbd{}local and global case is that the adele ring~\(\Adel_\GF\) and
the idele group~\(\Adel\inv_\GF\) are defined as \emph{restricted direct
  products}.  This is necessary to get locally compact groups.  Taking
this into account, we can define the same kind of structure as in the
\(S\)\nbd{}local case.

The field~\(\GF\) embeds diagonally into~\(\Adel_\GF\) and this is an
isomorphism onto a discrete and cocompact subfield of~\(\Adel_\GF\).  We
can again define a norm homomorphism on~\(\Adel\inv_\GF\) and
let~\(\Adel\one_\GF\) be its kernel.  The group~\(\GF\inv\) is a discrete
and cocompact subgroup of~\(\Adel\one_\GF\).  Equivalently, the
norm-\(1\)\brd{}subgroup~\(\ICL\one_\GF\) of the idele class group
\(\ICL_\GF\defeq \Adel\inv_\GF/\GF\inv\) is compact.

Since~\(\Adel\inv_\GF\) is the restricted direct product of the
groups~\(\GF\inv_v\), the space \(\CCINF(\Adel\inv_\GF)\) is the
restricted completed tensor product of the spaces \(\CCINF(\GF\inv_v)\)
with distinguished vectors \(1_{\MCR\inv}\in \CCINF(\GF\inv_v)\) for
finite places (see~\cite{Meyer:Smooth}).  The space
\(\Sch(\Adel_\GF)\) can be described similarly: it is the restricted
completed tensor product of the spaces \(\Sch(\GF_v)\) with
distinguished vectors \(1_{\MCR}\in\Sch(\GF_v)\) for finite places.

Let~\(S\) be a finite set of places containing all infinite places and
let \(\Cont{S}\defeq\Places(\GF)\setminus S\).  Define \(\GF_S\)
and~\(\GF\inv_S\) as in Section~\ref{sec:local_trace_conventions}.  Let
\[
\MCR_{\Cont{S}} \defeq \prod_{v\in\Cont{S}} \MCR_v,
\qquad
\MCR\inv_{\Cont{S}} \defeq \prod_{v\in\Cont{S}} \MCR\inv_v.
\]
These are the maximal compact subring of \(\Adel_{\Cont{S}}\) and the
maximal compact subgroup of \(\Adel\inv_{\Cont{S}}\), respectively.  We
also view \(\MCR\inv_{\Cont{S}}\) as a compact subgroup of~\(\ICL_\GF\).
There is an obvious injective group homomorphism
\(\ICL_\GF/\MCR\inv_{\Cont{S}}\to\ICL_S\).  We call a finite set of
places~\(S\) that contains all infinite places \emph{sufficiently large}
if this map is an isomorphism
\(\ICL_S\cong\ICL_\GF/\MCR\inv_{\Cont{S}}\).  Equivalently,
\(\Adel\inv_{\Cont{S}}/\MCR\inv_{\Cont{S}} \cong \GF\inv/\GF\inv_S\).
The finiteness of the ideal class group implies that there exists a
finite set of places~\(S_0\) such that all \(S\supseteq S_0\) are
sufficiently large in this sense.

Given a smooth representation \(\pi\colon \ICL_\GF\to\Aut(V)\), we let
\(V^S\subseteq V\) be the subspace of elements that are invariant under
\(\MCR\inv_{\Cont{S}}\) and call elements of~\(V^S\) \emph{unramified
  outside~\(S\)}.  The subgroups \(\MCR\inv_{\Cont{S}}\) for the
sufficiently large finite sets of places of~\(\GF\) form a fundamental
system of smooth compact subgroups.  Hence \(V=\varinjlim V^S\) for any
smooth representation.

We normalize the Haar measure on~\(\ICL_\GF\) as in the \(S\)\nbd{}local
case.  Since \(\GF\subseteq\Adel_\GF\) is discrete and cocompact, there
is \(\psi\in\widehat{\Adel_\GF}\) with \(\psi|_\GF=1\) and \(\psi\neq1\).
We call such characters \emph{normalized}.  This normalization insures
that the \emph{dual lattice}
\[
\GF^\bot \defeq
\{y\in\Adel_\GF\mid
\text{\(\psi(ya)=0\) for all \(a\in\GF\)}\}
\]
is equal to~\(\GF\) again, which is useful for the Poisson Summation
Formula.

Let~\(\psi\) be a normalized character.  We can write it as
\(\psi=\psi_0(\different\cdot x)\), where the local factors
\(\psi_{0,v}\in\widehat{\GF_v}\) of~\(\psi_0\) are normalized for all
places~\(v\) and where \(\different\in\Adel\inv_\GF\).  The
idele~\(\different\) is called a \emph{differental idele}.  It is easy
to see that~\(\psi\) is unique up to the action of~\(\GF\inv\).  Hence the
class of \(\different\) in~\(\ICL_\GF\) is independent of the choice
of~\(\psi\).  If~\(S\) is a sufficiently large set of places, we can
choose~\(\psi\) unramified outside~\(S\), that is, such that~\(\psi_v\) is
normalized for \(v\in\Cont{S}\).  We always assume~\(\psi\) unramified
outside~\(S\) if this is relevant.  The norm \(\abs{\different}^{-1}\) is
an important geometric invariant of the global field~\(\GF\).  For an
algebraic number field, \(\abs{\different}^{-1}\) is the
\emph{discriminant} of~\(\GF\).  For a global function field,
\(\abs{\different}^{-1}=g_\GF^{2g-2}\), where~\(g\) is the \emph{genus} of
the underlying curve of~\(\GF\).  See~\cite{Weil:Basic}*{page~113} for
this.

As in the \(S\)\nbd{}local case, we use~\(\psi\) to identify~\(\Adel_\GF\)
with its Pontrjagin dual and to view the Fourier transform as an
operator on \(\Sch(\Adel_\GF)\).  We normalize the Haar measure
on~\(\Adel_\GF\) so that the Fourier transform~\(\Fourier_\psi\) defined
by~\eqref{eq:def_Fourier} for normalized~\(\psi\) is unitary.  This Haar
measure on~\(\Adel_\GF\) is called \emph{self-dual}.  It can be
characterized by the condition that the compact group \(\Adel_\GF/\GF\)
has volume~\(1\).

An important difference between the global and the \(S\)\nbd{}local case
is that the Hilbert spaces \(L^2(\Adel_\GF,dx)\) and
\(L^2(\Adel\inv_\GF,\abs{x}\,d\inv x)\) become different.  The problem
is that~\(dx\) gives zero measure to the set of \(x\in\Adel_\GF\) with
\(\abs{x}\neq0\).  There are several other related problems that make it
difficult to analyze the global situation using Hilbert spaces.  Our
Fréchet space setup is not affected by these problems, however.

\subsection{The coinvariant space \(\Sch(\Adel_\GF)/\GF\inv\)}
\label{sec:Sch_global_coinvariants}

We are going to describe the coinvariant space \(\Hilm_+\defeq
\Sch(\Adel_\GF)/\GF\inv\) more explicitly.  Let~\(S\) be a sufficiently
large finite set of places containing all infinite places.  Let
\(\Sch(\Adel_\GF)^S\subseteq \Sch(\Adel_\GF)\) be the subspace of
\(\MCR\inv_{\Cont{S}}\)\brd{}invariants as above.

\begin{theorem}  \label{the:global_Sch_coinvariants}
  If~\(S\) is sufficiently large, then there are natural isomorphisms
  \[
  \C(1) \Lhot_{\CCINF(\GF\inv)} \Sch(\Adel_\GF)^S \cong
  \C(1) \Lhot_{\CCINF(\GF\inv_S)} \Sch(\Adel_S) \cong
  \Twist(\ICL_S)
  \]
  in the derived category \(\Der(\CCINF(\ICL_\GF))\) of smooth
  representations of~\(\ICL_\GF\).
\end{theorem}

\begin{proof}
  The space \(\Sch(\Adel_{\Cont{S}})\) is the restricted tensor product
  of the spaces \(\Sch(\GF_v)\) with distinguished vectors \(1_{\MCR_v}\)
  for \(v\in\Cont{S}\).  Hence
  \(\Sch(\Adel_{\Cont{S}})^{\MCR\inv_{\Cont{S}}}\) is the restricted
  tensor product of the spaces \(\Sch(\GF_v)^{\MCR\inv_v}\).
  Lemma~\ref{lem:Sch_local_finite} yields an equivariant isomorphism
  \(\Sch(\GF_v)^{\MCR\inv_v}\cong \CCINF(\GF\inv_v/\MCR\inv_v)\) that
  sends \(1_{\MCR_v}\) to \(1_{\MCR\inv_v}\).  These local isomorphisms
  combine to an isomorphism
  \[
  \Sch(\Adel_{\Cont{S}})^{\MCR\inv_{\Cont{S}}}
  \cong
  \CCINF(\Adel\inv_{\Cont{S}}/\MCR\inv_{\Cont{S}})
  \cong
  \C[\GF\inv/\GF\inv_S]
  \]
  because~\(S\) is sufficiently large.  Hence
  \begin{displaymath}
    \Sch(\Adel_\GF)^S
    \cong
    \Sch(\Adel_S) \hot \Sch(\Adel_{\Cont{S}})^{\MCR\inv_{\Cont{S}}}
    \cong
    \Sch(\Adel_S) \hot \C[\GF\inv/\GF\inv_S].
  \end{displaymath}
  The assertions now follows from
  Corollary~\ref{cor:Sch_local_coinvariants} and the Shapiro Lemma.
\end{proof}

The subspace of \(\MCR\inv_{\Cont{S}}\)\brd{}invariant elements of a
representation is the range of a projection~\(Q_S\) in the multiplier
algebra of \(\CCINF(\ICL_\GF)\).  Hence~\(Q_S\) acts on any chain complex
of smooth representations of~\(\ICL_\GF\).
Theorem~\ref{the:global_Sch_coinvariants} asserts
\[
Q_S(\C(1)\Lhot_{\CCINF(\ICL_\GF)} \Sch(\Adel_\GF))
\cong
\C(1)\Lhot_{\CCINF(\ICL_\GF)} \Sch(\Adel_\GF)^S
\cong
\Twist(\ICL_S).
\]
If \((S_n)_{n\in\N}\) is an increasing sequence of sufficiently large
finite sets of places with \(\bigcup S_n=\Places(\GF)\), then
\(\sum_{n\in\N} Q_{S_n}-Q_{S_{n-1}} = \ID\) on
\(\Der(\CCINF(\ICL_\GF))\).  Hence
\[
\C(1)\Lhot_{\CCINF(\ICL_\GF)} \Sch(\Adel_\GF)
\cong
\varinjlim \Twist(\ICL_S),
\]
that is, \(\Sch(\Adel_\GF)\) is acyclic for the coinvariant space
functor and
\[
\Hilm_+\defeq \Sch(\Adel_\GF)/\GF\inv
\cong
\varinjlim \Twist(\ICL_S).
\]
Moreover, the structure maps in this system are bornological
embeddings because
\[
\Hilm_+^S \defeq
\Hilm_+^{\MCR\inv_{\Cont{S}}} =
Q_S(\Hilm_+) \cong
\Twist(\ICL_S).
\]
It is remarkable that the subspaces~\(\Hilm_+^S\) do not depend on
the global structure of~\(\Adel_\GF\) but only on things happening
inside the finite set of places~\(S\).  This allows us to reduce global
computations to \(S\)\nbd{}local ones.

It remains to describe how the embedding \(\Hilm_+^S\to\Hilm_+^{S'}\)
for \(S\subseteq S'\) looks like as a map
\(\Twist(\ICL_S)\to\Twist(\ICL_{S'})\).  We have identified
\(\Sch(\Adel_S)/\GF\inv_S\cong\Twist(\ICL_S)\).  It is easy to see that
the isomorphism is given by the \emph{\(S\)\nbd{}local summation map}
\[
\SUM_S\colon \Sch(\Adel_S)/\GF\inv_S \congto \Twist(\ICL_S),
\qquad
f(x) \defeq \sum_{a\in\GF\inv_S} f(ax).
\]
Take \(f\in\Twist(\ICL_S)\cong\Sch(\Adel_S)/\GF\inv_S\) and represent
it by a function \(\bar{f}\in\Sch(\Adel_S)\).  The embedding
\(\Sch(\Adel_S)\to\Sch(\Adel_{S'})\) sends~\(\bar{f}\) to
\(1_{\MCR_{S'\setminus S}} \otimes\bar{f}\).  Thus the induced map
\(\Twist(\ICL_S)\to\Twist(\ICL_{S'})\) on coinvariants sends~\(f\) to
\(\SUM_{S'}(1_\MCR\otimes \bar{f})\).  To describe this more explicitly,
we choose for each \(v\in\Cont{S}\) an element \(p_v\in\GF_S\) with
\(\abs{p_v}_w=1\) for all \(w\in\Cont{S}\setminus\{v\}\) and
\(\abs{p_v}_v=q_v^{-1}\).  This is possible because~\(S\) is sufficiently
large.  We also view~\(p_v\) as an element of \(\Adel\inv_S\) by the usual
diagonal embedding \(\GF\inv\to\Adel\inv_S\).  We have
\[
\GF\inv_{S'}=\GF\inv_S\times \bigoplus_{v\in S'\setminus S} p_v^\Z.
\]
Represent \(\dot{x}\in\ICL_{S'}\) by \(x\in\Adel\inv_{S'}\) with
\(\abs{x}_v=1\) for all \(v\in S'\setminus S\).  This is possible
because~\(S\) is sufficiently large.  We have
\begin{multline*}
  \SUM_{S'} (1_\MCR\otimes \bar{f})(\dot{x}) =
  \sum_{a\in\GF\inv_{S'}} (1_\MCR\otimes\bar{f})(ax) =
  \sum_{a\in\GF\inv_S, n_v\ge0}
  \bar{f}(\prod_{v\in S'\setminus S} p_v^{n_v}\cdot a x)
  \\ =
  \sum_{n_v\ge0} \prod_{v\in S'\setminus S}
  \lambda_{p_v}^{-n_v} f(\dot{x})
  =
  \prod_{v\in S'\setminus S} \sum_{n_v=0}^\infty
  \lambda_{p_v}^{-n_v} f(\dot{x})
  =
  \prod_{v\in S'\setminus S} (\ID-\lambda_{p_v}^{-1})^{-1} f(\dot{x}).
\end{multline*}
That is, the structure map \(\Twist(\ICL_S)\to\Twist(\ICL_{S'})\) is
equal to \(\prod_{v\in S'\setminus S} (\ID-\lambda_{p_v}^{-1})^{-1}\)
composed with the identical embedding.

If we only wanted to construct some spectral interpretation for zeros
of \(L\)\nbd{}functions, then we could also take the spaces
\(\Twist(\ICL_S)\defeq\Hilm_+^S\) as our starting point and define
\(\Hilm_+\defeq \varinjlim \Twist(\ICL_S)\) with the above structure
maps.  This does produce \(\Sch(\Adel_\GF)/\GF\inv_S\) without the need
to compute coinvariant spaces and hence saves some work.  However, the
construction looks rather \emph{ad hoc}.  We prefer to arrive at this
construction starting from the natural object
\(\Sch(\Adel_\GF)/\GF\inv_S\).

\subsection{Some embeddings and projections}
\label{sec:summation}

We need some definitions.  Let \(I_\pp\defeq \ooival{1,\infty}\) and
\(I_\mm\defeq \ooival{-\infty,0}\).  Let \(\Sch(\ICL_\GF)_\pp\defeq
\Sch(\ICL_\GF)_{I_\pp}\), \(\Sch(\ICL_\GF)_\mm\defeq
\Sch(\ICL_\GF)_{I_\mm}\) and
\[
\Sch(\ICL_\GF)_\ppmm \defeq
\Sch(\ICL_\GF)_\pp\oplus \Sch(\ICL_\GF)_\mm.
\]
Finally, define
\[
L^2(\ICL_S)_\pp \defeq \bigcap_{\alpha>1} L^2(\ICL_S)_\alpha.
\]

\begin{lemma}  \label{lem:Sch_smoothen_Ltwo}
  \(\smooth_{\ICL_S} (L^2(\ICL_S)_\pp) = \Sch(\ICL_S)_\pp\).
\end{lemma}

\begin{proof}
  The norm homomorphism~\(\abs{x}\) is a proper map \(\ICL_S\to\R\inv_+\).
  Hence rapid decay on~\(\ICL_S\) is measured by the functions
  \((1+\abs{\ln{}\abs{x}})^\beta\) with \(\beta\in\R\).  These functions
  are dominated by \(\abs{x}^\epsilon+\abs{x}^{-\epsilon}\).  Hence
  \[
  L^2(\ICL_S)_\pp =
  \bigcap_{\alpha>1,\beta>0} L^2(\ICL_S,\abs{x}^{2\alpha}
  (1+\abs{\ln{}\abs{x}})^\beta \,d\inv x).
  \]
  It is easy to identify the smoothening of the latter
  representation with \(\Sch(\ICL_S)_\pp\).
\end{proof}

The \emph{global summation map} is defined by
\[
\SUM f(x) \defeq \sum_{a\in\GF\inv} f(ax)
\]
for \(f\in\Sch(\Adel_S)\), \(x\in\ICL_\GF\).

\begin{lemma}  \label{lem:summation}
  The summation map~\(\SUM\) defines a bounded linear map from~\(\Hilm_+\)
  to \(\Sch(\ICL_\GF)_\pp\).  For a sufficiently large finite set of
  places~\(S\), its restriction to \(\MCR\inv_{\Cont{S}}\)\brd{}invariants
  is isomorphic to the map
  \[
  \Lop_S \defeq \prod_{v\in \Cont{S}} (1-\lambda_{p_v}^{-1})^{-1}
  \colon \Twist(\ICL_S) \to \Sch(\ICL_S)_\pp.
  \]
  In addition, \(\Lop_S\) is a \(\ICL_S\)\nbd{}equivariant bornological
  isomorphism both on \(L^2(\ICL_S)_\pp\) and \(\Sch(\ICL_S)_\pp\) with
  inverse \(\Lopi_S \defeq \prod_{v\in \Cont{S}}
  (1-\lambda_{p_v}^{-1})\).
\end{lemma}

\begin{proof}
  We identify~\(\SUM\) with the map~\(\Lop_S\) on
  \(\Hilm_+^S\cong\Twist(\ICL_S)\) in the same way as we described the
  structure maps \(\Twist(\ICL_S)\to\Twist(\ICL_{S'})\) for \(S\subseteq
  S'\) at the end of Section~\ref{sec:Sch_global_coinvariants}.  It
  remains to check that \(\Lop_S\) and~\(\Lopi_S\) define bounded linear
  operators on \(\Sch(\ICL_S)_\pp\).
  
  Since \(1=\abs{a}=\abs{a}_S\cdot\abs{a}_{\Cont{S}}\) for all
  \(a\in\GF\inv\), we get \(\abs{p_v}_S=q_v\) and hence
  \[
  \norm{\lambda_{p_v}^{-1}}_{L^2(\ICL_S)_\alpha} = q_v^{-\alpha}
  \]
  for all \(v\in \Cont{S}\).  Therefore, the infinite products that
  define \(\Lop_S\) and~\(\Lopi_S\) converge absolutely in the operator
  norm on \(L^2(\ICL_S)_\alpha\) for any \(\alpha>1\).  This is equivalent
  to the convergence of the Euler product defining the
  \(\zeta\)\nbd{}function in the region \(\RE s>1\).  Thus \(\Lop_S\) and
  \(\Lopi_S\) define bounded linear operators on \(L^2(\ICL_S)_\pp\).
  Since smoothening is functorial, Lemma~\ref{lem:Sch_smoothen_Ltwo}
  also yields boundedness on \(\Sch(\ICL_S)_\pp\).
\end{proof}

\begin{lemma}  \label{lem:Hilm_p_embed}
  The map
  \[
  i_+ \colon \Hilm_+ \to \Sch(\ICL_\GF)_\ppmm,
  \qquad f \mapsto (\SUM f, J\SUM\Fourier^\ast f),
  \]
  is a \(\ICL_\GF\)\nbd{}equivariant bornological embedding.
  
  Let \(\phi\in\CCINF(\R_+)\) be a cut-off function as before and
  let~\(S\) be a sufficiently large finite set of places.  Then
  \begin{displaymath}
    p_+^S\colon \Sch(\ICL_S)_\ppmm \to \Twist(\ICL_S),
    \qquad
    (f_0,f_1)\mapsto
    M_{1-\phi}\Lopi_S f_0 + \Fourier M_{1-\phi} \Lopi_S J f_1.
  \end{displaymath}
  is a bounded \(\ICL\one_S\)\nbd{}equivariant linear operator.
\end{lemma}

\begin{proof}
  It suffices to check that~\(i_+\) restricts to an equivariant
  bornological embedding on~\(\Hilm_+^S\) for any sufficiently large
  finite set of places~\(S\).  We identify \(\Hilm_+^S\cong\Twist(\ICL_S)\)
  and the summation map with~\(\Lop_S\), which is a bornological
  isomorphism on \(\Sch(\ICL_S)_\pp\) by Lemma~\ref{lem:summation}.  The
  map \(J\colon \Sch(\ICL_S)_\pp\congto\Sch(\ICL_S)_\mm\) is a
  bornological isomorphism as well.  Hence the first assertion is
  equivalent to the statement that
  \[
  \Twist(\ICL_S)\to\Sch(\ICL_S)_\ppmm,
  \qquad
  f\mapsto (f, \EFoi f)
  \]
  is a bornological embedding.  This follows from
  Lemma~\ref{lem:maximum_principle}.  The map~\(\Lopi_S\) is an
  isomorphism on \(\Sch(\ICL_S)_\pp\) by Lemma~\ref{lem:summation},
  \(M_{1-\phi}\) is a bounded map from \(\Sch(\ICL_S)_\pp\) to
  \(\Sch(\ICL_S)_\cup\boin\Twist(\ICL_S)\) and~\(\Fourier\) is a bounded
  on \(\Twist(\ICL_S)\).  Thus~\(p_+^S\) is bounded.
\end{proof}

Hence we may identify~\(\Hilm_+\) with its image
\(i_+(\Hilm_+)\subseteq\Sch(\ICL_\GF)_\ppmm\).

We extend~\(p_+^S\) to an operator \(p_+^S\colon \Sch(\ICL_\GF)_\ppmm \to
\Hilm_+\) as follows.  We decompose all occurring spaces as direct sums
of \(\chi\)\nbd{}homogeneous subspaces for
\(\chi\in(\MCR\inv_{\Cont{S}})\sphat\).  For any such~\(\chi\), we
let~\(S(\chi)\) be the union of~\(S\) with the places where~\(\chi\) is
ramified.  We let \(p_+^S\defeq p_+^{S(\chi)}\) on the
\(\chi\)\nbd{}homogeneous subspace.

On~\(\Hilm_+^S\), we compute
\begin{multline*}
  p_+^S\circ i_+ =
  M_{1-\phi}\Lopi_S \Lop_S
  +\Fourier M_{1-\phi} \Lopi_S J^2 \Lop_S \Fourier^\ast
  \\ =
  M_{1-\phi} + \Fourier M_{1-\phi} \Fourier^\ast =
  \ID - (M_\phi - \EFou JM_{1-\phi}J \EFoi).
\end{multline*}
We call~\(\phi\) \emph{symmetric} if \(\phi(x)=1-\phi(x^{-1})\).  Then
\(JM_{1-\phi}J=M_\phi\) and hence
\[
p_+^S\circ i_+ = \ID - (M_\phi - \EFou M_\phi \EFoi).
\]
Thus \(\IN\lambda(f)\circ(\ID-p_+^Si_+)\) is a uniformly nuclear
operator on \(\Twist(\ICL_S)\) for \(f\in\CCINF(\ICL_S)\) and the local
trace formula~\eqref{eq:local_trace_formula_psi} computes its trace.
Therefore, we may think of~\(p_+^S\) as an approximate section for the
embedding~\(i_+\).  The other composition
\[
P_+^S \defeq
i_+\circ p_+^S =
\begin{pmatrix}
  \Lop_S M_{1-\phi} \Lopi_S & 
  \Lop_S \Fourier M_{1-\phi} \Lopi_S J \\
  J \Lop_S \Fourier^\ast M_{1-\phi} \Lopi_S & 
  J \Lop_S M_{1-\phi} \Lopi_S J
\end{pmatrix}
\]
is an approximate projection onto the range of~\(i_+\).

Our character computations compare~\(P_+^S\) with a more trivial
projection~\(P_-\), which we now define.  Let \(\Hilm_-\defeq
\Sch(\ICL_\GF)_\cup\).  We define maps
\begin{alignat*}{2}
  i_- &\colon \Hilm_- \to \Sch(\ICL_\GF)_\ppmm,
  &\qquad f &\mapsto (f,f),
  \\
  p_- &\colon \Sch(\ICL_\GF)_\ppmm \to \Hilm_-,
  &\qquad (f_0,f_1) &\mapsto M_{1-\phi} f_0+ M_\phi f_1.
\end{alignat*}
The map~\(i_-\) is a \(\ICL_\GF\)\nbd{}equivariant bornological embedding,
which we use to identify~\(\Hilm_-\) with a subspace of
\(\Sch(\ICL_\GF)_\ppmm\).  The operator~\(p_-\) is a well-defined
bounded \(\ICL\one_\GF\)\nbd{}equivariant linear operator as in the
proof of Lemma~\ref{lem:weighted_Sch_extension}.  The operator
\[
P_- \defeq
i_-\circ p_- =
\begin{pmatrix}
  M_{1-\phi} & 
  M_\phi \\
  M_{1-\phi} & 
  M_\phi
\end{pmatrix}.
\]
is a projection onto the range of~\(i_-\) because \(p_- i_- =
\ID_{\Hilm_-}\).

\subsection{The Poisson Summation Formula}
\label{sec:Poisson_summation}

In order to compare the operators \(P_+^S\) and~\(P_-\), we need the
Poisson Summation Formula, which is a special feature of the global
situation.  For any \(x\in\Adel\inv_\GF\), \(\GF\cdot x\subseteq
\Adel_\GF\) is a discrete cocompact subgroup.  Our normalization
of~\(\psi\) insures that the dual lattice \((x\cdot\GF)^\bot\) is equal to
\(x^{-1}\cdot\GF\).  Hence the \emph{Poisson Summation Formula} asserts
\[
\sum_{a\in\GF} f(ax) =
\abs{x}^{-1} \sum_{a\in\GF} \Fourier f(ax^{-1})
\]
for all \(f\in\Sch(\Adel_\GF)\), see~\cite{Weil:Basic}.  Define
\[
\Sing(f)\colon \ICL_S\to\C,
\qquad
x\mapsto f(0)-\abs{x}^{-1}\cdot \Fourier f(0),
\]
for \(f\in\Sch(\Adel_\GF)\).  This descends to a
\(\ICL_\GF\)\nbd{}equivariant bounded linear operator from~\(\Hilm_+\) to
the \(2\)\nbd{}dimensional space spanned by the constant function and
the function~\(\abs{x}^{-1}\).  The Poisson Summation Formula can be
rewritten as
\begin{equation}
  \label{eq:Poisson_summation}
  J \SUM \Fourier = \SUM + \Sing.
\end{equation}
When we restrict to~\(\Hilm_+^S\) for a sufficiently large finite set of
places~\(S\), we get \(J\Lop_S \Fourier = \Lop_S+\Sing\) on
\(\Twist(\ICL_S)\) because \(\Lop_S=\SUM\) on that subspace.

\begin{lemma}  \label{lem:Poisson_sumpm}
  The subspace \(i_+(\Hilm_+) + i_-(\Hilm_-)\) of \(\Sch(\ICL_\GF)_\ppmm\)
  is equal to
  \[
  \sumpm \defeq \{ (f_0,f_1)\in \Sch(\ICL_\GF)_\ppmm\mid
  \text{\(f_0-f_1 = c_0 -c_1\abs{x}^{-1}\) for some \(c_0,c_1\in\C\)} \}.
  \]
\end{lemma}

\begin{proof}
  It is clear that \(i_-(\Hilm_-)\subseteq \sumpm\) and that the
  quotient \(\sumpm/\Hilm_-\) is \(2\)\nbd{}dimensional.  The Poisson
  Summation Formula implies \(i_+(\Hilm_+)\subseteq \sumpm\) as well and
  yields that the remaining two basis vectors do arise because
  of~\(\Sing\).
\end{proof}

\begin{lemma}  \label{lem:global_perturbation}
  For \(\alpha>1\), both \(\IN\lambda(f) (P_- - P_+^S)\) and \((P_- -
  P_+^S) \IN\lambda(f)\) extend to uniformly nuclear operators
  \(L^2(\ICL_S)_\alpha\oplus L^2(\ICL_S)_{1-\alpha} \to
  \Sch(\ICL_\GF)_\ppmm\) for \(f\in\CCINF(\ICL_S)\).
\end{lemma}

\begin{proof}
  It suffices to check this on the subspaces of
  \(\MCR\inv_{\Cont{S}}\)\brd{}invariants and we may assume that~\(\phi\)
  is symmetric.  It is straightforward to verify the uniform
  nuclearity of \([\IN\lambda(f),P_-]= i_-[\IN\lambda(f),p_-]\) and
  \([\IN\lambda(f),P_+^S]= i_+[\IN\lambda(f),p_+^S]\) for all
  \(f\in\CCINF(\ICL_S)\).  Hence we only have to consider
  \begin{multline}  \label{eq:global_perturbation}
    \IN\lambda(f) (P_- - P_+^S)
    \\ =
    \begin{pmatrix}
      \IN\lambda(f) (\Lop_S M_{1-\phi} \Lopi_S - M_{1-\phi}) & 
      \IN\lambda(f) (\Lop_S \Fourier M_{1-\phi} \Lopi_S J - M_\phi) \\
      \IN\lambda(f) (J \Lop_S \Fourier^\ast M_{1-\phi} \Lopi_S -
      M_{1-\phi}) & 
      \IN\lambda(f) (J \Lop_S M_{1-\phi} \Lopi_S J - M_\phi)
    \end{pmatrix}.
  \end{multline}
  The map~\(J\) is a bornological isomorphism
  \[
  J\colon \Sch(\ICL_S)_\pp\congto\Sch(\ICL_S)_\mm,
  \qquad
  L^2(\ICL_S)_\alpha\congto L^2(\ICL_S)_{1-\alpha},
  \]
  we have \(\Fourier=\Fourier^\ast\) and \(JM_{1-\phi} J = M_\phi\),
  and \(J\IN\lambda(f)J=\IN\lambda(\check{f})\) is again of the same
  form.  Therefore, it suffices to check that
  \[
  \IN\lambda(f) (\Lop_S M_{1-\phi} \Lopi_S - M_{1-\phi})
  \quad\text{and}\quad
  \IN\lambda(f) (\Lop_S \Fourier M_{1-\phi} \Lopi_S - M_\phi J)
  \]
  extend to uniformly nuclear operators
  \(L^2(\ICL_S)_\alpha\to\Sch(\ICL_S)_\pp\).  Lemma~\ref{lem:summation}
  yields \(\Lop_S f\in\Sch(\ICL_S)_\pp\).  The same proof as
  for~\eqref{eq:Fourier_lambda} also yields \(\IN\lambda(f)\circ \Lop_S
  = \IN\lambda(\Lop_S f)\).  Hence
  \[
  \IN\lambda(f) (\Lop_S M_{1-\phi} \Lopi_S - M_{1-\phi})
  = [\IN\lambda(\Lop_S f),M_{1-\phi}] \Lopi_S
  - [\IN\lambda(f),M_{1-\phi}].
  \]
  Lemma~\ref{lem:mother_of_nuclearity} yields that this extends to
  a uniformly nuclear operator from \(L^2(\ICL_S)_\alpha\) to
  \(\Sch(\ICL_S)_\pp\) for \(f\in\CCINF(\ICL_S)\) for any \(\alpha>1\).  We
  have
  \begin{multline*}
    \IN\lambda(f) M_\phi J
    = [\IN\lambda(f),M_\phi] J
    + J M_{1-\phi} \IN\lambda(Jf) \Lop_S \Lopi_S
    \\ = [\IN\lambda(f),M_\phi] J
    + J [M_{1-\phi}, \IN\lambda(\Lop_S Jf)] \Lopi_S
    + J \Lop_S \IN\lambda(Jf) M_{1-\phi} \Lopi_S
  \end{multline*}
  and hence
  \begin{multline*}
    \IN\lambda(f) (\Lop_S \Fourier M_{1-\phi} \Lopi_S - M_\phi J)
    \\ =
    J (J\Lop_S\Fourier - \Lop_S) \IN\lambda(Jf) M_{1-\phi} \Lopi_S
    - [\IN\lambda(f),M_\phi] J
    - J [M_{1-\phi}, \IN\lambda(\Lop_S Jf)] \Lopi_S
    \\ =
    - [\IN\lambda(f),M_\phi] J
    - J([M_{1-\phi}, \IN\lambda(\Lop_S Jf)]
    + \IN\lambda(\Sing Jf) M_{1-\phi}) \Lopi_S,
  \end{multline*}
  where we used the Poisson Summation Formula.  The operator
  \([\lambda(f),M_\phi]J\) extends to a nuclear operator
  \(L^2(\ICL_S)_\alpha \to \Sch(\ICL_S)_\mm\) by
  Lemma~\ref{lem:mother_of_nuclearity}.  Since~\(\Lopi_S\) is bounded on
  \(L^2(\ICL_S)_\alpha\), it remains to prove that
  \[
  [M_{1-\phi}, \IN\lambda(\Lop_S Jf)] +
  \IN\lambda(\Sing Jf) M_{1-\phi}
  \]
  extends to a uniformly nuclear operator
  \(L^2(\ICL_S)_\alpha\to\Sch(\ICL_S)_\mm\).  The Poisson Summation
  Formula implies \(\Lop_S Jf= h-\phi(x)\cdot \Sing Jf\) with
  \(h\in\Sch(\ICL_S)_\cup\).  The operator \([M_{1-\phi},\IN\lambda(h)]\)
  extends to a nuclear operator from \(L^2(\ICL_S)_\alpha\) to
  \(\Sch(\ICL_S)_\cup\) by Lemma~\ref{lem:mother_of_nuclearity}.  We are
  left with
  \begin{multline*}
    [M_{1-\phi}, \IN\lambda(\phi\cdot \Sing Jf)] +
    \IN\lambda(\Sing Jf) M_{1-\phi}
    \\ =
    M_{1-\phi} \IN\lambda(\phi\cdot \Sing Jf) -
    \IN\lambda((1-\phi)\Sing Jf) M_{1-\phi}.
  \end{multline*}
  An examination of the integral kernels of \(M_{1-\phi}
  \IN\lambda(\phi\cdot \Sing Jf)\) and \(\IN\lambda((1-\phi)\Sing Jf)
  M_{1-\phi}\) shows that both define uniformly nuclear maps
  \(L^2(\ICL_S)_\alpha\to\Sch(\ICL_S)_\mm\) for \(f\in\CCINF(\ICL_S)\) as
  desired.
\end{proof}

\subsection{Construction of the global difference representation}
\label{sec:construct_global_difference}

Let
\begin{align*}
  \Hilm^0_+ &\defeq
  (\Hilm_++\Hilm_-)/\Hilm_- \cong \Hilm_+/(\Hilm_+\cap\Hilm_-),
  \\
  \Hilm^0_- &\defeq
  (\Hilm_++\Hilm_-)/\Hilm_+ \cong \Hilm_-/(\Hilm_+\cap\Hilm_-),
\end{align*}
where we take sums and intersections in \(\Sch(\ICL_\GF)_\ppmm\).
Let~\(\Hilm^0\) be the virtual representation
\(\Hilm^0_+\ominus\Hilm^0_-\).  This is the \emph{global difference
  representation} of~\(\ICL_\GF\).  We denote the canonical
representations of~\(\ICL_\GF\) on \(\Hilm^0_\pm\) and~\(\Hilm^0\) by
\(\pi_\pm\) and~\(\pi\), respectively.  We are going to show that~\(\pi\) is
summable.  The \emph{character} of~\(\pi\) is defined by \(\chi_\pi
\defeq \chi_{\pi_+} - \chi_{\pi_-}\).

\begin{definition}  \label{def:Weil_distribution}
  For any place~\(v\), define
  \[
  W_v\in\CCINF'(\ICL_\GF), \qquad W_v(f) \defeq \int'_{\GF\inv_v}
  \frac{f(x)\abs{x}}{\abs{1-x}} \,d\inv x,
  \]
  where~\(\int'_{\GF\inv_v}\) denotes the same principal value as in
  the local trace formula.  Let~\(\Discriminant_\GF\) be the discriminant
  of~\(\GF\) and let
  \[
  W(f) \defeq \sum_{v\in\Places(\GF)} W_v(f) - f(1)\cdot
  \ln{}\Discriminant_\GF
  \qquad
  \text{for \(f\in\CCINF(\ICL_\GF)\).}
  \]
  We call \(W\in\CCINF'(\ICL_\GF)\) the \emph{(raw) Weil distribution}.
\end{definition}

\begin{theorem}  \label{the:global_trace_formula}
  The virtual representation~\(\pi\) is smooth and summable and its
  character is equal to the Weil distribution~\(W\).  Thus
  \[
  W(f) =
  \chi_\pi(f) =
  \sum_{\omega\in\Irrep(\ICL_\GF)} \mult(\omega,\pi) \cdot
  \hat{f}(\omega)
  \]
  with uniform absolute convergence for~\(f\) in bounded subsets of
  \(\CCINF(G)\).
\end{theorem}

\begin{proof}
  Since smoothness is hereditary for subrepresentations and quotient
  representations, it is clear that~\(\pi_\pm\) and hence~\(\pi\) are
  smooth representations.  Recall the description of \(\sumpm\defeq
  \Hilm_+ + \Hilm_-\) in Lemma~\ref{lem:Poisson_sumpm}.  Since
  \(\sumpm/\Hilm_-=\Hilm_+^0\) is \(2\)\nbd{}dimensional, it is evidently
  summable.  Since all spaces that occur below are nuclear, a nuclear
  operator is automatically of order~\(0\).  Hence the restriction of a
  nuclear operator to an invariant subspace and the induced operator
  on a quotient space are again nuclear by
  Theorem~\ref{the:order_zero_restrict}.  Since the range of \(P_+^S\)
  is contained in \(\Hilm_+\subseteq\sumpm\) and the range of~\(P_-\) is
  contained in \(\Hilm_-\subseteq\sumpm\), all operators that we
  consider in the following map \(\Sch(\ICL_S)_\ppmm\) into the
  subspace~\(\sumpm\).  Thus they are nuclear as operators on~\(\sumpm\)
  if and only if they are nuclear as operators on \(\Sch(\ICL_S)_\ppmm\)
  and they have the same trace on both spaces.
  
  Let \(f\in\CCINF(\ICL_\GF)\).  There is \(\beta\ge1\) such that
  \(\abs{\supp f}\in \ooival{\beta^{-1},\beta}\) and there is a
  sufficiently large finite set of places~\(S\) such that~\(f\) is
  unramified outside~\(S\) and \(q_v\ge\beta\) for all \(v\in\Cont{S}\).
  Even more, we can choose \(\beta\) and~\(S\) uniformly for~\(f\) in a
  bounded subset of \(\CCINF(\ICL_\GF)\).  Consider the following
  operator on~\(\sumpm\):
  \[
  \IN\lambda(f) (\ID_\sumpm - P_+^S) =
  \IN\lambda(f) (\ID_\sumpm - P_-) +
  \IN\lambda(f) (P_- - P_+^S).
  \]
  The operator \(\IN\lambda(f) (\ID_\sumpm - P_-)\) on~\(\sumpm\)
  vanishes on the finite codimensional subspace~\(\Hilm_-\) and induces
  \(\IN\pi_+(f)\) on~\(\Hilm^0_+\).  Hence it has finite rank and its
  trace is \(\chi_{\pi_+}(f)\).  The operator \(\IN\lambda(f) (P_- -
  P_+^S)\) is of order~\(0\) by Lemma~\ref{lem:global_perturbation}.
  Hence the operator \(\IN\lambda(f)(\ID_\sumpm-P_+^S)\) is of order~\(0\)
  as well.  It leaves~\(\Hilm_+\) invariant and induces \(\IN\pi_-(f)\)
  on~\(\Hilm^0_-\) because the range of~\(P_+^S\) is contained
  in~\(\Hilm_+\).  Theorem~\ref{the:order_zero_restrict} now implies
  that~\(\pi_-\) is summable.

  The additivity of the trace in extensions yields
  \begin{multline*}
    \chi_{\pi_-}(f) =
    \tr \IN\lambda(f)(\ID_\sumpm-P_+^S)
    - \tr \IN\lambda(f)(\ID - P_+^S)|_{\Hilm_+}
    \\ =
    \tr \IN\lambda(f) (\ID_\sumpm - P_-)
    + \tr \IN\lambda(f) (P_- - P_+^S)
    - \tr (\ID - P_+^S)|_{\Hilm_+}.
  \end{multline*}
  Since \(\tr \IN\lambda(f) (\ID_\sumpm - P_-)=\chi_{\pi_+}(f)\), this
  implies
  \[
  \chi_\pi(f) =
  \tr \IN\lambda(f)(\ID - P_+^S)|_{\Hilm_+} -
  \tr \IN\lambda(f) (P_- - P_+^S).
  \]
  Since we embed~\(\Hilm_+\) into~\(\sumpm\) using~\(i_+\), the first
  term is
  \begin{multline*}
    \tr \IN\lambda(f)(\ID - P_+^S)|_{\Hilm_+}
    =
    \tr \IN\lambda(f)(\ID_{\Hilm_+} - p_+^S i_+)|_{\Twist(\ICL_S)}
    =
    \tr \IN\lambda(f)(M_\phi - \EFou M_\phi\EFoi)|_{\Twist(\ICL_S)}
    \\ =
    \tr \IN\lambda(f)(M_\phi - \EFou M_\phi\EFoi)|_{L^2(\ICL_S)_{1/2}}
    =
    \sum_{v\in S} W_v(f) - f(1) \ln{}\abs{\different}^{-1}_S
    =
    W(f).
  \end{multline*}
  We have assumed~\(\phi\) symmetric to simplify the formula and we
  have used the local trace formula~\eqref{eq:local_trace_formula_psi}
  and the fact that it makes no difference whether we compute this
  trace on \(L^2(\ICL_S)_{1/2}\) or on \(\Twist(\ICL_S)\).  Our choice
  of~\(S\) guarantees \(W_v(f)=0\) for \(v\in\Cont{S}\) and
  \(\abs{\different}_S^{-1}=\Discriminant_\GF\).
  
  It remains to show \(\tr \IN\lambda(f) (P_- - P_+^S)=0\).
  By~\eqref{eq:global_perturbation}, this amounts to
  \[
  \tr \IN\lambda(f) (\Lop_S M_{1-\phi} \Lopi_S - M_{1-\phi}) =
  \tr \IN\lambda(f) J(\Lop_S M_{1-\phi} \Lopi_S - M_{1-\phi})J = 0.
  \]
  Since
  \begin{multline*}
    \tr \IN\lambda(f) J(\Lop_S M_{1-\phi} \Lopi_S - M_{1-\phi})J =
    \tr J \IN\lambda(f) J(\Lop_S M_{1-\phi} \Lopi_S - M_{1-\phi})
    \\ =
    \tr \IN\lambda(Jf) (\Lop_S M_{1-\phi} \Lopi_S - M_{1-\phi}),
  \end{multline*}
  it suffices to treat \(\IN\lambda(f) (\Lop_S M_{1-\phi} \Lopi_S -
  M_{1-\phi})\).  The proof of Lemma~\ref{lem:global_perturbation}
  yields that this operator has an integral kernel in
  \(\Sch(\ICL_S)_\pp\hot\Sch(\ICL_S)_\mm\), so that its trace is equal
  to the integral of its integral kernel over the diagonal.  Our
  choice of~\(S\) and~\(\beta\) insures that \(q_v\ge\beta\) for all
  \(v\in\Cont{S}\), so that \(\Lopi_S-1\) and \(\Lop_S-1\) are sums of
  operators~\(\lambda_x\) with \(\abs{x}\le \beta^{-1}\).  The support
  condition on~\(f\) implies that the integral kernel of \(\IN\lambda(f)
  (\Lop_S M_{1-\phi} \Lopi_S - M_{1-\phi})\) vanishes identically near
  the diagonal.  Hence the trace also vanishes.
\end{proof}

Let~\(\chi\) be a quasi-character of~\(\ICL_\GF\).  Then there are unique
\(s\in\R\) and \(\chi_0\in\widehat{\ICL_\GF}\) such that
\(\chi=\abs{x}^s\cdot\chi_0\).  We call \(s\defeq \RE \chi\) the
\emph{real part} of~\(\chi\).

\begin{proposition}  \label{pro:difference_in_strip}
  The integrated form of~\(\pi\) extends to a bounded homomorphism
  \[
  \IN\pi \colon \varinjlim_{\epsilon\searrow0}
  \Sch(\ICL_\GF)_{[-\epsilon,1+\epsilon]} \to \ell^1(\Hilm^0).
  \]
  The spectrum of~\(\pi\) is contained in the region
  \(\{\omega\in\Irrep(\ICL_\GF)\mid \RE\omega\in[0,1]\}\), which we call
  the \emph{critical strip}.
\end{proposition}

\begin{proof}
  The assertion is clear for the positive part
  \(\Hilm^0_+\cong\C(1)\oplus\C(\abs{x})\), so that we restrict
  attention to~\(\Hilm^0_-\) in the following.  Fix a sufficiently
  large finite set of places~\(S\) and \(\epsilon>0\).
  Lemma~\ref{lem:global_perturbation} asserts that \((P_- -
  P_+^S)\IN\lambda(f)\) extends to a uniformly nuclear operator
  \(L^2(\ICL_S)_{1+\epsilon}\oplus L^2(\ICL_S)_{-\epsilon} \to
  \Sch(\ICL_S)_\ppmm\) for \(f\in\CCINF(\ICL_S)\).  The regular
  representations on \(L^2(\ICL_S)_{1+\epsilon}\) and
  \(L^2(\ICL_S)_{-\epsilon}\) integrate to a module structure over
  \(\Sch(\ICL_S)_{[-\epsilon,1+\epsilon]}\).  Hence the operators \((P_-
  - P_+^S)\IN\lambda(f_1)\IN\lambda(f_2)\) are uniformly nuclear
  operators on \(\Sch(\ICL_S)_\ppmm\) for \(f_1\in\CCINF(\ICL_S)\) and
  \(f_2\in\Sch(\ICL_S)_{[-\epsilon,1+\epsilon]}\).  Since
  \begin{displaymath}
    \Sch(\ICL_S)_{[-\epsilon,1+\epsilon]}
    \cong
    \CCINF(\ICL_S) \hot_{\CCINF(\ICL_S)}
    \Sch(\ICL_S)_{[-\epsilon,1+\epsilon]},
  \end{displaymath}
  we obtain that the map that sends \(f\in\CCINF(\ICL_S)\) to \((P_- -
  P_+^S)\IN\lambda(f)\) extends to a bounded linear map from
  \(\Sch(\ICL_S)_{[-\epsilon,1+\epsilon]}\) to
  \(\ell^1(\Sch(\ICL_S)_\ppmm)\).  This implies the corresponding
  statement for~\(\ICL_\GF\) instead of~\(\ICL_S\) because \(\IN\lambda(f)\)
  for smooth~\(f\) projects to \(\MCR\inv_{\Cont{S}}\)\brd{}invariants,
  anyway.  We realize~\(\Hilm^0_-\) as a subquotient of
  \(\Sch(\ICL_\GF)_\ppmm\) and \(\IN\pi_-(f)\) as the operator on this
  subquotient induced by \((P_- - P_+^S)\IN\lambda(f)\) as in the proof
  of Theorem~\ref{the:global_trace_formula}.
  Theorem~\ref{the:order_zero_restrict} yields that~\(\IN\pi_-\) extends
  to a bounded homomorphism
  \(\Sch(\ICL_\GF)_{[-\epsilon,1+\epsilon]}\to \ell^1(\Hilm^0_-)\).
  Since this holds for all \(\epsilon>0\), we get the assertion.

  It follows that any closed subrepresentation of~\(\pi\) is an essential
  module over \(\Sch(\ICL_\GF)_{[-\epsilon,1+\epsilon]}\) for all
  \(\epsilon>0\).  Since \(\C(\chi)\) has this property if and only if \(\RE
  \chi\in[0,1]\), we obtain the assertion about the spectrum of~\(\pi\).
\end{proof}

A more careful analysis shows that~\(\IN\pi\) also extends to a
representation of \(\Sch(\ICL_\GF)_{[0,1]}\).  We cannot shrink \([0,1]\)
any further because of the two quasi-characters \(\abs{x}^0\)
and~\(\abs{x}^1\) in \(\spec\pi_+\).

\subsection{Symmetries of the global difference representation}
\label{sec:symmetries}

It is rather hard to construct interesting operators on~\(\Hilm^0\).
There are, however, some obvious symmetries that correspond to
well-known symmetries of \(L\)\nbd{}functions.  The operators~\(\Fourier\)
on~\(\Hilm_+\) and~\(J\) on~\(\Hilm_-\) agree on \(\Hilm_+\cap\Hilm_-\) by the
Poisson Summation Formula.  Hence they yield an operator, again
called~\(J\), on~\(\Hilm^0\).  This operator is idempotent and satisfies
\(J\IN\lambda(f)J=\IN\lambda(Jf)\).  Complex conjugation on
\(\Sch(\ICL_\GF)_\ppmm\) leaves the subspaces~\(\Hilm_\pm\) invariant and
hence descends to~\(\Hilm^\pm_0\).  This operator~\(\gamma\) on~\(\Hilm_0\)
is conjugate linear and satisfies \(\gamma\IN\pi(f)\gamma =
\IN\pi(\conj{f})\) for all \(f\in\CCINF(\ICL_\GF)\), where
\(\conj{f}(x)=\conj{f(x)}\).  The covariance properties of \(J\)
and~\(\gamma\) imply that they correspond to the symmetries
\(\chi\mapsto\abs{x}\chi^{-1}\) and \(\chi\mapsto\conj{\chi}\) on the
spectral side.  Hence the spectrum of~\(\pi\) is invariant under these
two maps and
\[
\mult(\chi,\pi) =
\mult(\abs{x}\chi^{-1},\pi) =
\mult(\conj{\chi},\pi)
\qquad
\text{for all \(\chi\in\Irrep(\ICL_\GF)\).}
\]

There is another more geometric spectral interpretation for the poles
and zeros of \(L\)\nbd{}functions in the case of global function fields,
which involves the action of the Frobenius automorphism on an
appropriate étale cohomology.  The Galois group of the global function
field also acts on étale cohomology, and this action is of great
interest because it can be used to describe poles and zeros of Artin
\(L\)\nbd{}functions as well.  In our spectral interpretation we also
have an action of the Galois group, but it is much less interesting.
In particular, it has nothing to do with the Galois action on étale
cohomology for global function fields.  Let \(\tau\colon\GF\to\GF\) be
an automorphism of the field~\(\GF\).  All our constructions are natural
in an evident fashion, so that~\(\tau\) acts on~\(\Hilm^0\).  However,
this action is not very interesting.  The action on~\(\Hilm^0_+\) is
trivial and the action on~\(\Hilm^0_-\) is the quotient representation
of the natural representation on \(\Sch(\ICL_\GF)_\cup\) that comes from
the action of~\(\tau\) on~\(\ICL_\GF\).  Hence we obtain
\(\mult(\chi,\pi)=\mult(\chi\circ\tau_\ast,\pi)\) for all
\(\chi\in\Irrep(\ICL_\GF)\).  Since~\(\tau_\ast\) leaves the norm
on~\(\ICL_\GF\) invariant, its action on \(\Irrep(\ICL_\GF)\) merely
permutes the connected components.  Thus it has nothing to do with
Artin \(L\)\nbd{}functions.

By definition, \(\Hilm_-^0\) is a quotient of \(\Sch(\ICL_\GF)_\cup\) by
\(\Hilm_+\cap\Hilm_-\), which is an ideal with respect to convolution.
Hence~\(\IN\pi_-\) drops down to a bounded bilinear map
\[
\Hilm_-^0\times\Hilm_-^0 \to \Hilm_-^0,
\qquad
[f_1]\otimes [f_2]\mapsto [\IN\pi_-(f_1)(f_2)],
\]
which turns~\(\Hilm_-^0\) into an associative algebra.  Any character on
this algebra gives rise to a character on \(\Sch(\ICL_\GF)_\cup\) and
thus to a quasi-character of~\(\ICL_\GF\).  Thus the spectrum
of~\(\Hilm_-^0\) identifies with a subset of \(\Irrep(\ICL_\GF)\).  By
definition, this spectrum is exactly the spectrum of the
contragradient representation~\(\tilde\pi_-\).  It is easy to see
that~\(\pi_-\) and~\(\tilde\pi_-\) have the same spectrum and the same
spectral multiplicity function, see also
Theorem~\ref{the:summable_hereditary}.  The involution~\(J\gamma\)
turns~\(\Hilm_-^0\) into a \(\ast\)\nbd{}algebra and~\(\gamma\) gives an
additional real structure on~\(\Hilm_-^0\).  A character on~\(\Hilm_-^0\)
is a \(\ast\)\nbd{}character if and only if it corresponds to
\(\omega\in\Irrep(\ICL_\GF)\) with \(\RE\omega=1/2\).  Thus the
Generalized Riemann Hypothesis is equivalent to the assertion that
self-adjoint elements of the \(\ast\)\nbd{}algebra~\(\Hilm_-^0\) have real
spectrum.

\subsection{Spectral analysis of the global difference representation}
\label{sec:global_spectral_analysis}

\begin{lemma}  \label{lem:Lap_weighted_Sch}
  Let~\(I\) be an open interval.  The Fourier-Laplace transform maps
  \(\Sch(\ICL_\GF)_I\) isomorphically onto the space of holomorphic
  functions on the domain \(\{\omega\in\Irrep(\ICL_\GF)\mid \RE
  \omega\in I\}\) whose restriction to
  \(\abs{x}^\alpha\widehat{\ICL_\GF}\) is a Schwartz function for all
  \(\alpha\in I\).
\end{lemma}

\begin{proof}
  We use that the Fourier transform is an isomorphism
  \(\Sch(\ICL_\GF)\cong\Sch(\widehat{\ICL_\GF})\) and that
  \((f\cdot\abs{x}^s)\sphat(\omega)=\hat{f}(\omega\cdot\abs{x}^s)\) for
  all \(s\in\C\).  In particular, we may assume without loss of
  generality that \(0\in I\).  Multiplication by \(\abs{x}^{\ima t}\) for
  \(t\in\R\) leaves \(\Sch(\ICL_\GF)\) and \(\Sch(\ICL_\GF)_I\) invariant.
  There is only one possible way to extend the function \(\ima t\mapsto
  f\cdot\abs{x}^{\ima t}\) to a holomorphic function, namely, the
  function \(s\mapsto f\cdot\abs{x}^s\).  Hence a function in
  \(\Sch(\ICL_\GF)\) belongs to \(\Sch(\ICL_\GF)_I\) if and only if the
  function \(\ima t\mapsto f\cdot\abs{x}^{\ima t}\) extends to a
  holomorphic function from the strip \(\{s\in\C\mid \RE s\in I\}\) to
  \(\Sch(\ICL_\GF)\).  Applying the Fourier-Laplace transform, we get
  that \(\Sch(\ICL_\GF)_I\sphat\) consists of those
  \(f\in\Sch(\widehat{\ICL_\GF})\) for which the function \(\ima t\mapsto
  f(\omega\cdot\abs{x}^{\ima t})\) extends to a holomorphic function
  from the strip \(\{\RE s\in I\}\) to \(\Sch(\widehat{\ICL_\GF})\).
  Clearly, such a holomorphic extension is equivalent to a holomorphic
  extension of~\(f\) to a function on the strip \(\{\RE\omega\in I\}\)
  with the property that its restriction
  \(\abs{x}^\alpha\widehat{\ICL_\GF}\) is a Schwartz function for all
  \(\alpha\in I\).
\end{proof}

Suppose now that \(M\colon \Sch(\ICL_\GF)_I\to\Sch(\ICL_\GF)_I\) is a
bounded \(\ICL_\GF\)\nbd{}equivariant linear operator.  In the
Fourier-Laplace transformed picture, \(M\) must act as pointwise
multiplication by some holomorphic function
\[
\hat{M}\colon \{\omega\in\Irrep(\ICL_\GF)\mid \RE\omega\in I\} \to\C.
\]
If \(\alpha\in I\), then we have \(\Sch(\ICL_\GF)_I
\hot_{\Sch(\ICL_\GF)_I} \Sch(\ICL_\GF)_\alpha \cong
\Sch(\ICL_\GF)_\alpha\).  Hence~\(M\) gives rise to a bounded operator on
\(\Sch(\ICL_\GF)_\alpha\) for all \(\alpha\in I\).  Thus a holomorphic
function on the strip \(\{\RE\omega\in I\}\) multiplies
\(\Sch(\ICL_\GF)_I\sphat\) if and only if its restriction to
\(\widehat{\ICL_\GF}\cdot\abs{x}^\alpha\) multiplies
\(\Sch(\widehat{\ICL_\GF})\) for each \(\alpha\in I\).

The Poisson Summation Formula implies \(\SUM f = h + \phi\cdot\Sing f\)
with \(h\in\Sch(\ICL_\GF)_\cup\).  Thus the Fourier-Laplace transform
of~\(h\) is a holomorphic function on all of \(\Irrep(\ICL_\GF)\) and its
restriction to each vertical line is a Schwartz function.  It is easy
to see that a function of the form \(\widehat{\phi\chi}\) for
\(\chi\in\Irrep(\ICL_\GF)\) is a meromorphic function whose only pole is
a simple pole at~\(\chi^{-1}\) (we show this during the proof of
Theorem~\ref{the:boundary_spectrum} below).  Moreover,
\(\widehat{\phi\chi}\) is a Schwartz function on
\(\abs{x}^\alpha\widehat{\ICL_\GF}\) for all
\(\alpha\in\R\setminus\{-\RE\chi\}\).  We get a Schwartz function on
\(\chi^{-1}\widehat{\ICL_\GF}\) if we multiply with any function that
vanishes at~\(\chi^{-1}\) and is constant at~\(\infty\).  Hence \(\SUM(f)\)
is a meromorphic function whose only poles are simple poles at the
quasi-characters \(\abs{x}\) and~\(1\) and whose restriction to vertical
lines \(\abs{x}^\alpha\widehat{\ICL_\GF}\) is a Schwartz function for
all \(\alpha\in\R\) except for its two poles.  This implies the
meromorphic continuation of the \(L\)\nbd{}functions \(L_\GF(\chi,s)\)
because \(L_\GF(\chi,s)\) is \(\widehat{\SUM f_\chi}(\chi\abs{x}^s)\) for
suitable \(f_\chi\in\Sch(\Adel_\GF)\) for all
\(\chi\in\widehat{\ICL_\GF}\).

The Poisson Summation Formula is equivalent to
\begin{equation}  \label{eq:Poisson_spectral}
  \widehat{\SUM f}(\chi) = \widehat{\SUM \Fourier f}(\abs{x}\chi^{-1})
\end{equation}
for all \(\chi\in\Irrep(\ICL_\GF)\setminus \{1,\abs{x}\}\) and all
\(f\in\Sch(\Adel_\GF)/\GF\inv\).  This is evident if \(\Sing f=0\) and
extends to general~\(f\) because of the equivariance of~\(\SUM\).  When we
apply~\eqref{eq:Poisson_spectral} to the special functions~\(f_\chi\)
that give the \(L\)\nbd{}functions, we obtain the functional equations
for the \(L\)\nbd{}functions.  This proof method is due to John T.\ Tate
(\cite{Tate:Thesis}).

For global function fields, we can replace \(\Sch(\ICL_\GF)_\cup\) by
\(\Sch_\comp(\ICL_\GF)_\cup=\CCINF(\ICL_\GF)\).  The Fourier-Laplace
transform maps \(\CCINF(\ICL_\GF)\) onto the space of functions on
\(\Irrep(\ICL_\GF)\) whose restriction to each connected component is a
Laurent polynomial and which are supported in finitely many connected
components.  Thus we obtain that \(L(\chi,s)\) for \(\chi\neq1\) and
\(\zeta(s)\cdot(q_\GF^s-1)(q_\GF^{s-1}-1)\) are Laurent polynomials
in~\(q_\GF^s\).

Let~\(S\) be a sufficiently large set of places.  We now consider the
functions on \(\Irrep(\ICL_S)\) that correspond to
\(\Lop_S\) and~\(\EFou\).  The operator~\(\Lop_S\) corresponds to the
operator of pointwise multiplication by the function
\[
\widehat{\Lop_S}(\omega) \defeq
\prod_{v\in\Cont{S}} (1-\omega(p_v)^{-1})^{-1},
\]
where \(p_v\in\ICL_S\) for \(v\in\Cont{S}\) is as in the definition
of~\(\Lop_S\).  This infinite product converges in the region \(\RE
\omega>1\).  The same infinite product occurs in the definition of the
\(L\)\nbd{}function of~\(\GF\).  Namely, \(L_\GF(\omega)=
\widehat{\Lop_S}(\omega)\cdot L_{\GF,S}(\omega)\) for all
\(\omega\in\Irrep(\ICL_S)\) with \(\RE\omega>1\), where \(L_{\GF,S}\) is the
product of the local \(L\)\nbd{}factors for \(v\in S\).  Our analysis
of~\(\Lop_S\) shows that it is an invertible multiplier of
\(\Sch(\abs{x}^\alpha\widehat{\ICL_S})\) for \(\alpha>1\).  Similarly,
\(\EFou\) is a bounded operator on \(\Sch(\ICL_S)_+\), so that its
Fourier-Laplace transform is a holomorphic function in the region
\(\RE\omega>0\).  Its inverse~\(\EFoi\) is a bounded operator on
\(\Sch(\ICL_S)_-\), so that its Fourier-Laplace transform is holomorphic
in the region \(\RE \omega<1\).  Thus~\(\widehat{\EFou}\) extends to a
meromorphic function on all of \(\Irrep(\ICL_S)\), which has no poles
for \(\RE\omega>0\) and no zeros for \(\RE\omega<1\).

Next we describe the Fourier-Laplace transform of \(\SUM(\Hilm_+^S)\).  We
identify \(\Hilm_+^S\) with \(\Twist(\ICL_S)\) and the summation map~\(\SUM\)
with~\(\Lop_S\).  By definition, the Fourier-Laplace transform of
\(\Twist(\ICL_S)\) consists of those meromorphic functions~\(f\) on
\(\Irrep(\ICL_S)\) whose restriction to \(\abs{x}^\alpha\widehat{\ICL_S}\)
is Schwartz (and thus has no poles) for all \(\alpha>0\) and for which
\(\widehat{\EFoi}\cdot f\) is Schwartz on \(\abs{x}^\alpha\widehat{\ICL_S}\)
for \(\alpha<1\).  Both conditions are equivalent for \(0<\alpha<1\).  Thus
\(\SUM(\Hilm_+^S)\) corresponds to the space of meromorphic functions~\(f\)
on \(\Irrep(\ICL_S)\) for which \(f/\widehat{\Lop_S}\) is a Schwartz
function on \(\abs{x}^\alpha\widehat{\ICL_S}\) for \(\alpha>0\) and for
which \(f/\widehat{\EFou}\widehat{\Lop_S}\) is a Schwartz function on
\(\abs{x}^\alpha\widehat{\ICL_S}\) for \(\alpha<1\).  If \(0<\RE\alpha<1\),
then both conditions for \(\RE\omega=\alpha\) are equivalent.  If
\(\RE\alpha>1\), we can omit the factor \(\widehat{\Lop_S}\) because it is
invertible in this region.  If \(\RE\alpha<0\), we can also omit
\(\EFou\Lop_S\) because of the Poisson Summation Formula.  That is,
\(f\in(\SUM\Hilm_+^S)\sphat\) if and only if~\(f\) is a Schwartz function on
\(\abs{x}^\alpha\cdot\widehat{\ICL_S}\) for \(\alpha>1\) and \(\alpha<0\),
\(f/\widehat{\Lop_S}\) is a Schwartz function on
\(\abs{x}^\alpha\cdot\widehat{\ICL_S}\) for \(\alpha\in\ocival{0,1}\) and
\(f(\chi)/\widehat{\Lop_S}(\abs{x}\chi^{-1})\) is a Schwartz function for
\(\chi\in\widehat{\ICL_\GF}\).  This is exactly the information about
\(L\)\nbd{}functions that one obtains using only that they are of the form
\(\widehat{\SUM f_\chi}\) for some \(f_\chi\in\Sch(\Adel_\GF)\).

Finally, we relate \(\spec \pi_-\) to the zeros of \(L\)\nbd{}functions.
Let \(\ord(\omega,L_\GF)\) be the pole order of~\(L_\GF\) at~\(\omega\);
this is positive for poles and negative for zeros.

\begin{theorem}  \label{the:difference_spectrum}
  The spectrum of~\(\pi\) consists of the poles and zeros of the
  \(L\)\nbd{}function \(L_\GF\colon \Irrep(\ICL_\GF)\to\C\) and we have
  \(\mult(\omega,\pi) = \ord(\omega, L_\GF)\) for all
  \(\omega\in\Irrep(\ICL_\GF)\).
\end{theorem}

\begin{proof}
  The functions \(\mult(\omega,\pi)\) and \(\ord(\omega,L_\GF)\) are both
  invariant under \(\omega\mapsto \abs{x}\omega^{-1}\) by the Poisson
  Summation Formula.  Hence it suffices to consider
  \(\omega\in\Irrep(\ICL_\GF)\) with \(\RE\omega\ge1/2\).  It is clear
  that the two representations \(1\) and~\(\abs{x}\) that are contained in
  \(\pi_+\) are exactly the poles of~\(L_\GF\).  They also have correct
  multiplicity~\(1\).  Instead of describing the spectrum of~\(\pi_-\), it
  is easier to describe the spectrum of its
  contragradient~\(\tilde{\pi}_-\).  Both representations have the same
  character and hence the same spectral multiplicity function by
  Theorem~\ref{the:summable_hereditary}.  Since~\(\Hilm_-^0\) is a
  quotient of~\(\Hilm_-\), its contragradient is a subrepresentation of
  \((\Hilm_-)\sptilde\).  Clearly, a functional \(l\colon \Hilm_-\to\C\)
  is \(\omega\)\nbd{}homogeneous if and only if \(l(f)=0\) for all~\(f\)
  that satisfy \(\hat{f}(\omega)=0\).  Hence the only possible
  \(\omega\)\nbd{}homogeneous functionals are multiples of the
  evaluation map at~\(\omega\).  In order to get functionals on the
  quotient~\(\Hilm^0_-\), we need \(\hat{f}(\omega)=0\) for all
  \(f\in\Hilm_+\cap\Hilm_-\).  By our description of~\(\Hilm_+\), this is
  equivalent to \(\widehat{\Lop_S}(\omega)=0\) if \(\RE\omega>0\).  In
  this region, \(\widehat{\Lop_S}\) and~\(L_\GF\) have the same zeros.
  Thus~\(\omega\) is contained in~\(\pi_-\) if and only if it is a zero
  of~\(L_\GF\).  It is easy to see that the multiplicities agree as
  well.
\end{proof}

\section{The Prime Number Theorem}
\label{sec:PNT}

We are going to describe the asymptotic distribution of prime ideals
using our representation theoretic approach.  Actually, our argument
uses nothing more than the Explicit Formula in the form \(W(f) =
\sum_{\omega\in\Irrep(\ICL_\GF)} n(\omega) \hat{f}(\omega)\) with very
few properties of the function \(n\colon \Irrep(\ICL_\GF)\to\Z\).

\subsection{Prime ideals and the Weil distribution}
\label{sec:Primes_Weil}

Let~\(\GF\) be a global field and let~\(S\) be a sufficiently large finite
set of places of~\(\GF\).  We obtain a subring \(\GF_S\subseteq\GF\) and
its multiplicative group of units~\(\GF\inv_S\).  Let
\(\GF^\ast_S\defeq\GF_S\setminus\{0\}\).  We map
\(\GF^\ast_S\to\Adel\inv_S\to\ICL_S\), where the first map is the
diagonal map and the second one is the natural projection.  The image
of~\(\GF^\ast_S\) is a discrete subset of~\(\ICL_S\).  Two elements
in~\(\GF^\ast_S\) have the same image in~\(\ICL_S\) if and only if they
differ by a unit in~\(\GF\inv_S\).  Equivalently, they generate the same
ideal in~\(\GF_S\).  The prime ideals in~\(\GF_S\) correspond bijectively
to the places in~\(\Cont{S}\).  The element \(p_v\in\ICL_S\) for
\(v\in\Cont{S}\) that we used in the description of the summation map is
a generator for the prime ideal corresponding to the place~\(v\).  Thus
all prime ideals of~\(\GF_S\) are principal and the subset
\(\{p_v\}_{v\in\Cont{S}}\subseteq\ICL_S\) can be identified with the set
of prime ideals of~\(\GF_S\).  If \(\GF=\Q\) and~\(S\) consists just of the
infinite place, then the elements~\(p_v\) are exactly the prime numbers
in~\(\N\).  All this serves to justify the following

\begin{definition}  \label{def:prime_distribution}
  The \emph{prime ideal distribution} of~\(\GF_S\) is the distribution
  \[
  \Pi_S\colon \CCINF(\ICL_S)\to\C,
  \qquad
  \Pi_S(f)\defeq \sum_{v\in\Cont{S}} f(p_v).
  \]
\end{definition}

We are going to compare the prime ideal distribution~\(\Pi_S\) to the
Weil distribution~\(W\).  We restrict~\(W\) to
\(\CCINF(\ICL_S)\subseteq\CCINF(\ICL_\GF)\) and still denote it by~\(W\).
The normalization of the principal value in the definition of~\(W_v\)
yields
\begin{multline*}
  W_v(f) =
  \int_{\{x\in\GF\inv_v\mid \abs{x}\neq1\}}
  \frac{f(x)\abs{x}}{\abs{1-x}} \,d\inv x =
  \sum_{e=1}^\infty f(p_v^{-e}) q_v^{-e} \ln(q_v) +
  \sum_{e=1}^\infty f(p_v^e) \ln(q_v)
\end{multline*}
for all \(v\in\Cont{S}\).  This suggests that we define ``error terms''
\begin{align*}
  E_v(f) &\defeq
  \sum_{e=2}^\infty f(p_v^{-e}) q_v^{-e} \ln(q_v) +
  \sum_{e=2}^\infty f(p_v^e) \ln(q_v),
  \\
  E_S(f) &\defeq
  \sum_{v\in\Cont{S}} E_v + \sum_{v\in S} W_v -
  \ln(\Discriminant_\GF) \delta_1
\end{align*}
in such a way that
\[
W(f) = \Pi_S(f\cdot\ln{}\abs{x}) +
\Pi_S\bigl((Jf)\cdot\ln{}\abs{x}\bigr) +  E_S(f)
\]
for all \(f\in\CCINF(\ICL_S)\).

\begin{lemma}  \label{lem:estimate_Wtwo}
  The distribution~\(E_S\) extends to a bounded linear functional on
  \(\Sch(\ICL_S)_{1/2}\).  For any \(\epsilon>0\), the distribution
  \((\ln{}\abs{x})^{-1-\epsilon} \abs{x}^{-1/2}\cdot E_S|_{\abs{x}\ge
  1+\epsilon}\) is of the form \(f\mapsto \int f d\mu\) for some positive
  measure~\(\mu\) of finite total volume.
\end{lemma}

\begin{proof}
  We have already seen in Section~\ref{sec:local_trace_formula}
  that~\(W_v\) is bounded on \(\Sch(\GF\inv_v)_{1/2}\) for any place~\(v\).
  Since the distributions~\(E_v\) are \(J\)\nbd{}invariant and supported
  in the region \(\abs{x}\notin[1-\epsilon,1+\epsilon]\), the first
  assertion follows from the second one.  It is clear that
  \((\ln\abs{x})^{-1-\epsilon}\abs{x}^{-1/2}
  W_v|_{\abs{x}\ge1+\epsilon}\) comes from a positive measure of finite
  total volume for all places~\(v\).  The total volume of the measure
  that corresponds to \((\ln\abs{x})^{-1-\epsilon} \abs{x}^{-1/2}
  E_v|_{\abs{x}\ge1+\epsilon}\) is comparable to \(q_v^{-1}\cdot (\ln
  q_v)^{-1-\epsilon}\) for \(\abs{q_v}\to\infty\).  It is well-known that
  the sum of these terms converges.  Hence the sum also defines a
  measure of finite total volume.  (If we take~\(W_v\) instead of~\(E_v\),
  only \((\ln\abs{x})^{-1-\epsilon}\abs{x}^{-1}
  W|_{\abs{x}\ge1+\epsilon}\) has finite total volume.)
\end{proof}

The Fourier-Laplace transform identifies \(\Sch(\ICL_S)_{1/2}\) with the
space of Schwartz functions on the \emph{critical line}
\(\abs{x}^{1/2}\widehat{\ICL_S}\).  Lemma~\ref{lem:estimate_Wtwo} yields
that \(E_S(f)\) is given by pairing~\(\hat{f}\) with some tempered
distribution on the critical line.  We cannot hope for much better
error terms.  Since \(\Pi_S\) and \(\Pi_S\circ J\) are supported in the
disjoint regions \(\abs{x}>1\) and \(\abs{x}<1\), the asymptotics of
\(W\cdot\ln{}\abs{x}\) and~\(\Pi_S\) for \(\abs{x}\to\infty\) are
essentially equivalent.

\subsection{Proof of the Prime Number Theorem}
\label{sec:proof_PNT}

As usual, the first step is to show that \(\spec\pi_-\) is contained in
the interior of the critical strip.  The assertion that Riemann's
\(\zeta\)\nbd{}function does not vanish on the line \(\RE s=1\) is due to
Hadamard and de la Vallée Poussin.  A more abstract formulation of their
argument, due to Pierre Deligne (\cite{Deligne:Weil2}), also works for
rather general classes of \(L\)\nbd{}functions.  Our proof below is
motivated by the argument in~\cite{Deligne:Weil2}.

\begin{theorem}  \label{the:boundary_spectrum}
  Any \(\chi\in\spec\pi_-\) satisfies \(0<\RE\chi<1\).
\end{theorem}

\begin{proof}
  Due to the symmetry~\(J\), it suffices to rule out the existence of
  \(\chi\in\spec\pi_-\) with \(\RE\chi=0\), that is,
  \(\chi\in\widehat{\ICL_\GF}\).
  
  Let \(h\in\CCINF(\R_+)\) with \(h(x)=0\) for \(x\ge 1\) and \(h(x)=1\) for
  \(x\ll 1\) and consider the distributions \(h\abs{x}^\epsilon W\) for
  \(\epsilon>0\).  Since \(h(1)=0\), there is no need to take principal
  values for \(hW\) and the discriminant term in~\(W\) is also canceled.
  Thus \(h\abs{x}^\epsilon W\) is a positive measure on~\(\ICL_\GF\).  Its
  total volume is given formally by \(h\abs{x}^\epsilon
  W(1)=W(h\abs{x}^\epsilon)\).  We claim that this is finite.  We do a
  more general computation and consider also
  \(\chi\in\widehat{\ICL_\GF}\).  Since \(h\abs{x}^\epsilon\chi\in
  \Sch(\ICL_\GF)_{\ooival{-\epsilon,1+\epsilon}}\),
  Proposition~\ref{pro:difference_in_strip} yields
  \(\IN\pi(h\abs{x}^\epsilon\chi)\in \ell^{[0]}(\Hilm^0)\).  By
  Theorem~\ref{the:global_trace_formula}, \(W(h\abs{x}^\epsilon\chi)\)
  is the supertrace of this operator.  Hence it is finite and can be
  computed spectrally as
  \begin{displaymath}
    W(h\abs{x}^\epsilon\chi) =
    \sum_{\omega\in\Irrep(\ICL_\GF)} \mult(\omega, \pi) \cdot
    \hat{h}(\abs{x}^\epsilon\chi\omega).
  \end{displaymath}
  The positivity of \(h\abs{x}^\epsilon W\) gives constraints on the
  spectrum of~\(\pi\).  However, this is difficult to use because all
  quasi-characters interact.  Therefore, we consider the limit
  \(\epsilon h\abs{x}^\epsilon W\) for \(\epsilon\searrow0\), which
  suppresses the interior of the critical strip.

  We need the Bohr compactification to define this limit.
  Equip~\(\widehat{\ICL_\GF}\) with the discrete topology and let~\(B\) be
  its Pontrjagin dual.  This is a (non-metrizable) compact group
  called the \emph{Bohr compactification} of~\(\ICL_\GF\).  The space
  \(C(B)\) of continuous functions on~\(B\) is the closed linear span of
  the characters \(\chi\in\widehat{\ICL_\GF}\) in \(C_b(\ICL_\GF)\).  By
  restriction, \(\phi_\epsilon\defeq \epsilon h\abs{x}^\epsilon W\) is a
  positive linear functional on \(C(B)\) for all \(\epsilon>0\).  Its
  Fourier transform is \(\hat\phi_\epsilon(\chi)=\epsilon
  W(h\abs{x}^\epsilon\chi)\).  We claim that
  \[
  \lim_{\epsilon\searrow0} \hat\phi_\epsilon(\chi) = \mult(\chi,\pi)
  \qquad\text{for all \(\chi\in\widehat{\ICL_\GF}\).}
  \]
  Let \(h'\defeq dh/d\ln{}\abs{x}\) be the derivative of~\(h\) in the
  direction~\(\abs{x}\).  This function has compact support because~\(h\)
  is constant for \(\abs{x}\to0,\infty\).  Therefore, the set of
  functions \(\abs{x}^\epsilon\chi h'\) for \(\epsilon\in[0,1]\) and
  fixed~\(\chi\) is bounded in \(\CCINF(\ICL_\GF)\) and the sums
  \begin{equation}  \label{eq:W_restricted_sum}
    W(\abs{x}^\epsilon\chi h')
    = \sum_{\omega\in\Irrep(\ICL_\GF)} \mult(\omega, \pi) \cdot
    \hat{h}'(\abs{x}^\epsilon\chi\omega)
  \end{equation}
  are uniformly absolutely convergent for \(\epsilon\in[0,1]\) by
  Theorem~\ref{the:summable_spectrum}.  The functions \(\hat{h}\)
  and~\(\hat{h}'\) vanish outside the principal component
  \(\Irrep(\ICL_\GF/\ICL\one_\GF)\) and are related by
  \(\widehat{h'}(\abs{x}^s)= -s\cdot \hat{h}(\abs{x}^s)\).
  Thus~\(\hat{h}\) is a meromorphic function whose only pole is a single
  pole at the trivial representation.  Let \(U_\alpha \defeq
  \{\abs{x}^s \mid \abs{s}<\alpha\}\).  The uniform absolute
  convergence of the sum for~\(h'\) yields
  \[
  \sum_{
  \{\omega\in\Irrep(\ICL_\GF)\mid
  \omega\chi\abs{x}^\epsilon\notin U_\alpha\}}
  \bigl|
    \mult(\omega, \pi) \hat{h}(\abs{x}^\epsilon\chi\omega)
  \bigr|
  \le C/\alpha
  \]
  with some constant~\(C\) that only depends on \(\chi\) and~\(\alpha\) but
  not on~\(\epsilon\).  Hence
  \[
  \lim_{\epsilon\searrow0} \widehat{\phi_\epsilon}(\chi)
  =
  \lim_{\epsilon\searrow0}
  \sum_{\omega\in \chi^{-1}\abs{x}^{-\epsilon} U_\alpha}
  \mult(\omega,\pi) \cdot
  \epsilon\cdot \hat{h}(\abs{x}^\epsilon\chi\omega)
  \]
  for all \(\alpha>0\).  The spectrum of any summable representation
  is discrete.  Therefore, for sufficiently small \(\alpha\)
  and~\(\epsilon\) the only quasi-character in
  \(\chi^{-1}\abs{x}^{-\epsilon} U_\alpha\) that can possibly be
  contained in~\(\pi\) is~\(\chi^{-1}\).  This yields
  \begin{displaymath}
    \lim_{\epsilon\searrow0} \widehat{\phi_\epsilon}(\chi)
    =
    \lim_{\epsilon\searrow0}
    \mult(\chi^{-1},\pi) \cdot \epsilon
    \hat{h}(\abs{x}^\epsilon).
  \end{displaymath}
  The claim follows because \(\lim_{\epsilon\searrow0} \epsilon
  \hat{h}(\abs{x}^\epsilon) = -\widehat{h'}(1) = -\int_{\ICL_\GF}
  h'(x) \,d\inv x = h(0)-h(\infty) = 1\) and \(\mult(\chi^{-1},\pi) =
  \mult(\conj\chi,\pi) = \mult(\chi,\pi)\).

  Consequently, the positive linear functionals~\(\phi_\epsilon\) are
  uniformly bounded.  Since they converge when paired with the
  generators \(\chi\in\widehat{\ICL_\GF}\) of \(C(B)\), they converge in
  the weak-\(\ast\)\brd{}topology on \(C(B)\) towards a positive linear
  functional~\(\phi_0\), whose Fourier transform is \(\mult(\chi,\pi)\).
  By the Gelfand-Naimark-Segal (GNS) construction \(\phi_0\) gives rise
  to a \(\ast\)\nbd{}representation~\(\rho\) of \(C(B)\) with a cyclic
  vector~\(\xi\) such that \(\braket{\rho(f)\xi}{\xi}=\phi_0(f)\) for all
  \(f\in C(B)\).  Plugging in characters
  \(\chi\in\widehat{\ICL_\GF}\subseteq C(B)\), we get a unitary
  representation~\(\rho\) of~\(\widehat{\ICL_\GF}\) with
  \(\braket{\xi}{\rho_\chi\xi}=\mult(\chi,\pi)\) for all
  \(\chi\in\widehat{\ICL_\GF}\).
  
  Since \(\mult(1,\pi)=1\), we have \(\norm{\xi}=1\) and hence
  \(\abs{\mult(\chi,\pi)} = \abs{\braket{\xi}{\rho_\chi\xi}}\le 1\) for
  all \(\chi\in\ICL_\GF\).  However, the only possible values are \(0\) and
  \(\pm1\).  The value~\(+1\) only occurs for \(\chi=1\) because this is the
  only element of \(\widehat{\ICL_\GF}\cap\spec\pi_+\).  If
  \(\mult(\chi,\pi)=\mult(\chi',\pi)=-1\), then
  \(\rho_\chi\xi=-\xi=\rho_{\chi'}\xi\) because~\(-\xi\) is the only unit
  vector with \(\braket{\eta}{\xi}=-1\).  This yields
  \(\mult(\chi\chi',\pi)=+1\) and hence \(\chi\chi'=1\).  Therefore, there
  is at most one \(\chi\in\widehat{\ICL_\GF}\cap\spec\pi_-\), and
  this~\(\chi\) satisfies \(\chi^2=1\).
  
  Finally, it remains to show that no~\(\chi\) with \(\chi^2=1\) can
  belong to \(\spec\pi_-\).  I thank Peter Sarnak for pointing out to me
  how this can be proved using the Landau Lemma.  The following
  argument is a variant of this known trick.  We proceed slightly
  differently because we want to use as little function theory as
  possible.  Choose any \(\chi\in\Irrep(\ICL_\GF)\) with \(\chi^2=1\).
  Thus \(\chi(x)=\pm1\) and \(1+\chi(x)\ge0\) for all \(x\in\ICL_\GF\).  We
  have to prove that \(\chi\notin\spec\pi_-\).
  
  Let~\(h\) be as above and consider the function \(F(s)\defeq
  W(h\cdot(1+\chi)\cdot\abs{x}^s)\).  The same reasoning as above shows
  that this is well-defined for \(\RE s>0\) and that there is a positive
  measure~\(\mu\) on~\(\ICL_\GF\) with support in the region \(\abs{x}<1\)
  such that
  \begin{equation}  \label{eq:F_measure}
    F(s)= \int \abs{x}^s \,d\mu(x)
  \end{equation}
  for all \(s\in\C\) with \(\RE s>0\).
  We can also describe~\(F\) spectrally as
  \begin{multline}  \label{eq:F_spectral}
    F(s) = W(h\abs{x}^s) + W(h\chi\abs{x}^s)
    \\ = \sum_{\omega\in\Irrep(\ICL_\GF)} \bigl(\mult(\omega,\pi) +
    \mult(\chi\omega,\pi)\bigr) \hat{h}(\omega\abs{x}^s)
    =  \sum_{t\in\C} \rho(t) \hat{h}(\abs{x}^{s-t})
  \end{multline}
  with \(\rho(t)\defeq \bigl(\mult(\abs{x}^{-t},\pi) +
  \mult(\chi\abs{x}^{-t},\pi)\bigr)\).  The series converges absolutely
  for \(\RE s>0\).  Our next goal is to extend~\(F\) meromorphically
  to~\(\C\) using~\eqref{eq:F_spectral}.
  
  We have observed above that the function \(s\mapsto
  \hat{h}(\abs{x}^s)\) is a meromorphic function on~\(\C\) whose only
  pole is a simple pole of residue~\(1\) at \(s=0\).  Moreover, the same
  reasoning as for~\eqref{eq:W_restricted_sum} shows that the sum over
  \(\abs{s-t}\ge\alpha\) in~\eqref{eq:F_spectral} converges absolutely
  for any \(\alpha>0\), \(s\in\C\), with uniform absolute convergence
  for~\(s\) in a compact subset of~\(\C\).  Let \(s\in\C\).  Since the
  support of~\(\rho\) is discrete, there is \(\alpha>0\) such that
  \(\rho(t)=0\) for all \(t\in\C\setminus\{s\}\) with \(\abs{s-t}<2\alpha\).
  Hence
  \[
  F(s')\defeq
  \rho(s)\hat{h}(\abs{x}^{s'-s}) +
  \sum_{t\in\C, \abs{s'-t}\ge\alpha} \rho(t) \hat{h}(\abs{x}^{s'-t})
  \]
  defines a meromorphic function in the \(\alpha\)\nbd{}neighborhood
  of~\(s\).  Its only possible pole is a simple pole at~\(s\) with residue
  \(\rho(s)\).  There is no pole if \(\rho(s)=0\), of course.  These local
  definitions are compatible and extend~\(F\) to a meromorphic function
  on~\(\C\).

  Let \(s_0\in[-\infty,0]\) be the maximum of all \(s\in\R\) where~\(F\) has
  a pole, or~\(-\infty\) if~\(F\) has no poles on the real axis.  We have
  \(F(s)=\int \abs{x}^s\,d\mu(x)\) for \(s>0\) by~\eqref{eq:F_measure}.
  We claim that this continues to hold for \(s>s_0\); thus \(\abs{x}^s\)
  is in \(L^1(\mu)\) for all \(s\in\R\).  Before we prove this plausible
  claim, let us see how it finishes the proof of the theorem.  If
  \(s_0=-\infty\), then \(F(s)=\int \abs{x}^s\,d\mu(x)\) would be an
  entire function, hence \(\rho(t)=0\) for all \(t\in\C\) and \(F=0\)
  by~\eqref{eq:F_spectral}.  Since this is evidently false, \(s_0\in\R\)
  is a pole of~\(F\).  As any pole of~\(F\), it is simple and has residue
  \(\rho(s_0)\), which is non-zero because~\(s_0\) is a pole.  Thus
  \((s-s_0)F(s)\) converges to \(\rho(s_0)\) for \(s\to s_0\).  Since
  \(F(s)=\int \abs{x}^s\,d\mu(x)\ge0\) for \(s>s_0\), we get
  \(\rho(s_0)>0\).  Since \(\spec\pi_+=\{\abs{x}^0,\abs{x}^1\}\), this can
  only happen at \(0\) or~\(-1\).  Since
  \[
  \rho(0)=\rho(-1)
  = \mult(\abs{x}^0,\pi) - \mult(\chi,\pi)
  = 1 - \mult(\chi,\pi)
  \]
  we get \(\mult(\chi,\pi)=0\), that is, \(\chi\notin\spec\pi_-\), and we
  are done.

  It remains to prove \(F(s)=\int \abs{x}^s\,d\mu(x)\) for \(s>s_0\).
  Let~\(s_1\) be the infimum of the \(s\in\R\) for which \(\int
  \abs{x}^s\,d\mu(x)<\infty\).  The desired equality holds for \(s>s_1\)
  by analytic continuation.  We must show that \(s_1\le s_0\).  Assume
  the contrary.  Then~\(F\) is holomorphic in a
  \(4\epsilon\)\nbd{}neighborhood of~\(s_1\) for some \(\epsilon>0\).  Hence
  the spectral radius of the Taylor expansion of~\(F\) at \(s\defeq
  s_1+\epsilon\) is at least \(3\epsilon\).  We have
  \[
  \frac{d^n F}{dt^n}(t)
  = \int \frac{d^n \abs{x}^t}{dt^n} \,d\mu(x)
  = \int (\ln \abs{x})^n \abs{x}^t \,d\mu(x)
  \]
  for all \(t>s_1\) because~\eqref{eq:F_measure} holds in that region.
  Hence we get the estimate
  \[
  \frac{1}{n!}
  \left\lvert\int (\ln \abs{x})^n \abs{x}^s \,d\mu(x) \right\rvert
  \le C(2\epsilon)^{-n}
  \]
  for some constant~\(C\).  Since~\(\mu\) is supported in the region
  \(\abs{x}<1\), the integrand has constant sign on the support
  of~\(\mu\).  Hence we can move the absolute values into the integral
  in the above estimate.  Therefore, the power series in \(L^1(\mu)\)
  \[
  \abs{x}^t
  = \sum_{n=0}^\infty \frac{1}{n!} \abs{x}^s (\ln\abs{x})^n (t-s)^n
  \]
  converges absolutely for \(t> s-2\epsilon\), so that \(\abs{x}^t\in
  L^1(\mu)\).  Hence we arrive at the contradiction \(s-2\epsilon\ge
  s_1\), which finishes the proof.
\end{proof}

To cancel the factor \(\ln{}\abs{x}\) that occurs when comparing~\(W\)
and~\(\Pi_S\), let \(l_\xi(x)\) be \(1/\ln x\) for \(3/2\le x\le \xi\) and~\(0\)
otherwise.  We choose an isomorphism
\(\ICL_S\cong\ICL\one_S\times\R\inv_+\) and use this to view~\(l_\xi\) as
a distribution on~\(\ICL_S\) supported in \(\{1\}\times [3/2,\xi]\).  As
usual, we write \(a(\xi)\approx b(\xi)\) if \(\lim_{\xi\to\infty}
a(\xi)/b(\xi)=1\).

\begin{lemma}  \label{lem:PNT_one}
  For any \(f\in\CCINF(\ICL_S)\), we have
  \[
  W(f\ast l_\xi) \approx
  \frac{\xi}{\ln\xi} \hat{f}(\abs{x}^1) =
  \frac{\xi}{\ln\xi} \int_{\ICL_S} f(x) \abs{x}\,d\inv x.
  \]
\end{lemma}

\begin{proof}
  The global trace formula of Theorem~\ref{the:global_trace_formula}
  yields
  \begin{equation}  \label{eq:Wflxi}
    \lim_{\xi\to\infty} W(f\ast l_\xi) \frac{\ln\xi}{\xi}
    =
    \lim_{\xi\to\infty} \sum_{\chi\in\widehat{\ICL_S},\ s\in\C}
    \mult(\chi\abs{x}^s,\pi) \hat{f}(\chi\abs{x}^s)
    \widehat{l_\xi}(\abs{x}^s)
    \frac{\ln\xi}{\xi}.
  \end{equation}
  It is well-known that
  \[
  \widehat{l_\xi}(\abs{x}) =
  \int_{3/2}^\xi \frac{t\,d\inv t}{\ln t} =
  \int_{3/2}^\xi \frac{dt}{\ln t} \approx
  \frac{\xi}{\ln\xi}.
  \]
  The positivity of~\(l_\xi\) implies
  \(\norm{l_\xi}_{L^1(\R\inv_+,t^\alpha d\inv t)} =
  \widehat{l_\xi}(\abs{x}^\alpha) \le \widehat{l_\xi}(\abs{x})\) for
  all \(\alpha\le1\).  Since the \(L^1\)\nbd{}norm of~\(f\) controls the
  \(L^\infty\)\nbd{}norm of~\(\hat{f}\), we conclude that the functions
  \(\widehat{l_\xi}\cdot \ln\xi/\xi\) are uniformly bounded in the
  region \(\RE s\le1\).  Hence we can commute the summation
  over~\(\omega\) and the limit over~\(\xi\) in~\eqref{eq:Wflxi}.  We can
  also estimate
  \[
  \norm{l_\xi}_{L^1(\R\inv_+,t^\alpha\,d\inv t)} \le
  \int_{3/2}^\xi \frac{t^\alpha\,d\inv t}{\ln(3/2)} =
  \frac{\xi^\alpha-(3/2)^\alpha}{\alpha\ln(3/2)} =
  o(\xi/\ln\xi)
  \]
  for all \(\alpha<1\).  Arguing as above we get \(\lim_{\xi\to\infty}
  \widehat{l_\xi}(\abs{x}^s)\cdot \ln\xi/\xi = 0\) for all \(s\in\C\)
  with \(\RE s<1\).  Thus only \(\omega\in\spec\pi\) with \(\RE\omega=1\)
  survive the limit \(\xi\to\infty\).
  Theorem~\ref{the:boundary_spectrum} and the asymptotics for
  \(\widehat{l_\xi}(\abs{x})\) imply the assertion.
\end{proof}

Recall that \(W(f)= \Pi_S(f\cdot\ln{}\abs{x}) + \Pi_S\circ
J(f\cdot\ln{}\abs{x}) + E_S(f)\).  The same estimates as above yield
\(E_S(f\ast l_\xi)=O(\xi^{1/2}/\ln\xi)\) for \(\xi\to\infty\).  Since
\(\Pi_S\circ J((f\ast l_\xi) \cdot\ln{}\abs{x})\) remains bounded, we
get
\begin{equation}  \label{eq:asymptotics_Pi}
  \Pi_S\bigl(\ln{}\abs{x}\cdot(f\ast l_\xi)\bigr) \approx
  \int_{\ICL_S} f(x)\abs{x}\,d\inv x
  \cdot \frac{\xi}{\ln\xi}
  \qquad
  \text{for \(\xi\to\infty\)}
\end{equation}
for all \(f\in\CCINF(\ICL_S)\).  Let \(A\subseteq\ICL\one_S\) be an
open subset and let \(\pi_A(\xi)\) be the number of places
\(v\in\Cont{S}\) with \(p_v\in A\times [1,\xi]\subseteq\ICL_S\).

\begin{theorem}  \label{the:PNT}
  We have \(\pi_A(\xi) \approx \vol(A)\cdot \xi/\ln\xi\) for
  \(\xi\to\infty\).
\end{theorem}

\begin{proof}
  Decompose \(\ICL_S\cong\ICL_S\one\times\R\inv_+\) as above and
  let~\(f_0\) be the distribution \(\braket{f_0}{f_1} = \int_A
  f_1(x,1)\,d\inv x\) on~\(\ICL_S\).  We formally have
  \[
  \pi_A(\xi) =
  \braket{\Pi_S}{f_0\otimes 1_{[3/2,\xi]}} =
  \braket{\Pi_S}{\ln{}\abs{x}\cdot (f_0\ast l_\xi)}.
  \]
  If~\(f_0\) were a smooth function, the desired result would follow
  from~\eqref{eq:asymptotics_Pi}.  Therefore, we proceed by
  smoothening~\(f_0\) and checking that this does not change the
  asymptotics.  Fix \(\beta>1\).  There exist \(u_0,u_1\in\CCINF(\ICL_S)\)
  close to~\(f_0\) such that
  \[
  \beta^{-1}\vol(A) <
  \int_{\ICL_S} u_0(x)\abs{x}\,d\inv x,
  \int_{\ICL_S} u_1(x)\abs{x}\,d\inv x <
  \beta\vol(A)
  \]
  and such that \(\ln{}\abs{x} (u_0\ast l_\xi)\) is at least~\(1\) on
  \(A\times[\beta,\beta^{-1}\xi]\) and \(\ln{}\abs{x} (u_1\ast l_\xi)\) is
  supported in \(A\times[\beta^{-1},\beta\xi]\) and bounded by~\(1\)
  there.  Equation~\eqref{eq:asymptotics_Pi} applied to the smooth
  functions \(u_0\) and~\(u_1\) implies
  \[
  \pi_A(\beta^{-1}\xi) \le \beta \vol(A) \xi/\ln\xi,
  \qquad
  \pi_A(\beta\xi) \ge \beta^{-1} \vol(A) \xi/\ln\xi
  \]
  for sufficiently large~\(\xi\).  This implies the assertion for
  \(\beta\searrow1\).
\end{proof}

\end{document}